\documentclass[letterpaper,10pt]{article}

\usepackage{amsfonts,amsmath,verbatim,graphicx, subfigure,longtable}
\usepackage{amssymb,amsthm,amscd}
\usepackage{hyperref}
\usepackage{pstricks}
\usepackage{epstopdf}
\usepackage{color}
\usepackage{transparent}
\usepackage{tikz}
\usetikzlibrary{fit, positioning}
\newtheorem{thm}{Theorem}[section]
\newtheorem{cor}[thm]{Corollary}
\newtheorem{prop}[thm]{Proposition}
\newtheorem{lem}[thm]{Lemma}
\newtheorem{definition}[thm]{Definition}
\newtheorem{assumption}[thm]{Assumption}

\newtheorem{example}[thm]{Example}
\newtheorem{remark}[thm]{Remark}

\newenvironment{pf}{\noindent\emph{Proof \,}}{\mbox{}\qed}

\newcommand{\toL}{\,{\buildrel \mathcal{L} \over \longrightarrow}\,}
\newcommand{\equalL}{\,{\buildrel \mathcal{L} \over =}\,}

\numberwithin{equation}{section}
\def\R{{\mathbb R}}     
\def\P{{\mathbb P}}     
\def\E{{\mathbb E}}     
\def\D{{\mathcal E}}    
\def\F{{\mathcal F}}    
\def\H{{\mathcal H}}    %

\def\X{{\mathfrak X}}   
\def\Y{{\mathcal Y}}   

\def\ZZ{{\mathcal Z}}
\def\A{{\mathcal A}}   

\def\<{{\langle}}
\def\>{{\rangle}}
\def\1{{\bf 1}}         

\renewcommand{\bar}{\overline}

\begin{document}
\allowdisplaybreaks

\title{\Large \bf
Fluctuation limit for interacting diffusions \\ with partial annihilations through membranes\thanks{Research partially supported by NSF Grants DMS-1206276 and DMR-1035196.}
 }

\author{{\bf Zhen-Qing Chen} \quad and \quad {\bf Wai-Tong (Louis) Fan}}
\date{\today}
\maketitle

\begin{abstract}
We study  fluctuations of the empirical processes of  a non-equilibrium interacting particle system  consisting of two species over a domain that is 
recently introduced in \cite{zqCwtF14b} and establish its functional central limit theorem. This fluctuation limit is a distribution-valued Gaussian Markov process which can be represented as a mild solution of a stochastic partial differential equation. 
The drift of our fluctuation limit involves  
a new partial differential equation with nonlinear coupled term on the interface that characterized the hydrodynamic limit
of the system. The covariance structure of the Gaussian part consists two parts, one involving the spatial motion of the particles
inside the domain and other involving a boundary integral term that captures the  boundary interactions between two species.  
The key is to show that the Boltzman-Gibbs principle holds for our non-equilibrium system.  Our proof relies on generalizing the usual correlation functions to the join correlations at two different times.



\bigskip
\end{abstract}

\tableofcontents

\section{Introduction}

A new class  
of non-equilibrium  particle systems of two species that interact with each other along a  hypersurface
is recently introduced in \cite{zqCwtF13a} and \cite{zqCwtF14b}. The primary goal is to understand the connection between the microscopic transports of positive and negative charges in solar cells and the electric current generated. However, these models are flexible and general enough to provide insight to a variety of natural phenomena, such as the population dynamics of two segregated species under competition.

\begin{figure}[h]
	\begin{center}
    \vspace{-1em}
	\includegraphics[scale=.3]{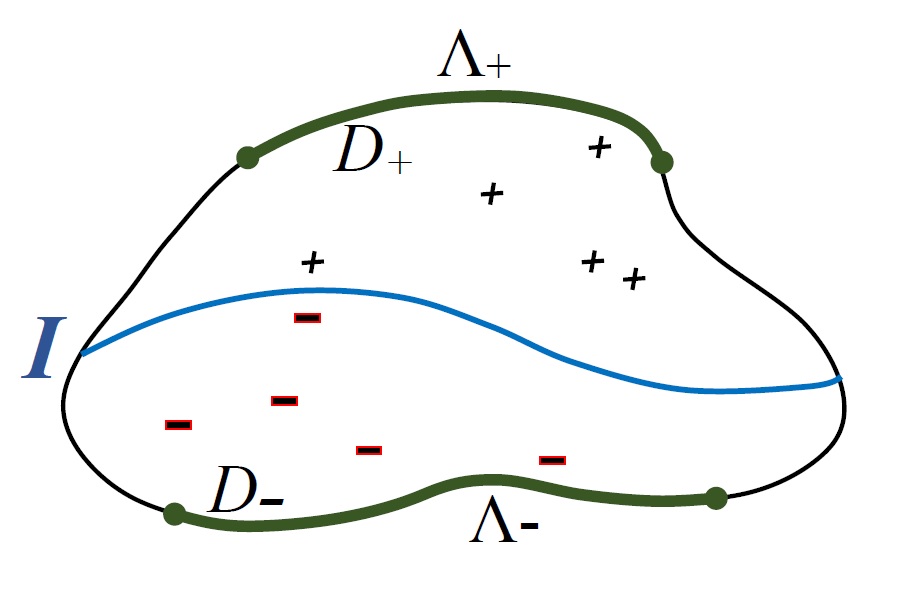}\label{fig:Model3Clean}
    \vspace{-1em}
    \caption{Annihilating diffusion model, where ${\color{blue}I}$ is the interface and ${\Lambda_{\pm}}$ are harvest sites}
    \vspace{-1em}
	\end{center}
	\end{figure}

Here is an informal description of the model introduced in \cite{zqCwtF14b}. A solar cell is modeled by a domain in $\R^d$ that is divided into two adjacent sub-domains $D_+$ and $D_-$ by an interface $I$, a $(d-1)$-dimensional Lipschitz hypersurface. Domains $D_+$ and $D_-$ represent the hybrid medium that confine the positive and the negative charges, respectively. An example to keep in mind is when $D_+=(0,1)^d$ and $D_-=(0,1)^{d-1}\times (0,-1)$ are two adjacent unit cubes. The interface is then $I=(0,1)^{d-1}\times\{0\}$. The particle system is indexed by $N$, the initial number of positive and negative charges in each of $D_+$ and $D_-$. At microscopic level, the motion of positive and negative charges are modeled by independent reflected diffusions (such as reflected Brownian motions) in $D_+$ and $D_-$, respectively. Besides, there is a harvest region $\Lambda_{\pm}\subset \partial D_{\pm}\setminus I$ that  absorbs (harvests) $\pm$ charges, respectively, whenever it is being visited. Furthermore, when two particles of different types are within a small distance $\delta_N$, they disappear at a certain rate\footnote{We say an event happens at rate $r$ if the time of occurrence is an exponential random variable of parameter $r$. In particular, the probability of occurrence in a short amount of time $t$ is $rt+o(t)$, where $o(t)/t \to 0$ as $t\to 0$.} $r_N$. This annihilation models the trapping, recombination and separation phenomena of the charges. We shall refer to the system just described as the \emph{annihilating diffusion model}. See Figure \ref{fig:Model3Clean} for an illustration.

Even though the boundary is fixed and there is no creation of particles, the interactions do affect the correlations among the particles: Whether or not a positive particle disappears at a given time affects the empirical distribution of the negative particles, which in term affects that of the positive particles. This challenge is reflected by the non-linearity of the macroscopic limit.
This challenge arises again in the study of its fluctuation limit, and it is further reflected by the boundary integral term
in the covariance of the Gaussian process; see \eqref{E:QuadVarMN_CLT} of Theorem \ref{T:Convergence_AnnihilatingSystem}.

In \cite{zqCwtF14b}, we established a \emph{functional law of large numbers} for the time trajectory of the particle densities. This is a first step in connecting the microscopic mechanism of the system with the macroscopic behaviors that emerge. More precisely, let $(\X^{N,+}_t,\,\X^{N,-}_t)$ be the pair of empirical measure of positive and negative charges at time $t$. We showed that, under a suitable scaling and appropriate conditions on the initial configurations, the random pairs of measures $(\X^{N,+}_t,\,\X^{N,-}_t)$ converge in distribution, as $N\to\infty$, to a pair of deterministic measures which are absolutely continuous with respect to the Lebesgue measure. Furthermore, the densities with respect to Lebesque measure satisfy a system of partial differential equations (PDEs) that has coupled nonlinear boundary conditions on the interface; see Theorem \ref{T:Conjecture_delta_N_CLT}. It is this nonlinear coupling effect near the interface that distinguishes this model from previously studied ones. The suitable scaling is  of order $(N\delta_N^{d+1})^{-1}$ and is rigorously formulated via the annihilation potential function in Assumption \ref{A:The annihilation potentialCLT}.

In the current work, we look at a finer scale of the annihilating diffusion model and establish the \emph{functional central limit theorem} in Theorem \ref{T:Convergence_AnnihilatingSystem}.  To focus on the fluctuation effect caused by the interaction on the interface $I$, we assume the harvest sites are empty in this paper. The fluctuations of the empirical measures from their mean (the coupled PDEs) is quantified by
\begin{equation}\label{Def:FluctuationField}
    \Y^{N,\pm}_t(\phi_{\pm}) :=  \sqrt{N}\,(\<\X^{N,\pm}_t,\phi_{\pm}\>-\E\<\X^{N,\pm}_t,\phi_{\pm}\>),
\end{equation}
where $\<\X^{N,\pm}_t,\phi_{\pm}\>$ is the integral of an observable (or test function) $\phi_{\pm}\in L^2(D_{\pm})$ with respect to the measure $\X^{N,\pm}_t$. Intuitively, if $\phi_{+}=\1_{K}$ is an indicator function of a subset $K\subset \bar{D}_{+}$, then $\<\X^{N,+}_t,\phi_{+}\>$ is the mass of particles in $K$ (which is the number of particles in $K$ divided by $N$). In this case, $\Y^{N,+}_t(\phi_{+})$ is the fluctuation of the mass of particles in $K$ at time $t$. Our main result in this paper, Theorem \ref{T:Convergence_AnnihilatingSystem}, asserts that the fluctuation limit (as $N\to\infty$) is a continuous Gaussian Markov process whose covariance structure is explicitly characterized. Roughly speaking, the limit  solves a stochastic partial differential equation (SPDE) which is a nonlinear version of the Langevin equation.

As a preliminary step to understand the fluctuation for the annihilating diffusion model, we consider in \cite{zqCwtF14c} a simpler single species model. In that paper, the particles move as i.i.d. reflected Brownian motions in a bounded domain $D\subset \R^d$ and are killed by a singular time-dependent potential which concentrates on the boundary $\partial D$. This is motivated by observation that we can view the positive charges as reflected diffusions in $D_+$ subject to killings by a time-dependent random potential. The techniques developed in \cite{zqCwtF14c} provides us with a functional analytic setting for our fluctuation processes $\Y^{N,\pm}$ and allow us to overcome some (but not all) challenges for the study of the fluctuation of the annihilating diffusion model. For the latter, we need two new ingredients, namely the \emph{asymptotic expansion of the correlation functions} and the \emph{Boltzman-Gibbs principle}. More precisely, by generalizing the approach of P. Dittrich \cite{pD88b}, we show that the correlation functions have the decomposition
        \begin{equation*}
            F^{N,(n,m)}_t(\vec{x},\vec{y})= \prod_{i=1}^nu_{+}(t,x_i)\,\prod_{j=1}^mu_{-}(t,y_j) +\frac{B_t^{N,(n,m)}(\vec{x},\vec{y})}{N}+ \frac{o(N)}{N},
        \end{equation*}
where $(u_+,\,u_-)$ is the hydrodynamic limit of the interacting diffusion system, $B_t^{N,(n,m)}$ is an explicit function and $o(N)$ is a term converging to zero as $N$ tends to infinity. See Theorem \ref{T:Asymptotic_nm_correlation_t} for the precise statement. This result implies propagation of chaos and allows explicit calculations of the covariance of the fluctuation process. The proof of Theorem \ref{T:Asymptotic_nm_correlation_t} is based on a comparison of the BBGKY \footnote{BBGKY stands for N. N. Bogoliubov, Max Born, H. S. Green, J. G. Kirkwood, and J. Yvon, who derived this type of  hierarchy of equations in the 1930s and 1940s in a series of papers.} hierarchy satisfied by the correlation functions with two other approximating hierarchies. On the other hand, the Boltzman-Gibbs principle, first formulated mathematically and proven for some zero range processes in equilibrium in \cite{BR84}, says that the fluctuation fields of non-conserved quantities change on a time scale much faster than the conserved ones, hence in a time integral only the component along those fields of conserved quantities survives. Although this principle is proved to hold for a few 
non-equilibrium situations (see \cite{BDP92} and the references therein), it is not known whether it holds in general. The validity of the principle for our annihilating diffusion model is far from obvious, since there is no conserved quantity.
An intuitive explanation for the validity here is as follow: Assumption \ref{A:ShrinkingRateCLT}  guarantees that the interaction near $I$ changes the occupation number of the particles at a slower rate with respect to diffusion (which conserves the particle number). In other words, the particle number is approximately conserved on the time scale that is relevant for the principle. Hence we are not far away from equilibrium fluctuation.

One of the earliest rigorous results on fluctuation limit was proven by It\^{o} \cite{kI80, kI83}, who considered a system of independent and identically distributed (i.i.d.) Brownian motions in $\R^d$ and showed that the limit is a $\mathcal{S}'$-valued Gaussian process solving a Langevin equation, where $\mathcal{S}'$ is the Schwartz space of tempered distributions.
Fluctuations for interacting diffusions in $\R^d$ are studied by various authors; see \cite{kO87, RR10} for examples of Gaussian fluctuations and \cite{XK04} for an example of non-Gaussian fluctuations. Sznitman \cite{asS84} studied the fluctuations of a conservative system of diffusions with normal reflected boundary conditions on smooth domains.

It is well known that the correlation method works well for certain stochastic particle systems modeling the reaction-diffusion equation
\begin{equation}\label{E:RDeqt}
\partial_t u(t,x)=\frac{1}{2}\Delta u(t,x) + R(u(t,x)),
\end{equation}
where $R(u)$ is the reaction term, such as a polynomial in $u$. See  \cite{pD88b, pK86,  pK88} for the continuous setting in which particles are diffusions  on the cube $[0,1]^d$ with linear or quadratic reaction terms.  See also \cite{DmFL86, MP91, BDP92, FNV11, DmPTV12, FGN13} for the discrete setting in which particles move on a one-dimensional lattice and correlations functions are called $\nu$-functions. The key message of the fluctuation results for  the reaction-diffusion equation is that, roughly speaking, the fluctuation limit $\Y$ solves the following {\it stochastic partial differential equation} in a distributional Hilbert space:
\begin{equation*}
d\Y_t = \Big(\frac{1}{2}\Delta\Y_t+R'(u(t))\Y_t\Big)\,dt + dM_t,
\end{equation*}
where $u(t,x)$ solves equation \eqref{E:RDeqt}, $R'(u)$ is the derivative of $R(u)$ (e.g. $-2u$ when $R(u)=-u^2$) and is viewed as a multiplicative operator, $M$ is a Gaussian martingale with independent increment and covariance structure $$\E[(M_t(\phi))^2]=\int_0^t\<|\nabla\phi|^2,\,u(s)\>+\<\phi^2,\,|R(u(s))|\>\,ds.$$ Here $\<\cdot\,,\,\cdot\>$ is the $L^2$ inner product in the spatial variable and $|R(u)|$ is the polynomial obtained by putting an absolute sign to each coefficient in $R(u)$. So the  transportation component (or drift) $\frac{1}{2}\Delta\Y_t+R'(u(t))\Y_t$ of the fluctuation limit involves the derivative of the reaction term.

For our annihilating diffusion model, it is not clear which function should one ``differentiate'' with respect to for the nonlinear term in the hydrodynamic limit,  since it involves two functions. It turns out that the transportation component of our fluctuation limit $\ZZ$ in Theorem \ref{T:Convergence_AnnihilatingSystem} is described by (\ref{E:FluctuationPDE_limit}), which is another coupled PDE which can be view as a `linearization' of the hydrodynamic equation in Theorem \ref{T:Conjecture_delta_N_CLT}. 
Our fluctuation results hold for all dimensions $d\geq 1$. A distinct feature is that 
the covariance structures of our fluctuation limits have boundary integral terms that capture the boundary interactions at the fluctuation level, 
in addition to the usual energy terms that describe the spatial motion of the particles.

\section{Notations}

For the convenience of the reader, we list our notations used in this paper.

\begin{longtable}{ll}
$\mathcal{B}(E)$ & Borel measurable functions on $E$ \\
$\mathcal{B}_b(E)$& bounded Borel measurable functions on $E$ \\
$\mathcal{B}^+(E)$& non-negative Borel measurable functions on $E$\\
$C(E)$ & continuous functions on $E$\\
$C_b(E)$ & bounded continuous functions on $E$\\
$C^+(E)$ & non-negative continuous functions on $E$\\
$C_c(E)$ & continuous functions on $E$ with compact support\\
$W^{1,2}(E)$ & $\{f\in L^2(E; dx): |\nabla f| \in L^2(E; dx)\}$ Sobolev space of order $(1,2)$\\
$D([0,\infty),\,E)$ & space of $c\grave{a}dl\grave{a}g$ paths from $[0,\infty)$ to $E$ \\
&\quad equipped with the Skorokhod metric (see \cite{pB99} or \cite{EK86})\\
$\|\,\cdot\,\|$& uniform norm (unless otherwise stated)\\
$\mathcal{H}^m$ &  $m$-dimensional Hausdorff measure\\
$M_+(E)$ (or $M_{\geq 0}(E)$) & space of finite non-negative Borel measures on $E$ with the weak topology\\
$\{\F^X_t:\,t\geq 0\}$ &
augmented filtration
 induced by the process $(X_t)_{t\geq 0}$, i.e. $\F^X_t=\sigma(X_s,\,s\leq t)$\\
$\1_x$ &  indicator function at $x$ or the Dirac measure at $x$,\\
$\toL$ & convergence in law of random variables (or processes)\\
$\equalL$ & equal in law\\
$:=$ & is defined as\\
$x\vee y$ & $\max\{x,\,y\}$\\
$x\wedge y$ & $\min\{x,\,y\}$\\
$I^{\delta}$ & $\{(x,y)\in D_+\times D_-:\, |x-z|^2+|y-z|^2<\delta^2 \,\text{ for some }z\in I\}$, \\
$c_{d}$    & volume of the unit ball in $\R^{d}$\\
$\ell_{\delta}(x,y)$ & annihilating potential function in Assumption \ref{A:The annihilation potentialCLT}\\
$\<\varphi,\,\mu^+(dx) \otimes \mu^-(dy)\>$ & $\frac{1}{N^2}\sum_{i}\sum_{j}\varphi(x_i,y_j)$, when $\mu=(\frac{1}{N}\sum_{i} \1_{x_i},\,\frac{1}{N}\sum_{j}\1_{y_j})$\\
$\ZZ^N := \Y^{N,+}\oplus\Y^{N,-}$ & fluctuation process defined in (\ref{Def:ZZ^N})\\
$\mathbf{H}_{-\alpha}$ & the Hilbert space $\{\mu^+ \oplus \mu^-:\; \mu^{\pm}\in \H_{-\alpha}^{\pm} \}$ defined in (\ref{Def:mathbfH})\\
$\{\phi^{\pm}_{k}\}$ & complete orthonormal system of $\A^{\pm}:=\frac{1}{2\,\rho_{\pm}}\,\nabla\cdot(\rho_{\pm}\,\textbf{a}_{\pm}\nabla)$ \\
                        &\quad in $L^2(D_{\pm},\,\rho_{\pm})$ consisting of Neumann eigenfunctions\\
$\lambda^{\pm}_{k}$ & the eigenvalue corresponding to $\phi^{\pm}_{k}$ such that $\A^{\pm}\phi^{\pm}_k=-\lambda^{\pm}_k$ \\
$\<\phi,\,\psi\>_{\rho_{\pm}}$ & $\int_{D_{\pm}} \phi(x)\psi(x)\,\rho_{\pm}(x)dx$, the inner product of $L^2(D_{\pm},\,\rho_{\pm}(x)dx)$\\
$F^{(n,m)}_t=F^{N,(n,m)}_t$ & $(n,m)$-correlation function at time $t$ in Definition \ref{Def:CorrelationFcn}\\
$F^{(n,m),(p,q)}_{s,t}=F^{N,(n,m),(p,q)}_{s,t}$ & generalized correlation function in Definition \ref{Def:GeneralizedCorrelationFcn} \\
$E^{(n,m),(p,q)}_{u,r}$ & $F^{(n,m),(p,q)}_{u,u+r}- F^{(n,m)}_{u}\cdot F^{(p,q)}_{u+r}$ defined in (\ref{Def:Enmpq})
\end{longtable}

A constant $C$ that depends only on $D$ and $T$ will sometimes be written as $C(D,T)$. The exact value of the constant may vary from line to line. We also use the following abbreviations:

\begin{tabular}{ll}
a.s. & almost surely\\
$c\grave{a}dl\grave{a}g$
(or rcll)
& right continuous with left limits\\
LDCT & Lebesque dominated convergence theorem\\
LHS & left hand side\\
PDE & partial differential equation\\
RBM & reflected Brownian motion\\
RHS & right hand side\\
SPDE & stochastic partial differential equation\\
WLOG & without loss of generality\\
w.r.t. & with respect to\\
\end{tabular}

\medskip

\begin{definition}
A Borel subset $E$ of $\R^d$ is called $\mathcal{H}^m$\textbf{-rectifiable} if $E$ is a countable union of Lipschitz images of bounded subsets of $\R^m$ with $\mathcal{H}^m(E)<\infty$ (As usual, we ignore sets of $\mathcal{H}^m$ measure 0). Here $\mathcal{H}^m$ denotes the $m$-dimensional
Hausdorff measure.
\end{definition}

\begin{definition}
A \textbf{bounded Lipschitz domain} $D\subset\R^d$ is a bounded connected open set such that for any $\xi\in \partial{D}$, there exits $r_{\xi}>0$ such that $B(\xi,r_{\xi})\cap D$ is represented by $B(\xi,r_{\xi}) \cap  \{(y',y^{d})\in \R^d: \phi_{\xi}(y')<y^d\}$ for some coordinate system centered at $\xi$ and a Lipschitz function $\phi_{\xi}$ with Lipschitz constant $M_D$, where $M_D>0$ does not depend on $\xi$.
\end{definition}

\section{Model description and assumptions}

In this section, we recall the basic assumptions and the definition of the annihilating diffusion model introduced in
 \cite{zqCwtF14b}. To focus on the fluctuation effect caused by the interaction on the interface $I$, we assume the harvest sites are empty in this paper.

We first recall the basic notion and properties of reflected diffusions in a domain.
 Let $D\subset \R^d$ be a bounded Lipschitz domain, and $ W^{1,2}(D)=\{f\in L^2(D; dx): \nabla f \in L^2(D; dx)\}$. Consider the bilinear form on $W^{1,2}(D)$ defined by
$$\D (f,g) := \dfrac{1}{2} \int_{D}\nabla f(x)\cdot \, \textbf{a} \nabla g(x)\,\rho(x) \,dx := \frac{1}{2}\int_D\sum_{i,j=1}^da^{ij}(x)\,\frac{\partial f}{\partial x_i}(x)\,\frac{\partial g}{\partial x_i}(x)\,\rho(x)\,dx ,$$
where $\rho \in W^{1,2}(D)$ is a positive function on $D$ which is bounded away from zero and infinity, $\textbf{a}=(a^{ij})$ is a symmetric bounded uniformly elliptic $d\times d$ matrix-valued function such that $a^{ij}\in W^{1,2}(D)$ for each $i,\,j$. Since $D$ is Lipschitz boundary, $(\D, W^{1,2}(D))$ is a regular symmetric Dirichlet form on $L^2(D; \rho (x) dx)$ and hence has a unique (in law) associated $\rho$-symmetric strong Markov process $X$ (cf.  \cite{zqC93}). Denote by $\vec{n}$ the unit inward  normal vector of $D$ on $\partial D$. The Skorokhod representation of $X$ (cf.  \cite{zqC93}) tells us that  $X$ behaves like a diffusion process associated to the elliptic operator
\begin{equation*}
\mathcal{A} := \dfrac{1}{2\,\rho}\,\nabla\cdot(\rho\,\textbf{a}\nabla):= \dfrac{1}{2\,\rho}\,\sum_{i,j=1}^d \dfrac{\partial}{\partial x_i}\Big(\rho\,a^{ij}\,\dfrac{\partial}{\partial x_j}\Big)
\end{equation*}
in the interior of $D$, and is instantaneously reflected at the boundary in the inward conormal direction $\vec{\nu}:=\textbf{a}\vec{n}$.

\begin{definition}\label{Def:ReflectedDiffusion}
Let $(\textbf{a},\,\rho)$  and $X$ be as in the preceding paragraph. We call $X$ an \textbf{$(\textbf{a},\,\rho)$-reflected diffusion} or an \textbf{$(\A,\,\rho)$-reflected diffusion}.  A special but important case is when $\textbf{a}$ is the identity matrix, in which $X$ is called a reflected Brownian motion with drift $\frac{1}{2}\,\nabla (\log \rho)$. If in addition $\rho=1$, then $X$ is called a reflected Brownian motion (RBM).
\end{definition}

Assumptions \ref{A:SettingCLT} to \ref{A:ShrinkingRateLLN} below are in force throughout this paper.
\begin{assumption}\label{A:SettingCLT}(Geometric setting)
$D_{\pm}$ are given adjacent bounded Lipschitz domains in $\R^d$ such that $I:= \bar{D}_{+}\cap \bar{D}_{-}=\partial D_+\cap \partial D_-$ is a finite union of disjoint connected $\mathcal{H}^{d-1}$-rectifiable sets.
$\rho_{\pm}\in W^{1,2}(D_{\pm})\cap C(\bar{D}_{\pm})$ are given strictly positive functions, $\textbf{a}_{\pm}=(a^{ij}_{\pm})$ are symmetric, bounded, uniformly elliptic $d\times d$ matrix-valued functions such that $a^{ij}_{\pm}\in W^{1,2}(D_{\pm})$ for each $i,\,j$.
\end{assumption}
We let $X^{\pm}$  be an $(\textbf{a}_{\pm},\,\rho_{\pm})$-\textbf{reflected diffusion} in $D_{\pm}$. Under Assumption \ref{A:SettingCLT}, $X^{\pm}$ is a continuous strong Markov process with symmetrizing measure $\rho_{\pm}$ and has infinitesimal generator $\A^{\pm}:= \frac{1}{2\,\rho_{\pm}}\,\nabla\cdot(\rho_{\pm}\,\textbf{a}_{\pm}\nabla)$. An example to keep in mind is when $D_+=(0,1)^d$ and $D_-=(0,1)^{d-1}\times (0,-1)$ are two adjacent unit cubes, the functions $\rho_{\pm}=1$ are constants and $\textbf{a}_{\pm}$ are the identity matrices. The interface is then $I=(0,1)^{d-1}\times\{0\}$, and we have $\A^{\pm}=\frac{1}{2}\Delta$.

\begin{assumption}\label{A:ParameterAnnihilationCLT}(Parameter of annihilation)
Suppose $\lambda\in C_+(I)$ is a given non-negative continuous function on $I$.
Let  $\widehat{\lambda}\in C(\bar{D}_+\times \bar{D}_-)$  be an arbitrary extension of $\lambda$ in the sense that
$\widehat{\lambda}(z,z)=\lambda(z)$ for all $z\in I$.  Note that such $\widehat{\lambda}$ always exists.
\end{assumption}

\begin{figure}[h]
	\begin{center}
	\includegraphics[scale=.25]{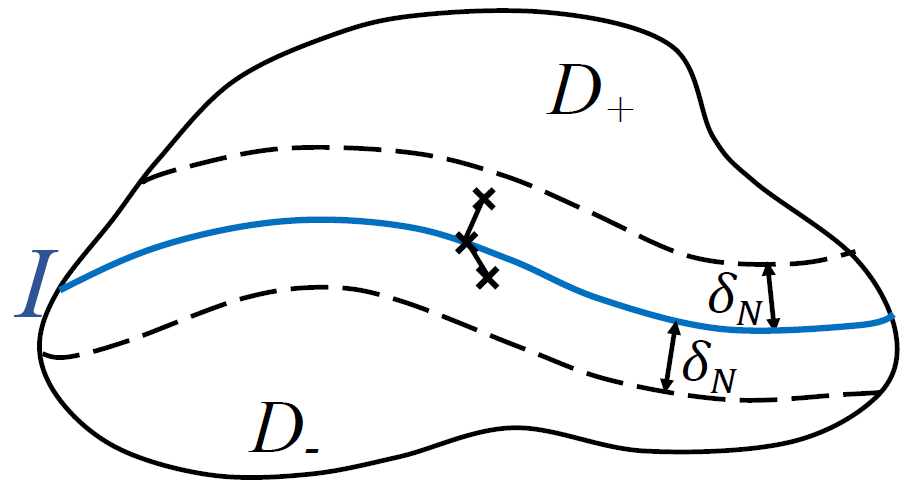}
    \vspace{-1em}
    \caption{A pair of points $(x,y)\in I^{\delta_N}$}\label{fig:AnnihilationDist}
    \vspace{-1em}
	\end{center}
	\end{figure}

\begin{assumption}\label{A:The annihilation potentialCLT}(The annihilation potential)
We choose $\{\ell_{\delta}:\,\delta>0\} \subset  C_+(\bar{D}_+\times \bar{D}_-)$ in such a way that $\ell_{\delta}(x,y)\leq \frac{\widehat{\lambda}(x,y)}{c_{d+1}\,\delta^{d+1}}\1_{I^{\delta}}(x,y)$ on $D_+\times D_-$ and
\begin{equation}
\lim_{\delta\to 0} \Big\|\ell_{\delta}-\frac{\widehat{\lambda}}{c_{d+1}\,\delta^{d+1}}\1_{I^{\delta}}
\Big\|_{L^2(D_+\times D_-)} =0,
\end{equation}
where $I^{\delta} :=  \{(x,y)\in D_+\times D_-:\, |x-z|^2+|y-z|^2<\delta^2 \,\text{ for some }z\in I\}$
 (see Figure \ref{fig:AnnihilationDist}) and $c_{d+1}$ is the volume of the unit ball in $\R^{d+1}$.
\end{assumption}

The motivation for the definition of $\ell_{\delta}$ is the fact (see \cite[Theorem 3.2.39]{hF69}) that
 \begin{equation}\label{E:Federer_MinkowSki_Content}
            \lim_{\delta\to 0}\frac{\mathcal{H}^{2d}(I^{\delta})}{c_{d+1}\,\delta^{d+1}}=\mathcal{H}^{d-1}(I).
 \end{equation}
In the proof of the main theorem, we will need a strengthened version of (\ref{E:Federer_MinkowSki_Content}) which is stated in Lemma \ref{L:MinkowskiContent_I_k_CLT}. Intuitively, if $N$ is the initial number of particles, then $\delta=\delta_N$ is the annihilation distance and $I^{\delta}$ controls the frequency of interactions. As remarked in the Introduction of \cite{zqCwtF14b}, we need to assume that the annihilation distance $\delta_N$ does not shrink too fast. This is formulated in Assumptions \ref{A:ShrinkingRateLLN} and \ref{A:ShrinkingRateCLT}.

\begin{assumption}\label{A:ShrinkingRateLLN}(Annihilation distance for functional LLN)
$\{\delta_N\}\subset (0,\infty)$ converges to 0 as $N\to\infty$
and $\liminf_{N\to\infty}N\,\delta_N^{d} \in (0,\infty]$.
\end{assumption}

Let
 $(\X^{N,+},\,\X^{N,-})$ be the normalized empirical measures for the annihilating diffusion system described in the Introduction and rigorously constructed in \cite{zqCwtF14b}. The main result of \cite{zqCwtF14c} implies
the following.

\begin{thm}\label{T:Conjecture_delta_N_CLT} \textbf{(Hydrodynamic Limit) }
Suppose Assumptions \ref{A:SettingCLT} to \ref{A:ShrinkingRateLLN} in the above hold. If
$(\X^{N,+}_0,\,\X^{N,-}_0) \toL (u^+_0(x)\rho_+(x)dx,\,u^-_0(y)\rho_-(y)dy)$ in $M_+(\bar{D}_+)\times M_+(\bar{D}_-)$ as $N\to \infty$, where $u^{\pm}_0\in C(\bar{D}_{\pm})$, then
\begin{equation*}
(\X^{N,+},\X^{N,-}) \toL  (u_+(t,x)\rho_+(x)dx,\,u_-(t,y)\rho_-(y)dy) \quad \text{in }D([0,T],\,M_+(\bar{D}_+)\times M_+(\bar{D}_-))
\end{equation*}
for any $T>0$, where $(u_+,\,u_-)$ is the probabilistic solution (see Remark \ref{Rk:ProbabilisticSol_CLT}) of the following coupled heat equations
    \begin{equation}\label{E:coupledpde:+++}
        \left\{\begin{aligned}
        \dfrac{\partial u_+}{\partial t}           &= \A^+ u_+     & &\qquad\text{on } (0,\infty)\times D_+  \\
        \dfrac{\partial u_+}{\partial \vec{\nu_+}} &=\frac{\lambda}{\rho_+}\,u_+u_-\,\1_{\{I\}}  & &\qquad\text{on }  (0,\infty)\times \partial D_+
        \end{aligned}\right.
    \end{equation}
and
    \begin{equation}\label{E:coupledpde:---}
        \left\{\begin{aligned}
        \dfrac{\partial u_-}{\partial t}           &= \A^- u_-   & &\qquad\text{on } (0,\infty)\times D_-  \\
        \dfrac{\partial u_-}{\partial \vec{\nu_{-}}} &= \frac{\lambda}{\rho_-}\,u_+u_-\,\1_{\{I\}}  & &\qquad\text{on }  (0,\infty)\times \partial D_-,
        \end{aligned}\right.
    \end{equation}
with initial value $(u^+_{0},u^-_{0})$, where $\vec{\nu}_{\pm}:=\textbf{a}_{\pm}\vec{n}_{\pm}$ is the inward conormal vector field of $\partial D_{\pm}$.
\end{thm}

\begin{remark}\label{Rk:ProbabilisticSol_CLT}\rm
The notion of probabilistic solution in Theorem \ref{T:Conjecture_delta_N_CLT} follows that in \cite{zqCwtF13a, zqCwtF14b}. Precisely, $(u_+,\,u_{-})$ is the unique element in $C([0,\infty)\times\bar{D}_+)\times C([0,\infty)\times\bar{D}_-)$ satisfying
    \begin{equation}\label{E:ProbabilisticRep_CoupledPDE_CLT}
    \left\{\begin{aligned}
        u_+(t,x) &= \E^{x} \Big[\,u^+_0(X^+_t)\,\exp{\Big(-\int^t_0(\lambda\,u_-)(t-s,X^+_s)\,dL^+_s \Big)}\,\Big]\\
        u_-(t,y) &= \E^{y}\Big[\,u^-_0(X^-_t)\,\exp{\Big(-\int^t_0 (\lambda\,u_+)(t-s,X^-_s)\,dL^-_s \Big)}\,\Big],
    \end{aligned}\right.
    \end{equation}
where $L^{\pm}$ is the boundary local time of the reflected diffusion $X^{\pm}$ on the interface $I$. The validity of the previous assertion can be verified by the same argument for Proposition 2.19 in \cite{zqCwtF13a}. In this chapter, $(u_+,\,u_{-})$ always denote the probabilistic solution of the coupled PDEs (also known as hydrodynamic limit) in Theorem \ref{T:Conjecture_delta_N_CLT}. \qed
\end{remark}

\section{Fluctuation process}

Our object of study in this paper is the \textbf{fluctuation process} defined by
\begin{equation}\label{Def:ZZ^N}
    \ZZ^N := \Y^{N,+}\oplus\Y^{N,-}=(\Y^{N,+}_t\oplus\Y^{N,-}_t)_{t\geq 0},
\end{equation}
where $\Y^{N,+}_t\oplus\Y^{N,-}_t\,(\phi_+,\phi_-)\,:=\,\Y^{N,+}_t(\phi_+)+\Y^{N,-}_t(\phi_-)$ and $\Y^{N,\pm}$ is the fluctuation field in $\bar{D}_{\pm}$
as defined in (\ref{Def:FluctuationField}).

\textbf{Functional analytic framework: } As in \cite{zqCwtF14c}, it is nontrivial to describe the state space of $\ZZ^N$ in which we have weak convergence. For this, we adopt the functional analytic setting developed in \cite{zqCwtF14c} to each of $D_+$ and $D_-$. Let $\{\phi^{\pm}_{k}\}$ be a complete orthonormal system (CONS) of $\A^{\pm}$ in $\H^{\pm}_0:= L^2(D_{\pm},\,\rho_{\pm})$ consisting of Neumann eigenfunctions, and $-\lambda^{\pm}_{k}$ the eigenvalue corresponding to $\phi^{\pm}_{k}$ (i.e. $\A^{\pm}\phi^{\pm}_{k}=- \lambda^{\pm}_k\,\phi^{\pm}_{k}$), with $0<\lambda^{\pm}_{1} < \lambda^{\pm}_2 \leq \lambda^{\pm}_3 \leq \cdots$. Moreover, for $\gamma\in\R$,
 let $\H_{\gamma}^{\pm}$ be
the separable Hilbert space with inner product $\<\,,\,\>^{\pm}_{\gamma}$ constructed as in \cite{zqCwtF14c}, which has CONS $\{h^{(\gamma),\pm}_{k}:= (1+\lambda^{\pm}_{k})^{-\gamma/2}\,\phi^{\pm}_{k}; k\geq 1\}$. Now for $\alpha\geq 0$ and $(\mu^{+},\,\mu^-)\in \H_{-\alpha}^+\times \H_{-\alpha}^-$, we define $\mu^+\oplus \mu^-$ by
$$\mu^+\oplus \mu^-\,(\phi_+,\phi_-):= \<\mu^+,\phi_+\>_{+}+ \<\mu^-,\phi_-\>_{-},$$
where $\<\,,\,\>_{\pm}$ is the dual paring extending $\<\,,\,\>^{\pm}_{0}$. Equip
\begin{equation}\label{Def:mathbfH}
\mathbf{H}_{-\alpha}: =\{\mu^+ \oplus \mu^-:\; \mu^{\pm}\in \H_{-\alpha}^{\pm} \}
\end{equation}
with the inner product $\<\mu^+\oplus \mu^-,\,\nu^+\oplus\nu^-\>_{-\alpha} := \<\mu^+,\nu^+\>_{-\alpha}^+ + \<\mu^-,\nu^-\>_{-\alpha}^-$. Then $\mathbf{H}_{-\alpha}$ is a separable Hilbert space which has CONS $\{(h^{(-\alpha),+}_{k},\;0)\} \cup \{(0,\;h^{(-\alpha),-}_{k})\}$ and hence has norm given by
\begin{equation}\label{E:norm_MinusAlpha_AnnihilatingSystem}
|\mu^+\oplus \mu^-|^2_{-\alpha} := \sum_{k} \left( \frac{1}{(1+\lambda^{+}_{k})^{\alpha}}\<\mu^+,\phi^{+}_{k}\>^2_{+} +
\frac{1}{(1+\lambda^{-}_{k})^{\alpha}}\<\mu^-,\phi^{-}_{k}\>^2_{-} \right).
\end{equation}

\begin{remark}\rm
\begin{description}
\item{(i)} We do not lose any information (in terms of finite dimensional distributions) by considering $\Y^{N,+}\oplus\Y^{N,-}$ rather than $(\Y^{N,+},\,\Y^{N,-})$. This is because the distribution of
$$\Big(\,\big(\Y^{N,+}_{t_1}(f_1),\,\Y^{N,-}_{t_1}(g_1)\big),\;\cdots,\;\big(\Y^{N,+}_{t_k}(f_k),\,\Y^{N,-}_{t_k}(g_k)\big)\,\Big)\,\in\,(\R^2)^{k},$$
is determined by that of
$$\Big(\,\ZZ^N_{t_1},\,\ZZ^N_{t_2},\,\cdots,\,\ZZ^N_{t_k}\,\Big)\,\in\,(\mathbf{H}_{-\alpha})^k,$$
where $k\in \mathbb{N}$, $\{f_i\}\subset \H^+_{\alpha}$ and $\{g_j\}\subset \H^-_{\alpha}$ are arbitrary.

\item{(ii)}As a matter of fact, $\mathbf{H}_{-\alpha}$ is equal to the set
 of continuous
 linear functionals on $\H_{\alpha}^+\times \H_{\alpha}^-$, where $\H_{\alpha}^+\times \H_{\alpha}^-$ is equipped with the natural linear structure inherited from $\H^{\pm}_{\alpha}$.
\end{description}
\qed
\end{remark}

For a general bounded Lipschitz domain $D\subset \R^d$, the Weyl's asymptotic law for the Neumann eigenvalues  holds (see \cite{NS05}). That is, the number of eigenvalues (counting their multiplicities) less than or equal to $x$, denoted by $\sharp\,\{k:\,\lambda_{k}\leq x\}$, satisfies
\begin{equation}\label{E:WeylLaw_CLT}
\lim_{x\to\infty}\frac{\sharp\,\{k:\,\lambda_{k}\leq x\}}{x^{d/2}}=C \quad \text{ for some constant }C=C(d,D)>0.
\end{equation}
Moreover, we have the following bounds for the eigenfunctions proved in \cite[Lemma 2.2]{zqCwtF14c}:
\begin{equation}\label{E:EigenfcnUpperBound_CLT}
\|\phi_{k}\|  \leq C\,\lambda_{k}^{d/4} \quad \hbox{ and }
\quad \int_{\partial D} \phi_{k}^2\,d\sigma \leq C\,(\lambda_{k}+1)
\end{equation}
for some $C=C(d,D)>0$.

For our fluctuation result (Theorem \ref{T:Convergence_AnnihilatingSystem}) to hold, we need the following assumption on $\{\delta_N\}$ which is stronger than Assumption \ref{A:ShrinkingRateLLN}. Roughly speaking, we require $\delta$ to decrease at a slower rate so that the fluctuations in $D_+$ propagate through $D_-$. This is a high density assumption for the particles.

\begin{assumption}\label{A:ShrinkingRateCLT}(Annihilation distance for functional CLT)
$\{\delta_N\}\subset (0,\infty)$ converges to 0 as $N\to\infty$ and $\liminf_{N\to\infty}N\,\delta_N^{2d} \in (0,\infty]$.
\end{assumption}

The following lemma tells us the space in which the fluctuation processes $\ZZ^N$ live.

\begin{lem}\label{L:StateSpace_AnnihilatingSystem}
Suppose that Assumption \ref{A:ShrinkingRateCLT} holds and that the initial position of particles in each of $\bar{D}_{\pm}$ are i.i.d with distribution
$u^{\pm}_0 (x) dx$, where
$u^{\pm}_0\in C(\bar{D}_{\pm})$. Then for any $\alpha>d$, $t\geq 0$ and $N\geq 1$  we have $\ZZ^N_t \in \mathbf{H}_{-\alpha}$.
\end{lem}

\begin{pf}
Fix any integer $N\geq 1$ and $t\geq 0$. We have, by definition,
\begin{eqnarray*}
\E \left[\Y^{N,+}_t(\phi))^2\right] = N\left( \E \left[\<\phi,\X^{N,+}_t\>^2 \right]-\Big(\E \left[ \<\phi,\X^{N,+}_t\> \right]\Big)^2\,\right) \leq N\,\|\phi\|^2.
\end{eqnarray*}
Using the definition of the norm $|\,\cdot\,|_{-\alpha}$ is defined in (\ref{E:norm_MinusAlpha_AnnihilatingSystem}), we have $\E[|\ZZ^N_t|_{-\alpha}^2]<\infty $ provided that
$$\sum_{k\geq 1} \left( \frac{\|\phi^{+}_{k}\|^2}{(1+\lambda^{+}_{k})^{\alpha}} +
\frac{\|\phi^{-}_{k}\|^2}{(1+\lambda^{-}_{k})^{\alpha}} \right) <\infty,$$
which is true if $\alpha>d$, using the Weyl's law (\ref{E:WeylLaw_CLT}) and the bound (\ref{E:EigenfcnUpperBound_CLT}).
\end{pf}

Suppose the initial position of particles in each of $\bar{D}_{\pm}$ are i.i.d with distribution $u^{\pm}_0\in C(\bar{D}_{\pm})$. It is easy to check that if $\alpha>d/2$, then $\ZZ^N_0\in \mathbf{H}_{-\alpha}$; furthermore,
\begin{equation}\label{E:Limit_ZZ0}
\ZZ^N_0 \toL \ZZ_0:= \Y^{+}_0 \oplus\Y^-_0 \quad \text{in }\mathbf{H}_{-\alpha},
\end{equation}
where $\Y^{\pm}_0$ is the centered Gaussian random variable in $\H^{\pm}_{-\alpha}$ with covariance
    $$\tilde{\E}\left[\Y^{\pm}_0(\phi)\Y^{\pm}_0(\psi) \right]
= \<\phi\,\psi,\,u^{\pm}_0\>_{\rho_{\pm}} - \<\phi,u^{\pm}_0\>_{\rho_{\pm}}\,\<\psi,u^{\pm}_0\>_{\rho_{\pm}}.$$
Here $\<\phi,\,\psi\>_{\rho_{\pm}}:= \int_{D_{\pm}} \phi(x)\psi(x)\,\rho_{\pm}(x)dx$ is the inner product of $L^2(D_{\pm},\,\rho_{\pm}(x)dx)$.
The main goal of this paper
is to show that the sequence of processes $\{\ZZ^N\}$ converges as $N\to\infty$, and to characterize the limit.

\section{Main results and key ideas}

\subsection{Main results}

Before stating the fluctuation result, we first define an evolution operator (see \cite{rfChZ95}) $\{Q_{s,t}\}_{s\leq t}$ as follows:
Fix any $\phi_{\pm}\in C(\bar{D}_{\pm})$ and $t>0$. Consider the following system of backward heat equations for $(v^+(s,x),\,v^-(s,y))$
for $s\in(0,t)$ with terminal data $v^{\pm}(t)=\phi_{\pm}$ and  with \emph{nonlinear and coupled} boundary conditions:
    \begin{equation}\label{E:FluctuationPDE_limit}
        \left\{\begin{aligned}
        -\,\dfrac{\partial v^+(s,x)}{\partial s}           &= \A^+ v^+(s,x)    & &\qquad\text{on } (0,t)\times D_+  \\
        \dfrac{\partial v^+(s,z)}{\partial \vec{n}_+}           &=  \lambda(z)\,\Big(v^+(s,z)+v^-(s,z)\Big)\,u_-(s,z)\,\rho_-(z)\,\1_{\{I\}}   & &\qquad\text{on } (0,t)\times \partial D_+  \\
        -\,\dfrac{\partial v^-(x,y)}{\partial s}           &= \A^- v^-(s,y)    & &\qquad\text{on } (0,t)\times D_-\\
        \dfrac{\partial v^-(s,z)}{\partial \vec{n}_-}           &=  \lambda(z)\,\Big(v^+(s,z)+v^-(s,z)\Big)\,u_+(s,z)\,\rho_+(z)\,\1_{\{I\}}  & &\qquad\text{on } (0,t)\times \partial D_-,  \\
        \end{aligned}\right.
    \end{equation}
where $(u_+,\,u_{-})$ is the hydrodynamic limit in Theorem \ref{T:Conjecture_delta_N_CLT}, $\vec{n}_{\pm}$ is the inward unit normal of $D_{\pm}$ and $\1_{\{I\}}$ is the indicator function on the interface $I$. Let $Q_{s,t}(\phi_+,\phi_-):=(v^+(s),v^-(s))$ be the solution\footnote{See Proposition \ref{prop:Integral_BackwardCoupleEqtQ} for the existence and uniqueness of solution for (\ref{E:FluctuationPDE_limit}) in $C([0,t]\times \bar{D}_+)\times C([0,t]\times \bar{D}_-)$.} for (\ref{E:FluctuationPDE_limit}) and define
$$(\mathbf{U}_{(t,s)}\mu)(\phi_+,\phi_-) := \mu\Big(Q_{s,t}(\phi_+,\phi_-)\Big)$$
for $\alpha\geq 0$, $\mu\in \mathbf{H}_{-\alpha}$ and $(\phi_+,\phi_-)$ whenever $Q_{s,t}(\phi_+,\phi_-) \in \H_{\alpha}^+\times \H_{\alpha}^-$. 
This equation is in a sense  a `linearization' of the hydrodynamic equation in Theorem \ref{T:Conjecture_delta_N_CLT}. It is the transportation component of $\ZZ$ in Theorem \ref{T:Convergence_AnnihilatingSystem}.

We are now in the position to state our main result in this chapter.
\begin{thm}\label{T:Convergence_AnnihilatingSystem}\textbf{(Fluctuation limit) }
    Suppose that Assumptions \ref{A:SettingCLT} to \ref{A:The annihilation potentialCLT} hold, and that Assumption \ref{A:ShrinkingRateCLT} holds.
    Suppose the initial position of particles in $\bar{D}_{\pm}$ are i.i.d with distribution $u^{\pm}_0(x) dx$, where
$u^{\pm}_0\in C(\bar{D}_{\pm})$. Then for any $T>0$, there exists a constant $C=C(D_+,D_-,T)>0$ such that
    $$\ZZ^N \toL \ZZ \quad \text{in  } D([0,\,T_0],\,\mathbf{H}_{-\alpha})$$
    for $\alpha>d+2$, where $T_0 :=\, T\,\wedge \,(\|u^+_0\|\vee \|u^-_0\|)^{-2}\,C$ and $\ZZ$ is the generalized Ornstein-Uhlenbeck process given by
    \begin{equation}\label{E:GeneralizedOUformula_AnnihilatingSystem}
        \ZZ_t \equalL \mathbf{U}_{(t,0)}\ZZ_0+\int_0^t \mathbf{U}_{(t,s)}\,dM_s \quad\text{in  } D([0,\,T_0],\;\mathbf{H}_{-\alpha}).
    \end{equation}
    In (\ref{E:GeneralizedOUformula_AnnihilatingSystem}), $M$ is a (unique in distribution) continuous, square integrable, $\mathbf{H}_{-\alpha}-$valued Gaussian martingale with independent increments and covariance functional characterized by
    \begin{eqnarray}\label{E:QuadVarMN_CLT}
     \<M^{N}(\phi_+,\phi_-)\>_t   &=&  \int_0^t \bigg(\<\textbf{a}_+\nabla \phi_+ \cdot \nabla \phi_+,\;u_+(s)\>_{\rho_+} + \<\textbf{a}_-\nabla \phi_- \cdot \nabla \phi_-,\;u_-(s)\>_{\rho_-}  \notag\\
        && \qquad \qquad +\int_{I}\lambda\,(\phi_++\phi_-)^2\,u_+(s)u_-(s)\,\rho_+\,\rho_-\,d\sigma \bigg)\;ds,
    \end{eqnarray}
    where $\<M(\phi_+,\phi_-)\>_t$ is the predictable quadratic variation of the real martingale $M_t(\phi_+,\phi_-)$, the pair
    $(u_+(s),\,u_-(s))$ is the hydrodynamic limit given by Theorem \ref{T:Conjecture_delta_N_CLT}, and $\ZZ_0:= \Y^{+}_0 \oplus\Y^-_0$ is the centered Gaussian random variable in (\ref{E:Limit_ZZ0}) defined on the same probability space $(\tilde{\Omega},\,\tilde{\mathcal{F}},\,\tilde{\mathcal{F}}_{t},\,\tilde{\P})$ as $M$, with $\{M,\,\Y^+_0,\,\Y^-_0\}$ being independent.
    \end{thm}

\begin{remark}\label{Rk:TransportationZZ}\rm
Observe that the representation (\ref{E:GeneralizedOUformula_AnnihilatingSystem}) of $\ZZ$ tells us that $\ZZ$ is the sum of two independent Gaussian processes,
hence is Gaussian.
 The covariance structure of $\ZZ$ is completely characterized; hence the distribution of $\ZZ$ in $ D([0,\,T_0],\,\mathbf{H}_{-\alpha})$ is uniquely determined. Moreover, the coupled PDE (\ref{E:FluctuationPDE_limit}) describes the `transportation' for the fluctuation limit $\ZZ$, and $M$ defined above describes the `driving noise'. Formally, (\ref{E:FluctuationPDE_limit}) is obtained from  (\ref{E:BackwardCoupleEqtQN}), and (\ref{E:QuadVarMN_CLT}) is obtained from (\ref{E:QuadVarMN_CLT_N}), both by letting $N\to\infty$. \qed
\end{remark}

As mentioned in Remark \ref{Rk:TransportationZZ}, the limiting process $\ZZ$ is a Gaussian. Moreover, we obtain the following properties for the limiting process.

    \begin{thm}\label{T:Properties_limit_AnnihilatingSystem}\textbf{(Properties of $Z$)}
        The fluctuation limit $\ZZ$ in Theorem \ref{T:Convergence_AnnihilatingSystem} is a continuous Gaussian Markov process which is uniquely determined in distribution, and $\ZZ$ has a version in $C^{\gamma}([0,T_0],\,\mathbf{H}_{-\alpha})$ (i.e. H\"older continuous with exponent $\gamma$) for any $\gamma\in(0,1/2)$.
    \end{thm}

We omit the proof of Theorem \ref{T:Properties_limit_AnnihilatingSystem},
 since it follows from that of Theorem 4.11 in \cite{zqCwtF14c} and the covariance structure of $\ZZ$ given by Theorem \ref{T:Convergence_AnnihilatingSystem}. Roughly speaking, the Makov property follows from the evolution property of $\mathbf{U}_{(t,s)}$ and the independent increments of the differentials. In particular, the exponent of the H\"older continuity for $\ZZ$ follows from Lemma \ref{L:ConvergenceOfUNZN0} and Theorem \ref{T:ConvergenceStochInt_CLT}.

\begin{remark}\rm
\begin{description}
\item{(i)}
        Observe that the limiting process $\ZZ=\Y^{+}\oplus  \Y^{-}$ for some processes $\Y^{\pm}$ taking values in
$\H^{\pm}_{-\alpha}$ when $\alpha$ is large enough, since it has state space   $\mathbf{H}_{-\alpha}$. Theorem \ref{T:Convergence_AnnihilatingSystem} implies that $\{\Y^{+}_t(\phi_+)+\Y^{-}_t(\phi_-):\;t\geq 0,\,\phi_+\in \H^{+}_{\alpha},\phi_-\in\H^{-}_{\alpha}\}$ is a Gaussian system. Since we can choose $\phi_{\pm}$ to be identically 0, we can strengthen the previous statement to be:
        $$\left(\,\Y^{+}_{s_1}(\phi^+_1),\cdots,\,\Y^{+}_{s_k}(\phi^+_k),\,\Y^{-}_{t_1}(\phi^-_1),\cdots,\,\Y^{-}_{t_{\ell}}(\phi^-_{\ell})\,\right)$$
        is a centered Gaussian vector in $\R^{k+\ell}$ for any $k,\ell\in \mathbb{N}$, $\{s_i\}_{i=1}^k \subset [0,T]$, $\{t_j\}_{j=1}^{\ell} \subset [0,T]$, $\{\phi^+_i\}_{i=1}^k\subset \H^{+}_{\alpha}$ and $\{\phi^-_j\}_{j=1}^{\ell}\subset \H^{-}_{\alpha}$.
\item{(ii)}
        Moreover, $M$ can be decomposed as
        $$ M \equalL M^+\oplus M^- \;in\, C([0,T_0],\,\mathbf{H}_{-\alpha}),$$
        where $M^{\pm}=(M^{\pm}_t)_{t\geq 0}$ is a continuous $\H^{\pm}_{-\alpha}$-valued Gaussian martingale with independent increment and with covariance functionals
        \begin{eqnarray*}
        \< M^{+}(\phi)\>_t &=& \int_0^{t}\<\textbf{a}_+\nabla \phi \cdot \nabla \phi,\;u_+(r)\>_{\rho_+}\,+\,\int_{I}\lambda\,\phi^2\,u_+(r)u_-(r)\,\rho_+\,\rho_-\,d\sigma\,dr,\\
        \< M^{-}(\psi)\>_t &=& \int_0^{t} \<\textbf{a}_-\nabla \psi \cdot \nabla \psi,\;u_-(r)\>_{\rho_-} \,+\,\int_{I}\lambda\,\psi^2\,u_+(r)u_-(r)\,\rho_+\,\rho_-\,d\sigma\,dr \quad\text{and}\\
        \E[M^{+}_s(\phi)\,M^{-}_t(\psi)]&=& \int_0^{s\wedge t}\int_{I}\lambda\,\phi\psi\,u_+(r)u_-(r)\,\rho_+\,\rho_-\,d\sigma\,dr.
        \end{eqnarray*}
\end{description}
\qed
\end{remark}

\subsection{Idea of proof}

Our starting point for the study of fluctuation is the following result proved in \cite{zqCwtF14b}. Let us recall it here for the convenience of the readers.
\begin{lem}\label{L:KeyMtgAnnihilatingDiffusionModel}
For any $\phi_{\pm}\in Dom(\A^{\pm})$, we have
    \begin{eqnarray*}
        && \<\phi_+,\X^{N,+}_t\>+\<\phi_-,\X^{N,-}_t\> - \<\phi_+,\X^{N,+}_0\>
  - \<\phi_-,\X^{N,-}_0\>\\
     && -\int_0^t
        \<\A^+\phi_+,\,\X^{N,+}_s\>+\<\A^-\phi_-,\,\X^{N,-}_s\>- \<\ell_{\delta_N}(\phi_++\phi_-),\,\X^{N,+}_s\otimes \X^{N,-}_s\>\,ds\\
    \end{eqnarray*}
is an $\F^{(\X^{N,+},\X^{N,-})}_t$-martingale with predictable quadratic variation
    \begin{eqnarray*}
        &&  \frac{1}{N}\int_0^t \<\textbf{a}_+\nabla \phi_+ \cdot \nabla \phi_+,\,\X^{N,+}_s\> + \<\textbf{a}_-\nabla \phi_- \cdot \nabla \phi_-,\,\X^{N,-}_s\>    +\<\ell_{\delta_N}(\phi_++\phi_-)^2,\,\X^{N,+}_s\otimes\X^{N,-}_s\>\;ds.
    \end{eqnarray*}
Here   $\<\varphi(x,y),\,\mu^+(dx) \otimes \mu^-(dy)\>:= \frac{1}{N^2}\sum_{i}\sum_{j}\varphi(x_i,y_j)$ when $\mu=(\frac{1}{N}\sum_{i} \1_{x_i},\,\frac{1}{N}\sum_{j}\1_{y_j})$.
\end{lem}

\medskip

Recall that $\ZZ^N_t := \Y^{N,+}_t \oplus \Y^{N,-}_t$. Hence Lemma \ref{L:KeyMtgAnnihilatingDiffusionModel} reads as
\begin{equation}\label{E:Ch6_KeyStep_1}
\ZZ^N_t-\ZZ^N_0=
\int_0^t \left( \mathbf{A} \ZZ^N_s -  K^N_s\right) ds + M^N_t,
\end{equation}
where
\begin{eqnarray*}
\mathbf{A}\mu(\phi_+,\phi_-) &:=& \mu(\A^{+}\phi_+,\,\A^{-}\phi_-),\\
K^N_s(\phi_+,\phi_-) &:=&  \sqrt{N}\Big(\,\<\ell_{\delta_N}(\phi_++\phi_-),\,\X^{N,+}_s\otimes \X^{N,-}_s\>- \E[\<\ell_{\delta_N}(\phi_++\phi_-),\,\X^{N,+}_s\otimes \X^{N,-}_s\>]\,\Big).
\end{eqnarray*}
and $M^{N}_t(\phi_+,\phi_-)$ is a real valued $\F^{(\X^{N,+},\,\X^{N,-})}_t$-martingale with predictable quadratic variation
    \begin{eqnarray}\label{E:QuadVarMN_CLT_N}
     \<M^{N}(\phi_+,\phi_-)\>_t   &=&  \int_0^t \<\textbf{a}_+\nabla \phi_+ \cdot \nabla \phi_+,\,\X^{N,+}_s\> + \<\textbf{a}_-\nabla \phi_- \cdot \nabla \phi_-,\,\X^{N,-}_s\>   \notag\\
        && \qquad \qquad +\<\ell_{\delta_N}(\phi_++\phi_-)^2,\,\X^{N,+}_s\otimes\X^{N,-}_s\>\;ds.
    \end{eqnarray}

The \emph{key idea} is to rewrite (\ref{E:Ch6_KeyStep_1}) as
\begin{equation}\label{E:KeyStepProof}
\ZZ^N_t-\ZZ^N_0= \int_0^t (\mathbf{A}- \mathbf{B}^N_s)\ZZ^N_s\,ds + \int_0^t (\mathbf{B}^N_s\ZZ^N_s- K^N_s)\,ds + M^N_t,
\end{equation}
in which $\mathbf{B}^N_s\mu(\phi_+,\phi_-)$ is defined as
\begin{equation}\label{Def:mathbf{B}^N_s}
\<\ell_{\delta_N}(\phi_+ + \phi_-),\,f^{N,+}_s\>_{\rho_+} \,
\mu \big(\<\ell_{\delta_N}(\phi_+ + \phi_-),\,f^{N,+}_s\>_{\rho_-},\; \<\ell_{\delta_N}(\phi_+ + \phi_-),\,f^{N,-}_s\>_{\rho_+} \big).
\end{equation}
This expression is inserted to the right-hand side of \eqref{E:KeyStepProof} to, roughly speaking, project $K^N_s$ onto the image of $\ZZ^N_s$. Here $(f^+,f^-):=(f^{N,+},f^{N,-})$ is defined to be the unique element in $C([0,\infty)\times \bar{D}_+)\times C([0,\infty)\times \bar{D}_-)$ satisfying the coupled integral equations
    \begin{equation}\label{E:CoupledIntegral_fN}
    \left\{\begin{aligned}
    f^+_t(x) &= P^+_tu^+_0(x)-\int_0^tP^+_{t-r}\Big(f^+_r(\cdot)\,\int_{D_-} \ell_{\delta_N}(\cdot,y)f^-_r(y)dy\Big)(x)\,dr  \\
    f^-_t(y) &= P^-_tu^-_0(y)-\int_0^tP^-_{t-r}\Big(f^-_r(\cdot)\,\int_{D_+} \ell_{\delta_N}(x,\cdot)f^+_r(x)dx\Big)(y)\,dr,
    \end{aligned}\right.
    \end{equation}
where $\{P^{\pm}_t\}_{t\geq 0}$ is the transition semi-group for $X^{\pm}$ and $P^{\pm}_{t-r}$ acts on the dot variable. The existence and uniqueness of $(f^+,f^-)$ can be checked by the same fixed point argument as that for $(u_+,u_-)$ in Proposition 2.19 in \cite{zqCwtF13a}. We will show in Lemma \ref{L:UConvergencefNtoU} that $(f^{N,+},\,f^{N,-})$ converges to $(u_+,\,u_-)$ as $N\to\infty$. Intuitively, both $(f^{N,+},f^{N,-})$ and $(u_+,u_-)$ are approximations to $(\X^{N,-},\,\X^{N,-})$, but $(f^{N,+},f^{N,-})$ is
the better one.

One of the most challenging task in the proof is to show that, in an appropriate sense,
\begin{equation*}
\int_0^t (\mathbf{B}^N_s\ZZ^N_s- K^N_s)\,ds \rightarrow 0 \text{ when }N\to\infty;
\end{equation*}
that is, we can replace $\int K^N_s\,ds$ by $\int \mathbf{B}^N_s\ZZ^N_s\,ds$ in (\ref{E:Ch6_KeyStep_1}). This is basically step 6 in the `Outline of proof' below.

We discovered the formula (\ref{Def:mathbf{B}^N_s}) of $\mathbf{B}^N_s$, roughly speaking, by projecting $K^N_s$ onto the image of $\ZZ^N_s$. This inspiration comes from the well-known \emph{Boltzman-Gibbs principle} in mathematical physics. The principle says that the fluctuation fields of non-conserved quantities change on a time scale much faster than the conserved ones, hence in a time integral only the component along those fields of conserved quantities survive. This idea leads us to reasonably hope
 that  $\int_0^t (\mathbf{B}^N_s\ZZ^N_s- K^N_s)\,ds \rightarrow 0$, which is confirmed in Theorem \ref{T:BGprinciple}. Analytically, the proof of $\int_0^t (\mathbf{B}^N_s\ZZ^N_s- K^N_s)\,ds \rightarrow 0$ stems from a `magical cancelation' (see (\ref{E:MagicCancel_1}) and (\ref{E:MagicCancel_2}) in the proof of Theorem \ref{T:BGprinciple}) for the first \emph{two} terms of the asymptotic expansion of the correlation functions.

The rest part of the paper is devoted to the proof of Theorem \ref{T:Convergence_AnnihilatingSystem}.  

\noindent
\textbf{Outline of proof: }
We prove Theorem \ref{T:Convergence_AnnihilatingSystem} through the following
six steps.

\begin{enumerate}
\item[] Step 1\;  $\ZZ^N$ satisfies the following stochastic integral equation:
                \begin{eqnarray*}
                        \ZZ^{N}_t= \mathbf{U}^N_{(t,0)}\ZZ^{N}_0+ \int_0^t \mathbf{U}^N_{(t,s)}\,dM^N_s + \int_0^t \mathbf{U}^N_{(t,s)}(\mathbf{B}^N_s\ZZ^N_s- K^N_s)\,ds\;\text{ a.s.},
                \end{eqnarray*}
                \qquad \qquad where $\mathbf{U}^N_{(t,s)}$ is the evolution system generated by $\mathbf{A}-\mathbf{B}^N_s$; see Theorem \ref{T:3.3_CLT}.
\item[]  Step 2\; $\ZZ^N$ is tight in $D([0,T_0],\,\mathbf{H}_{-\alpha})$; see Theorem \ref{T:Tightness_Z}.
\item[]  Step 3\; $M^N \toL M$ in $D([0,T_0],\,\mathbf{H}_{-\alpha})$; see Theorem \ref{T:ConvergenceOfM^N_Annihilation}.
\item[]  Step 4\; $\mathbf{U}^N_{(t,0)}\ZZ^N_0 \toL \mathbf{U}_{(t,0)}\ZZ_0$ in $D([0,T_0],\,\mathbf{H}_{-\alpha})$; see Lemma \ref{L:ConvergenceOfUNZN0}.
\item[]  Step 5\; $\int_0^t \mathbf{U}^N_{(t,s)}\,dM^N_s \toL  \int_0^t \mathbf{U}_{(t,s)}\,dM_s$ in $D([0,T_0],\,\mathbf{H}_{-\alpha})$; see Theorem \ref{T:ConvergenceStochInt_CLT}.
\item[]  Step 6\; $\int_0^t \mathbf{U}^N_{(t,s)}(\mathbf{B}^N_s\ZZ^N_s- K^N_s)\,ds \;\to\; 0 $ in $D([0,T_0],\,\mathbf{H}_{-\alpha})$; see Theorem \ref{T:BGprinciple}.
\end{enumerate}

This rough outline is the same as that for the single species model in \cite{zqCwtF14c}. In fact, with all the preliminary estimates in Section 2, and with the asymptotic expansion of the correlation functions (Theorem \ref{T:Asymptotic_nm_correlation_t}) proved in Section 3, all the steps except Step 2 and Step 6 can be treated using the method in \cite{zqCwtF14c}.
Some of the main efforts are directed toward Step 2 and Step 6 which require asymptotic analysis of the correlation functions (section 6.2) and the generalized correlation functions (section 6.3) respectively.

\section{Proofs}

\textbf{Convention: } To avoid unnecessary complications, we assume, from now on, that $\lambda=\widehat{\lambda}=1$ and that the underlying motion of the particles are reflected Brownian motions (i.e.  $\rho_{\pm}=1$ and $\textbf{a}_{\pm}$ are identity matrices). However, our arguments work for general symmetric reflected diffusions as in \cite{zqCwtF14c} and for any continuous functions $\lambda(z)\in C(I)$ as in \cite{zqCwtF14b}.
When there is no danger of confusion, for each fixed $N$, we
write $\ell(x,y)$ in place of $\ell_{\delta_N}(x,y)$ for simplicity. The constant $C_0$ is always equal to $C_0:=\|u^+_0\|\vee \|u^-_0\|$. The
minimal augmented filtration $\F^{(\X^{N,+},\,\X^{N,-})}_t$ generated by
 the annihilating diffusion process will be abbreviated as $\F^N_t$. Assumptions \ref{A:SettingCLT} to \ref{A:ShrinkingRateLLN} are in force throughout
the rest of the paper,
and we will indicate explicitly whenever Assumption \ref{A:ShrinkingRateCLT} is invoked.

\subsection{Preliminaries}
\subsubsection{Transition densities of reflected diffusions}

It is well known (cf. \cite{BH91, GSC11}) that the $(\A,\,\rho)$-reflected diffusion $X$ in Definition \ref{Def:ReflectedDiffusion} has a transition density $p(t,x,y)$ with respect to  the symmetrizing measure $\rho(x)dx$ (i.e., $\P_x(X_t\in dy)=p(t,x,y)\,\rho(y)dy$ ) satisfying $p(t,x,y)=p(t,y,x)$, that $p$ is locally H\"older continuous and hence $p\in C((0,\infty)\times \bar{D}\times \bar{D})$, and that we have two-sided Gaussian bounds: for any $T>0$, there are constants $c_1,\,c_2\geq 1$ such that
\begin{equation}\label{E:Gaussian2SidedHKE}
\dfrac{1}{c_1 t^{d/2}}\,\exp\left(- \frac{c_2 |y-x|^2}{t}\right)
\leq p(t,x,y) \leq \dfrac{c_1}{t^{d/2}}\,\exp\left(- \frac{|y-x|^2}{c_2\,t}\right)
\end{equation}
for every $(t,x,y)\in(0,T]\times \bar{D}\times \bar{D}$. Using (\ref{E:Gaussian2SidedHKE}) and the Lipschitz assumption
of $D$,
we can check that
\begin{eqnarray}
\sup_{x\in\bar{D}}\,\sup_{0<\delta\leq \delta_0}\,\frac{1}{\delta}\int_{D^{\delta}}p(t,x,y)\,dy &\leq& \frac{C}{\sqrt{t}}
\quad \hbox{for } t\in(0,T]\quad \text{and} \label{E:boundary_strip_boundedness}\\
\sup_{x\in\bar{D}}\,\int_{\partial D}p(t,x,y)\,\sigma(dy) &\leq& \frac{C}{\sqrt{t}} \quad \hbox{for } t\in(0,T], \label{E:Surface_integral_boundedness}
\end{eqnarray}
where $C,\,\delta_0>0$ are constants that depend only on $d$, $T$, the Lipschitz characteristics of $D$, the ellipticity of $\textbf{a}$
and the lower and upper bound of $\rho$. Here $D^\delta:=\{x\in D: \hbox{dist} (x, \partial D)<\delta\}$. Therefore, under Assumption \ref{A:SettingCLT}, the transition density $p^{\pm}(t,x,y)$ of $X^{\pm}$ (with respect to $\rho_{\pm}$) satisfies \eqref{E:Gaussian2SidedHKE}, \eqref{E:boundary_strip_boundedness} and \eqref{E:Surface_integral_boundedness}.  Observe that
$$\int_{z\in D_+}\int_{y\in D_-}\ell(z,y)\,p^+(s,x,z)\,dy\,dz \leq  \frac{1}{\delta^{d+1}}\int_{D_-^{\delta}}\int_{B(z,\delta)\cap D_+}p^+(s,x,z)\,dz\,dy,$$
where $B(z,\delta)$ is the ball of radius $\delta$ centered at $z$. Hence by \eqref{E:boundary_strip_boundedness}, we have
\begin{equation}\label{E:Idelta_bound}
\sup_{x\in \bar{D}_+}\int_{z\in D_+}\int_{y\in D_-}\ell(z,y)\,p^+(s,x,z)\,dy\,dz \leq  \frac{C(d,D_+,T)}{\sqrt{s}}\quad\text{for }s\in(0,T],
\end{equation}
whenever $N\geq N_0(d,D_+)$. A similar inequality holds for $p^-$.

\subsubsection{Minkowski content}

We will make extensive use of the following result about Minkowski content of the interface $I$. It is
established  in  \cite{zqCwtF14b} and
restated here for the convenience of the reader.

\begin{lem}\label{L:MinkowskiContent_I_k_CLT}
Suppose Assumptions \ref{A:SettingCLT}, \ref{A:ParameterAnnihilationCLT} and \ref{A:The annihilation potentialCLT} hold. Suppose $k\in\mathbb{N}$ and $\F\subset  C\big((\bar{D}_+\times \bar{D}_-)^k\big)$ is an equi-continuous and uniformly bounded family of functions on $(\bar{D}_+\times \bar{D}_-)^k$. Then as $\delta\to 0$, we have
\begin{eqnarray*}
&& \int_{(x_1,y_1)\in D_+\times D_-}\cdots \int_{(x_k,y_k)\in D_+\times D_-} f(x_1,y_1,\cdots,x_k,y_k)\prod_{i=1}^k\ell_{\delta}(x_i,y_i)\;d(x_1,y_1,\cdots,x_k,y_k)\\
&\to& \int_{z_1\in I}\cdots \int_{z_k\in I} f(z_1,z_1,\cdots,z_k,z_k)\,)\prod_{i=1}^k \lambda(z_i)\,d\sigma(z_1)\cdots d\sigma(z_k)
\end{eqnarray*}
uniformly for $f\in \F$.
\end{lem}

\subsubsection{Three sets of coupled equations}\label{subsection:ThreeCoupledEqts}

Recall that $(f^+,f^-)=(f^{N,+},f^{N,-})$ is the deterministic pair solving (\ref{E:CoupledIntegral_fN}). In this subsection, we will construct two more coupled integral equations that is the core in the study of fluctuations of the annihilating diffusion system. For each $N\in\mathbb{N}$, the solutions of them will be denoted by $(G^{N},\,G^{N,+},\,G^{N,-})$ and $(g^{N,+},\,g^{N,-})$ respectively. We will suppress the notation $N$ and write  $(g^+,g^-)$ in place of $(g^{N,+},\,g^{N,-})$, etc.

We first prove that $(f^+,f^-)$ is a good approximation to $(u_+,\,u_-)$.
\begin{lem}\label{L:UConvergencefNtoU}
$|f^{N,\pm}|$ is uniformly bounded above by $\|u^{\pm}_0\|$. Moreover, For each $t \geq 0$, we have $f^{N,\pm}_t $converges uniformly on $\bar{D}_{\pm}$ to $u_{\pm}(t)$, as $N\to\infty$.
\end{lem}

\begin{pf}
Clearly, $\sup_{(t,x)}\sup_{N}|f^{N,\pm}(t,x)|\leq \|u^{\pm}_0\|$. This can be seen, for example, by the probabilistic representations of $(f^+_t(x),\,f^-_t(y))$ given by
    \begin{equation}\label{E:ProbabilisticRep_fN}
    \left\{\begin{aligned}
    f^+_t(x) &= \E^{x} \Big[\,u^+_0(X^+_t)\,\exp{\Big( -\int^t_0\int_{D_-} \ell(X^+_s,y)\,f^-_{t-s}(y)\,dy\,ds \Big)}\,\Big] \\
    f^-_t(y) &= \E^{y} \Big[\,u^-_0(X^-_t)\,\exp{\Big( -\int^t_0\int_{D_+} \ell(x,X^-_s)\,f^+_{t-s}(x)\,dx\,ds \Big)}\,\Big].
    \end{aligned}\right.
    \end{equation}
We now show that $\{(f^{N,+}_t,\,f^{N,-}_t)\}_{N\geq1}$ is an equi-continuous sequence in $C(\bar{D}_+)\times C(\bar{D}_-)$. This can be achieved by using \eqref{E:Idelta_bound} and the H\"older continuity of $p^{\pm}$ (cf. \cite[Chapter 3]{GSC11}) as follows.

Fix $\epsilon>0$. By \eqref{E:Idelta_bound}, there exists $t_{\ast}\in(0,t)$ and $N_0(d,D_+)$such that
$$ 2\,\|u^+_0\|\,\|u^-_0\|\int_0^{t_{\ast}}  \Big(\sup_{x \in \bar{D}_+}\int_{z\in D_+}\int_{y\in D_-}\ell(z,y)\, p^+(t-r,x,z)\,dy\,dz\Big)\,dr < \frac{\epsilon}{2}\quad\text{for } N\geq N_0.$$
From (\ref{E:CoupledIntegral_fN}), we have, for any $x_1,\,x_2\in \bar{D}_+$,
\begin{eqnarray*}
&& |f^+_t(x_1)-f^+_t(x_2)| \\
&=& \bigg|\int_0^t \int_{D_+}\Big(p^+(t-r,x_1,z)-p^+(t-r,x_2,z)\Big)\, \left(f^+_r(x)\,\int_{D_-} \ell(z,y)f^-_r(y)\,dy\right)\,dz\,dr\bigg|\\
&\leq& \|u^+_0\|\,\|u^-_0\|\int_0^t \int_{z\in D_+}\int_{y\in D_-}\ell(z,y)\,\Big| p^+(t-r,x_1,z)-p^+(t-r,x_2,z)\Big|\,dy\,dz\,dr\\
&\leq& \frac{\epsilon}{2}\,+\,\|u^+_0\|\,\|u^-_0\|\int_{t_{\ast}}^t \int_{z\in D_+}\int_{y\in D_-}\ell(z,y)\,\Big| p^+(t-r,x_1,z)-p^+(t-r,x_2,z)\Big|\,dy\,dz\,dr\\
&\leq&  \frac{\epsilon}{2}\,+\,\|u^+_0\|\,\|u^-_0\|\,\frac{\mathcal{H}^{2d}(I^{\delta})}{c_{d+1}\,\delta^{d+1}}\,\int_{t_{\ast}}^t  \frac{|x_1-x_2|^{\alpha}}{(t-r)^{\beta}}\,dr\qquad\text{for some }\alpha(d,D_+,T),\,\beta(d,D_+,T)>0\\
&\leq& \frac{\epsilon}{2}\,+\,\|u^+_0\|\,\|u^-_0\|\,(\sigma(I)+1)\int_{t_{\ast}}^t  \frac{|x_1-x_2|^{\alpha}}{(t-r)^{\beta}}\,dr\qquad\text{for } N\geq N_1(d,D_+),\,\text{ by } \eqref{E:Federer_MinkowSki_Content}.
\end{eqnarray*}
We have used the H\"older continuity of $p^{+}$ (cf. \cite[Chapter 3]{GSC11}) in the second to the last inequality. It is now clear that $\{f^{N,+}_t\}_{N\geq 1} \subset C(\bar{D}_+)$ is equi-continuous  for any $t>0$. A similar calculation applies to $f^-_t$. Hence $\{(f^{N,+}_t,\,f^{N,-}_t)\}_{N\geq1}$ is equi-continuous for any $t>0$. Finally, by comparing the probabilistic representations of $(u_+(t,x),\,u_-(t,y))$ in (\ref{E:ProbabilisticRep_CoupledPDE_CLT}) and that of $(f^+_t(x),\,f^-_t(y))$ in (\ref{E:ProbabilisticRep_fN}), we can check that any subsequential limit of $f^{N,+}_t$ is equal to $u_+(t)$ by using Lemma \ref{L:MinkowskiContent_I_k_CLT}.
\end{pf}

\medskip

Next we define $(G,\,G^+,\,G^-)=(G^N,\,G^{N,+},\,G^{N,-})$ to be the unique solution in
$C([0,\infty)\times \bar{D}_+ \times \bar{D}_-)\times C([0,\infty)\times \bar{D}_+\times \bar{D}_+)\times C([0,\infty)\times \bar{D}_-\times \bar{D}_-)$ to the coupled integral equations.
\begin{eqnarray*}
G_t(x,y)&=&-\int_0^t P^{(1,1)}_{t-r}\bigg\{\,G_r(\tilde{x},\tilde{y})\,\left(\int_{D_-}\ell(\tilde{x},w)f^{-}_r(w)dw + \int_{D_+}\ell(z,\tilde{y})f^{+}_r(z)dz\right)\\
&&\qquad\qquad +\int_{D_+}G^+_r(\tilde{x},z)\ell(z,\tilde{y})f^-_r(w)\,dz + \int_{D_-}G^-_r(\tilde{y},w)\ell(\tilde{x},w)f^+_r(\tilde{x})\,dw\\
&&\qquad\qquad -f^+_r(\tilde{x})f^-_r(\tilde{y})\ell(\tilde{x},\tilde{y})\,\bigg\}(x,y)\,dr,
\end{eqnarray*}
\begin{eqnarray*}
G^+_t(x_1,x_2)&=&-\int_0^t P^{(2,0)}_{t-r}\bigg\{\,G^+_r(\tilde{x_1},\tilde{x_2})\,\int_{D_-}[\ell(\tilde{x_1},w)+\ell(\tilde{x_2},w)]f^-_r(w)dw\\
&&\qquad\qquad +\int_{D_-}f^+_r(\tilde{x_1})\ell(\tilde{x_1},w)G_r(\tilde{x_2},w)+f^+_r(\tilde{x_2})\ell(\tilde{x_2},w)G_r(\tilde{x_1},w)\,dw\,\bigg\}(x,y)\,dr
\end{eqnarray*}
and
\begin{eqnarray*}
G^-_t(y_1,y_2)&=&-\int_0^t P^{(0,2)}_{t-r}\bigg\{\,G^-_r(\tilde{y_1},\tilde{y_2})\,\int_{D_+}[\ell(z,\tilde{y_1})+\ell(z,\tilde{y_2})]f^+_r(z)dz\\
&&\qquad\qquad +\int_{D_+}f^-_r(\tilde{y_1})\ell(z,\tilde{y_1})G_r(z,\tilde{y_2})+f^-_r(\tilde{y_2})\ell(z,\tilde{y_2})G_r(z,\tilde{y_1})\,dz\,\bigg\}(x,y)\,dr,
\end{eqnarray*}
where the semigroup $P^{(i,j)}_t$ acts on the variables with a `tilde'.
\begin{remark}\rm
It is clear from the definition that $G^{\pm}$ is symmetric; that is, $G^{+}(x_1,x_2)=G^{+}(x_2,x_1)$ and $G^{-}(y_1,y_2)=G^{-}(y_2,y_1)$. The term $f^+_r(\tilde{x})f^-_r(\tilde{y})\ell(\tilde{x},\tilde{y})$ in the equation for $G$ guarantees that $(G,\,G^+,\,G^-)$ cannot be constantly zero, even though they are zero when $t=0$. This non-negative term contributes to the creation of fluctuation near the $I$.
\end{remark}

Finally, $(g^+,\,g^-)=(g^{N,+},\,g^{N,-})$ is defined to be the unique solution in $C([0,\infty)\times \bar{D}_+)\times C([0,\infty)\times \bar{D}_-)$  to the following coupled integral equations:
\begin{eqnarray*}
g^+_t(x)&=&-\int_0^t P^{+}_{t-r}\bigg\{\int_{D_-}\ell(\tilde{x},w)\,[\,g^+_r(\tilde{x})f^-_r(w)+g^-_r(w)f^+_r(\tilde{x})+G_r(\tilde{x},w)\,]\,dw\,\bigg\}(x)\,dr\\
g^-_t(y)&=&-\int_0^t P^{-}_{t-r}\bigg\{\int_{D_+}\ell(z,\tilde{y})\,[\,g^+_r(z)f^-_r(\tilde{y})+g^-_r(\tilde{y})f^+_r(z)+G_r(z,\tilde{y})\,]\,dz\,\bigg\}(y)\,dr,
\end{eqnarray*}
where the semigroups $P^{+}_{t-r}$ and $P^{-}_{t-r}$ act on $\tilde{x}$ and $\tilde{y}$ respectively.

\begin{remark}\rm
The functions  $(G^{N},\,G^{N,+},\,G^{N,-})$ and $(g^{N,+},\,g^{N,-})$ appear in the second order term in the asymptotic expansion of the correlation functions in Theorem \ref{T:Asymptotic_nm_correlation_t}. Their definitions are motivated and justified by the hierarchy (\ref{E:BBGKY_nm_alpha_t}).
It turns out that, as in Corollary \ref{cor:Asymptotic_nm_correlation_t_2}, the covariance structure of $\ZZ^N$ involves $(G^{N},\,G^{N,+},\,G^{N,-})$ but not $(g^{N,+},\,g^{N,-})$. \qed
\end{remark}

Although for fixed $N$, the supremum norms for $G$ and $G^{\pm}$  are finite,
unlike the cases in \cite{pD88b, pK86, pK88}, these norms
become unbounded as $N\to\infty$\footnote{In fact, we can check, using the probabilistic representation of $G$ (cf. the proof of the lemma below) and a simple exit time estimate, that $G\to\infty$ as $N\to\infty$ on the set $\{(z,z):\,z\in I\}$, provided $\inf\,u^{\pm}_0>0$.}.
Fortunately, we still have the following bounds.
\begin{lem}\label{L:BoundsForGGG}
For any $T>0$, there exist $C=C(D_+,D_-,T)>0$ and an integer $N_0=N_0(D_+,D_-)$ such that
\begin{eqnarray}
\int \ell(\tilde{x},y)|G_t(x,y)|\,d(x,\tilde{x},y)+ \int \ell(x,\tilde{y})|G_t(x,y)|\,d(x,\tilde{y},y) &\leq & (C_0\,C)^2\,\sqrt{t}\\
\int \ell(x_1,\tilde{y})|G^+_t(x_1,x_2)|\,d(x_1,\tilde{y},x_2)&\leq & (C_0\,C)^3\,t\\
\int \ell(\tilde{x},y_1)|G^-_t(y_1,y_2)|\,d(y_1,\tilde{x},y_2)&\leq & (C_0\,C)^3\,t\\
\int |G_t(x,y)|\,d(x,y)&\leq & (C_0\,C)^2\,t\\
\int |G^+_t(x_1,x_2)|\,d(x_1,x_2)&\leq & (C_0\,C)^3\,t^{3/2}\\
\int |G^-_t(y_1,y_2)|\,d(y_1,y_2)&\leq & (C_0\,C)^3\,t^{3/2}
\end{eqnarray}
for all $t\in [0,\,T\wedge(C_0\,C)^{-2}]$ and $N\geq N_0$.
\end{lem}

\begin{pf}
Since each of $G$, $G^{+}$ and $G^{-}$ is the probabilistic solution of a heat equation, they have the following probabilistic representations (see  Proposition 2.19 of \cite{zqCwtF13a}):
\begin{eqnarray*}
G_t(x,y) &=& \int_{\theta=0}^t\E^{(x,y)}\bigg[
\bigg(
f^+_{t-\theta}(X_{\theta})f^-_{t-\theta}(Y_{\theta})\ell(X_{\theta},Y_{\theta})- \int_{D_+}G^{+}_{t-\theta}(X_{\theta},z)\ell(z,Y_{\theta})f^-_{t-\theta}(Y_{\theta})dz \\
&& \qquad \qquad \qquad - \int_{D_-}G^{-}_{t-\theta}(Y_{\theta},w)\ell(X_{\theta},w)f^+_{t-\theta}(X_{\theta})dw
\bigg) \\
&& \qquad \qquad  \cdot \exp{\Big(-\int_{s=0}^{\theta}\int_{D_-}\ell(X_s,w)f^-_{t-s}(w)dw + \int_{D_+}\ell(z,Y_s)f^+_{t-s}(z)dz\,ds\Big)}
\bigg]\,d\theta,
\end{eqnarray*}
\begin{eqnarray*}
G^+_t(x_1,x_2) &=& -\,\int_{\theta=0}^t\E^{(x_1,x_2)}\bigg[
\bigg(
f^+_{t-\theta}(X^{1}_{\theta})\int_{D_-}\ell(X^1_{\theta},w)G_{t-\theta}(X^2_{\theta},w)dw \\
&& \qquad \qquad \qquad + f^+_{t-\theta}(X^{2}_{\theta})\int_{D_-}\ell(X^2_{\theta},w)G_{t-\theta}(X^1_{\theta},w)dw
\bigg) \\
&& \qquad \qquad  \cdot \exp{\Big(-\int_{s=0}^{\theta}\int_{D_-}[\ell(X^1_s,w)+\ell(X^2_s,w)]f^-_{t-s}(w)dw\,ds\Big)}
\bigg]\,d\theta,
\end{eqnarray*}
\begin{eqnarray*}
G^-_t(y_1,y_2) &=& -\,\int_{\theta=0}^t\E^{(y_1,y_2)}\bigg[
\bigg(
f^-_{t-\theta}(Y^{1}_{\theta})\int_{D_+}\ell(z,Y^1_{\theta})G_{t-\theta}(z,Y^2_{\theta})dz \\
&& \qquad \qquad \qquad + f^-_{t-\theta}(Y^{2}_{\theta})\int_{D_+}\ell(z,Y^2_{\theta})G_{t-\theta}(z,Y^1_{\theta})dz
\bigg) \\
&& \qquad \qquad  \cdot \exp{\Big(-\int_{s=0}^{\theta}\int_{D_+}[\ell(z,Y^1_s)+\ell(z,Y^2_s)]f^+_{t-s}(z)dz\,ds\Big)}
\bigg]\,d\theta,
\end{eqnarray*}
where $\{X,\,X^1,\,X^2\}$ are independent RBMs on $D_+$ and $\{Y,\,Y^1,\,Y^2\}$ are independent RBMs on $D_-$
that are also independent of $\{X,\,X^1,\,X^2\}$. Here
$\E^{(x,y)}$ denotes the expectation w.r.t. the law of $(X,Y)$ starting at $(x,y)$, etc.

Since $f^{\pm}\geq 0$ and $\|f^{\pm}\|\leq \|u^{\pm}_0\|\leq C_0$, the three formulae above give rise to the following point-wise bounds:
\begin{eqnarray*}
|G_t(x,y)|&\leq& C_0\int_{\theta=0}^t\E^{(x,y)}\bigg[
C_0\,\ell(X_{\theta},Y_{\theta})+\Big|\int_{D_+}G^+_{t-\theta}(X_{\theta},z)\ell(z,Y_{\theta})dz\Big|
\\
&& \qquad\qquad \qquad\qquad \qquad +\Big|\int_{D_-}G^-_{t-\theta}(Y_{\theta},w)\ell(X_{\theta},w)dw\Big|\,
\bigg],\\
|G^+_t(x_1,x_2)|&\leq& C_0\int_{\theta=0}^t\E^{(x_1,x_2)}\bigg[\,
\Big|\int_{D_-}G_{t-\theta}(X^2_{\theta},w)\ell(X^1_{\theta},w)dw\Big|
+\Big|\int_{D_-}G_{t-\theta}(X^1_{\theta},w)\ell(X^2_{\theta},w)dw\Big|\,
\bigg],\\
|G^-_t(y_1,y_2)|&\leq& C_0\int_{\theta=0}^t\E^{(y_1,y_2)}\bigg[\,
\Big|\int_{D_+}G_{t-\theta}(z,Y^2_{\theta})\ell(z,Y^1_{\theta})dz\Big|
+\Big|\int_{D_+}G_{t-\theta}(z,Y^1_{\theta})\ell(z,Y^2_{\theta})dz\Big|\,
\bigg].
\end{eqnarray*}

Plug in the bound for $|G^+|$ and $|G^-|$ into that of $|G|$, we have
\begin{eqnarray}\label{E:BoundsForGGG}
|G_t(x,y)|&\leq& C^2_0\,\E^{(x,y)}\int_{\theta=0}^t
\ell(X_{\theta},Y_{\theta})\, +\,C^2_0\,\E^{(x,y)}\int_{\theta=0}^t\int_{\alpha=0}^{t-\theta}\int_{z\in D_+}\int_{w\in D_-} \bigg\{ \\
&& \quad \ell(z,Y_{\theta})\,\E^{(X_{\theta},z)}
\bigg(\Big|G_{t-\theta-\alpha}(X^2_{\alpha},w)\Big|\ell(X^1_{\alpha},w)+ \Big|G_{t-\theta-\alpha}(X^1_{\alpha},w)\Big|\ell(X^2_{\alpha},w) \bigg) \notag\\
&& \quad +\ell(X_{\theta},w)\,\E^{(Y_{\theta},w)}
\bigg(\Big|G_{t-\theta-\alpha}(z,Y^2_{\alpha})\Big|\ell(z,Y^1_{\alpha})+ \Big|G_{t-\theta-\alpha}(z,Y^1_{\alpha})\Big|\ell(z,Y^2_{\alpha}) \bigg)
\,\bigg\}. \notag
\end{eqnarray}

Define $$\phi(t):= \phi^{(N)}(t):= \int \ell(\tilde{x},y)|G_t(x,y)|\,d(x,\tilde{x},y)+ \int \ell(x,\tilde{y})|G_t(x,y)|\,d(x,\tilde{y},y),$$
which serves
as an approximation to
$$\int_{x\in D_+}\int_{y\in I} |G_t(x,y)|\,dx\,d\sigma(y) + \int_{x\in I}\int_{y\in D_-} |G_t(x,y)|\,d\sigma(x)\,dy.$$

Simplifying the RHS of (\ref{E:BoundsForGGG}) using Chapman Kolmogorov equation and then applying (\ref{E:boundary_strip_boundedness}), we obtain, for $N\geq N_0(D_{\pm})$,
\begin{eqnarray*}
\phi(t) &\leq& (C_0\,C)^2\,\left(\;\sqrt{t}+ \int_{\theta=0}^t\int_{\alpha=0}^{t-\theta}\Big(\frac{1}{\theta\,\alpha}+\frac{1}{\sqrt{(\theta+\alpha)\,\alpha}}\Big)\,\phi(t-\theta-\alpha) \;\right)\\
&=&(C_0\,C)^2\,\left(\;\sqrt{t}+ (\pi+2) \int_{\alpha=0}^t\phi(t-\alpha) \;\right) \quad\text{by Fubinni's Theorem.}\\
\end{eqnarray*}
By Gronwall's inequality,
\begin{equation*}
\phi(t) \leq (C_0\,C)^2\,\sqrt{t}\,\exp{((C_0\,C)^2\,t)}
\end{equation*}
for all $t\in [0,T]$ and $N\geq N_0$. Hence the first inequality in Lemma \ref{L:BoundsForGGG} are established. The remaining inequalities in the lemma then follow by the same argument, using point-wise upper bound for $|G|$ and $|G^{\pm}|$ we obtained.
\end{pf}

\begin{remark}\rm
It can be shown that $(G^N,\,G^{N,+},\,G^{N,-})$ converges uniformly on compact subsets of $[0,\infty)\times (\bar{D}_+\times \bar{D}_-\setminus \mathfrak{I})$, $[0,\infty)\times (\bar{D}_+\times \bar{D}_+\setminus I\times I)$ and $[0,\infty)\times (\bar{D}_-\times \bar{D}_-\setminus I\times I)$ respectively, where $\mathfrak{I}:= \{(z,z)\in \bar{D}_+\times \bar{D}_-:\;z\in I\}$. Furthermore, the limit  $(G^{\infty},\,G^{\infty,+},\,G^{\infty,-})$ is the unique continuous solution to the following couple integral equations.

\begin{eqnarray*}
G^{\infty}_t(x,y)&=&-\,\int_0^t \int_{I}\int_{D_-} p(t-r,(x,y),(z,\tilde{y}))\,(G^{\infty}_r(z,\tilde{y})f^-_r(z)+ G^{\infty,-}_r(\tilde{y},z)f^+_r(z))\,d\tilde{y}\,d\sigma(z)\\
&&\quad +\int_{I}\int_{D_+}p(t-r,(x,y),(\tilde{x},z))\,(G^{\infty}_r(\tilde{x},z)f^+_r(z)+ G^{\infty,+}_r(\tilde{x},z)f^-_r(z))\,d\tilde{x}\,d\sigma(z)\\
&&\quad -\int_{I}p(t-r,(x,y),(z,z))f^+_r(z)f^-_r(z)\,d\sigma(z)\;dr,
\end{eqnarray*}
\begin{eqnarray*}
G^{\infty,+}_t(x_1,x_2)&=&-\,\int_0^t \int_{I}\int_{D_+} p(t-r,(x_1,x_2),(z,\tilde{x}))\,(G^{\infty,+}_r(z,\tilde{x})f^-_r(z)+ G^{\infty}_r(\tilde{x},z)f^+_r(z))\\
&& \quad + \, p(t-r,(x_1,x_2),(\tilde{x},z))\,(G^{\infty,+}_r(\tilde{x},z)f^-_r(z)+ G^{\infty}_r(\tilde{x},z)f^+_r(z))\,d\tilde{x}\,d\sigma(z)\,dr,
\end{eqnarray*}
\begin{eqnarray*}
G^{\infty,-}_t(y_1,y_2)&=&-\,\int_0^t \int_{I}\int_{D_-} p(t-r,(y_1,y_2),(z,\tilde{y}))\,(G^{\infty,-}_r(z,\tilde{y})f^+_r(z)+ G^{\infty}_r(z,\tilde{y})f^-_r(z))\\
&& \quad +\,p(t-r,(y_1,y_2),(\tilde{y},z))\,(G^{\infty,-}_r(\tilde{y},z)f^+_r(z)+ G^{\infty}_r(z,\tilde{y})f^-_r(z))\,d\tilde{y}\,d\sigma(z)\,dr.
\end{eqnarray*}
\end{remark}

\subsubsection{Evolution operators $Q^N_{s,t}$ and $\mathbf{U}^N_{(t,s)}$}

We fix $N\in \mathbb{N}$ and consider the following coupled backward PDE for $(v^{+}_N,\,v^{-}_N)$, with Neumann boundary conditions and terminal data $v^{\pm}_{N}(t)=\phi_{\pm}\in L^2(D_{\pm})$:
    \begin{equation}\label{E:BackwardCoupleEqtQN}
        \left\{\begin{aligned}
        -\,\dfrac{\partial v^+_N}{\partial s}           &= \frac{1}{2}\Delta v^+_N -\<\ell(v^+_N+v^-_N),\,f^-\>_-    & &\qquad\text{on } (0,t)\times D_+  \\
        -\,\dfrac{\partial v^-_N}{\partial s}           &= \frac{1}{2}\Delta v^-_N -\<\ell(v^+_N+v^-_N),\,f^+\>_+    & &\qquad\text{on } (0,t)\times D_-,
        \end{aligned}\right.
\end{equation}
where $(f^{+},\,f^{-})=(f^{N,+},\,f^{N,-})$ is defined in (\ref{E:CoupledIntegral_fN}) and $\<\phi,\,\psi\>_{\pm}:= \int_{D_{\pm}}\phi\,\psi$.

Note that each of the two equations in (\ref{E:BackwardCoupleEqtQN}) is of the form $\frac{-\partial v}{\partial s} = \frac{1}{2}\Delta v -kv-h$ where $k(s,x)\geq 0$ is a killing potential and $h(s,x)$ (not necessarily non-negative) is an external perturbation. This is because we can rewrite
\begin{eqnarray}
\<\ell(v^++v^-),\,f^-\>_- &\text{ as }& k^+v^++h^+ := \<\ell,\,f^-\>_{-}\,v^+ + \<\ell v^-,\,f^-\>_- \label{Def:khPlus} \quad \text{and}\\
\<\ell(v^++v^-),\,f^+\>_+ &\text{ as }& k^-v^-+h^- := \<\ell,\,f^+\>_{+}\,v^- + \<\ell v^+,\,f^+\>_+\,. \label{Def:khMinus}
\end{eqnarray}

By the same proof as that of Proposition 2.19 in \cite{zqCwtF13a}, we have the following.
\begin{prop}\label{prop:Integral_BackwardCoupleEqtQN}
For $N\in\mathbb{N}$ large enough, $t>0$ and $\phi_{\pm}\in C(\bar{D}_{\pm})$. There is a unique element $(v^+,v^-)=(v^{N,+},\,v^{N,-})$ in $C([0,t]\times \bar{D}_+)\times C([0,t]\times \bar{D}_-)$ which satisfies the following coupled integral equations:
\begin{eqnarray*}
    v^+(s,x)&=&P^+_{t-s}\phi_+(x)-\dfrac{1}{2}\int_0^{t-s} P^+_{\theta}\left(k^+(s+\theta)v^+(s+\theta)+h^+(s+\theta)\right)(x)\,d\theta\\
    v^-(s,y)&=&P^-_{t-s}\phi_-(y)-\dfrac{1}{2}\int_0^{t-s} P^-_{\theta}\left(k^-(s+\theta)v^-(s+\theta)+h^-(s+\theta)\right)(y)\,d\theta,
\end{eqnarray*}
where $k^{\pm}$ and $h^{\pm}$ (which are functions indexed by $N$) are defined in (\ref{Def:khPlus}) and (\ref{Def:khMinus}). Moreover, $(v^+,v^-)$ has the following probabilistic representations:
\begin{eqnarray*}
v^+(s,x) &=&
\E\left[\phi_+(X^+_{t-s})e^{-\int_0^{t-s}k^+(s+r,X^+_r)dr}-\int_0^{t-s} h^+(s+\theta,X^+_{\theta})e^{-\int_0^{\theta}k^+(s+r,X^+_r)dr}\,d\theta\,\Big|\,X^+_0=x\right]\\
v^-(s,y) &=&
\E\left[\phi_-(X^-_{t-s})e^{-\int_0^{t-s}k^-(s+r,X^-_r)dr}-\int_0^{t-s} h^-(s+\theta,X^-_{\theta})e^{-\int_0^{\theta}k^-(s+r,X^-_r)dr}\,d\theta\,\Big|\,X^-_0=y\right].
\end{eqnarray*}
We call this $(v^+,v^-)=(v^{N,+},\,v^{N,-})$ the \textbf{probabilistic solution} of the coupled PDE (\ref{E:BackwardCoupleEqtQN}) with Neumann boundary conditions and terminal data $\phi_{\pm}$.
\end{prop}

\begin{definition}\label{Def:U^N_ts}
For $0 \leq s\leq t$ and $\phi_{\pm}\in C(\bar{D}_{\pm})$, we define
$$Q^N_{s,t}(\phi_+,\phi_-):=  (v^{N,+}(s),\,v^{N,-}(s))$$
to be the probabilistic solution given by Proposition \ref{prop:Integral_BackwardCoupleEqtQN}. Clearly, $Q^N_{s,t}:\,C(\bar{D}_+)\times C(\bar{D}_-)\rightarrow C(\bar{D}_+)\times C(\bar{D}_-)$ and $Q^N_{s,u}\circ \,Q^N_{u,t}= Q^N_{s,t}$ for $0\leq s\leq u\leq t$. Now we define
\begin{equation}
\<\mathbf{U}^N_{(t,s)}\mu,\; (\phi_+,\phi_-)\> := \<\mu,\;Q^N_{s,t}(\phi_+,\phi_-)\>
\end{equation}
for $\alpha\geq 0$, $\mu\in \mathbf{H}_{-\alpha}$ and $(\phi_+,\phi_-)$ whenever $Q^N_{s,t}(\phi_+,\phi_-) \in \H_{\alpha}^+\times \H_{\alpha}^-$.
\end{definition}

\subsubsection{Evolution operators $Q_{s,t}$ and $\mathbf{U}_{(t,s)}$}

Formally, if we let $N\to\infty$ in (\ref{E:BackwardCoupleEqtQN}), we obtain
    \begin{equation}\label{E:BackwardCoupleEqtQ}
        \left\{\begin{aligned}
        -\,\dfrac{\partial v^+}{\partial s}           &= \frac{1}{2}\Delta v^+ -(v^++v^-)\,u_-\,d\sigma\big|_{I}    & &\qquad\text{on } (0,t)\times D_+  \\
        -\,\dfrac{\partial v^-}{\partial s}           &= \frac{1}{2}\Delta v^- -(v^++v^-)\,u_+\,d\sigma\big|_{I}   & &\qquad\text{on } (0,t)\times D_-,
        \end{aligned}\right.
    \end{equation}
where $(u_+,\,u_{-})$ is the hydrodynamic limit
of the interacting diffusion systems.
 Observe that this equation is equivalent to (\ref{E:FluctuationPDE_limit}) with $\lambda=1$ and $\A^{\pm}=\frac{1}{2}\Delta$. Note the difference between this coupled PDEs and that for the hydrodynamic limit.

Recall that $k^{\pm}$, $h^{\pm}$, $f^{\pm}$, and $v^{\pm}$ are functions indexed by $N$.
Heuristically, as $N\to \infty$,  we have
\begin{eqnarray*}
&& P^+_{\theta}\left(k^+(s+\theta)v^{+}(s+\theta)+h^+(s+\theta)\right)(x)  \\
&=& \int_{D_+}p^+(\theta,x,z)\left(\<\ell,\,f^-\>_{-}(s+\theta,z)v^{+}(s+\theta,z)+ \<\ell v^{-},\,f^-\>_-(s+\theta,z)\right)\,dz\\
&\rightarrow& \int_{I}p^+(\theta,x,z)\,\big[v^+(s+\theta,z)+v^-(s+\theta,z)\big]\,u_-(s+\theta,z)\,d\sigma(z) \\
&:=& {\cal G}^+_{\theta}\Big(\,[v^+(s+\theta)+v^-(s+\theta)]\,u_-(s+\theta)\,\Big)(x).
\end{eqnarray*}
The abbreviation in the last term is based on the notation ${\cal G}^{\pm}_{\theta}\varphi(x):= \int_I p^{\pm}(\theta,x,z)\,\varphi(z)\,d\sigma(z)$.
The following result is analogous to Proposition \ref{prop:Integral_BackwardCoupleEqtQN} and can be proved as in the same way.

\begin{prop}\label{prop:Integral_BackwardCoupleEqtQ}
Fix $t>0$ and $\phi_{\pm}\in C(\bar{D}_{\pm})$. There is a unique element $(v^+,v^-)$ in $C([0,t]\times \bar{D}_+)\times C([0,t]\times \bar{D}_-)$ which satisfies the following coupled integral equations:
\begin{eqnarray*}
    v^+(s,x)&=&P^+_{t-s}\phi_+(x)-\dfrac{1}{2}\int_0^{t-s} {\cal G}^+_{\theta}\Big(\,[v^+(s+\theta)+v^-(s+\theta)]\,u_-(s+\theta)\,\Big)(x)\,d\theta\\
    v^-(s,y)&=&P^-_{t-s}\phi_-(y)-\dfrac{1}{2}\int_0^{t-s} {\cal G}^-_{\theta}\Big(\,[v^+(s+\theta)+v^-(s+\theta)]\,u_+(s+\theta)\,\Big)(y)\,d\theta,
\end{eqnarray*}
where ${\cal G}^{\pm}_{\theta}\varphi(x):= \int_I p^{\pm}(\theta,x,z)\,\varphi(z)\,d\sigma(z)$. Moreover, $(v^+,v^-)$ has the following probabilistic representations:
\begin{eqnarray*}
v^+(s,x) &=&
\E\Big[\phi_+(X^+_{t-s})e^{-\int_0^{t-s}u_-(s+r,X^+_r)dL^+_r}\\
&& \qquad  -\int_0^{t-s} (v^-\cdot u_-)(s+\theta,X^+_{\theta})\,e^{-\int_0^{\theta}u_-(s+r,X^+_r)dL^+_r}\,dL^+_{\theta}\,\Big|\,X^+_0=x\Big]\\
v^-(s,y) &=&
\E\Big[\phi_-(X^-_{t-s})e^{-\int_0^{t-s}u_+(s+r,X^-_r)dL^-_r}\\
&& \qquad  -\int_0^{t-s} (v^+ \cdot u_+)(s+\theta,X^-_{\theta})\,e^{-\int_0^{\theta}u_+(s+r,X^-_r)dL^-_r}\,dL^-_{\theta}\,\Big|\,X^-_0=y\Big],
\end{eqnarray*}
where $L^{\pm}_t$ is the boundary local time of the RBM $X^{\pm}$ on $I$. We call this $(v^+,v^-)$ the \textbf{probabilistic solution} of the coupled PDE (\ref{E:BackwardCoupleEqtQ}) with terminal data $\phi_{\pm}$.
\end{prop}
We stress that the
right hand side of the above formula is well-defined; for instance,
$\int_0^{t-s}u_-(s+r,X^+_r)dL^+_r$ is well-defined since the value of $u_-$ at $(s+r,X^+_r)$ is picked up only when $X^+$ hits $I$ (which is a subset of $\bar{D}_-$).

\begin{definition}\label{Def:Q_ts_U_ts}
For $0 \leq s\leq t$ and $\phi_{\pm}\in C(\bar{D}_{\pm})$, we define
$$Q_{s,t}(\phi_+,\phi_-):=  (v^{+}(s),\,v^{-}(s))$$
to be the probabilistic solution given by Proposition \ref{prop:Integral_BackwardCoupleEqtQ}. Clearly, $Q_{s,t}:\,C(\bar{D}_+)\times C(\bar{D}_-)\rightarrow C(\bar{D}_+)\times C(\bar{D}_-)$ and $Q_{s,u}\circ \,Q_{u,t}= Q_{s,t}$ for $0\leq s\leq u\leq t$. To stress the dependence in $t$, we sometimes write $Q_{s,t}(\phi_+,\phi_-)$ as $(v^{+}_t(s),\,v^{-}_t(s))$ for $(\phi_+,\phi_-)$ fixed. Now we define, for $\alpha>0$ and $\mu\in \H^+_{-\alpha}\oplus \H^-_{-\alpha}$,
\begin{equation}
\<\mathbf{U}_{(t,s)}\mu,\; (\phi_+,\phi_-)\> := \<\mu,\;Q_{s,t}(\phi_+,\phi_-)\>.
\end{equation}
\end{definition}

\subsubsection{Key estimates for evolution operators}

On the space $C(\bar{D}_+)\times C(\bar{D}_-)$, we let $(\psi_+,\psi_-)-(\phi_+,\phi_-)= (\psi_+-\phi_+,\,\psi_--\phi_-)$ and denote by $\|(\psi_+,\psi_-)\|:= \|\psi_+\| + \|\psi_-\|$ the sum of the sup-norm of its components. The following uniform bound and uniform convergence are useful in many places in this
paper.

\begin{lem}
For all $\phi_{\pm}\in C(\bar{D}_{\pm})$ and $0\leq s\leq t \leq T$, we have
\begin{equation}\label{E:ContractionQNQ_AnnihilatingSystem}
\sup_{0\leq s\leq t\leq T} \left(\sup_{N\geq N_0}\|Q^N_{s,t}(\phi_+,\phi_-)\|\vee \|Q_{s,t}(\phi_+,\phi_-)\|\right) \leq \widehat{c}\,\|(\phi_+,\phi_-)\|
\end{equation}
for some positive integer $N_0=N_0(D_+,D_-)$ and $\widehat{c}=\widehat{c}(d,D_+,D_-,T)>0$. Moreover,
\begin{equation}\label{E:UConvergenceQNQ_AnnihilatingSystem}
\lim_{N\to\infty} \Big\|Q^N_{s,t}(\phi_+,\phi_-)- Q_{s,t}(\phi_+,\phi_-)\Big\|=0.
\end{equation}
\end{lem}

\begin{pf}
Recalling the probabilistic representations of $Q^N$ and $Q$ in Proposition \ref{prop:Integral_BackwardCoupleEqtQN} and Proposition \ref{prop:Integral_BackwardCoupleEqtQ} respectively, we see that (\ref{E:ContractionQNQ_AnnihilatingSystem}) follows from the non-negativity of $f^{N,\pm}$ and $u_{\pm}$. To prove (\ref{E:UConvergenceQNQ_AnnihilatingSystem}), we fix $t<T$ and let $(v^{N,+}(s),\,v^{N,-}):=Q^N_{s,t}(\phi_+,\phi_-)$ and  $(v^{+}(s),\,v^{-}):= Q_{s,t}(\phi_+,\phi_-)$ for $s\in[0,t]$. We look at the RHSs of the integral equations satisfied by $v^{N,+}(s)$ and $v^{+}(s)$, in Proposition \ref{prop:Integral_BackwardCoupleEqtQN} and Proposition \ref{prop:Integral_BackwardCoupleEqtQ}, respectively. The proof is the same as that of Lemma \ref{L:UConvergencefNtoU}, with the uniform bound for $f^{N,\pm}$ replaced by the bound (\ref{E:ContractionQNQ_AnnihilatingSystem}).
\end{pf}

\begin{lem}\label{L:Bound_QtMinusQs_AnnihilatingSystem}
There exists a constant $\bar{c}>0$ such that for any $0\leq s\leq t\leq T$ and $k\in \mathbb{N}$, we have
\begin{eqnarray*}
&& \sup_{r\in[0,s]}\big\|Q_{(r,t)}(\phi^+_{k},0)\,-\,Q_{(r,s)}(\phi^+_{k},0)\big\| \notag\\
&\leq& (\widehat{c}\vee 1)\,\bar{c}\,\|\phi^+_{k}\|
\left(\,\lambda^+_{k}\,(t-s)
+\widehat{c}\,(c^++c^-)\,(\|u^+_0\|\vee\|u^-_0\|)\,(t-s)^{1/2}\,\right)
\end{eqnarray*}
and
\begin{eqnarray*}
&& \sup_{r\in[0,s]}\big\|Q_{(r,t)}(0,\phi^-_{k})\,-\,Q_{(r,s)}(0,\phi^-_{k})\big\| \notag\\
&\leq& (\widehat{c}\vee 1)\,\bar{c}\,\|\phi^-_{k}\|
\left(\,\lambda^-_{k}\,(t-s)
+\widehat{c}\,(c^++c^-)\,(\|u^+_0\|\vee\|u^-_0\|)\,(t-s)^{1/2}\,\right).
\end{eqnarray*}
Here $c^{\pm}=C(d,D_{\pm},T)$ is the constant in (\ref{E:Surface_integral_boundedness}) applied to $D_{\pm}$ and $\widehat{c}=\widehat{c}(d,D_+,D_-,T)$ is the constant in (\ref{E:ContractionQNQ_AnnihilatingSystem}). Furthermore, there exists $N_0=N_0(D_+,D_-)$ such that for $N\geq N_0$, the above two inequalities remain valid with $\{Q^N_{s,t}\}$ in replace of $\{Q_{s,t}\}$ and $c^{\pm}=C(d,D_{\pm},T)$ being the constant in (\ref{E:boundary_strip_boundedness}).
\end{lem}

\begin{pf}
We fix $(\phi^+_{k},0)$ and only prove the first inequality, since the second inequality follows from the same argument.
 Recall the definition of $Q_{r,t}(\phi^+_{k},0)\in C(\bar{D}_+)\times C(\bar{D}_+)$ in Definition \ref{Def:Q_ts_U_ts}. Suppose $Q_{r,t}(\phi^+_{k},0)=(v_t^+(r),\,v_t^-(r))$. Then
\begin{eqnarray*}
    v^+_t(r,x)&=&P^+_{t-s}\phi^+_{k}(x)-\dfrac{1}{2}\int_0^{t-r} {\cal G}^+_{\theta}\Big(\,[v^+_t(r+\theta)+v^-_t(r+\theta)]\,u_-(r+\theta)\,\Big)(x)\,d\theta , \\
    v^-_t(r,y)&=&0 -\dfrac{1}{2}\int_0^{t-r} {\cal G}^-_{\theta}\Big(\,[v^+_t(r+\theta)+v^-_t(r+\theta)]\,u_+(r+\theta)\,\Big)(y)\,d\theta,
\end{eqnarray*}
where ${\cal G}^{\pm}_{\theta}\varphi(x):= \int_I p^{\pm}(\theta,x,z)\,\varphi(z)\,d\sigma(z)$.
Hence for any $0\leq r\leq s\leq t$, we have
\begin{eqnarray*}
&& |v^+_t(r,x)-v^+_s(r,x)| \\
&\leq & \big|\left(e^{-\lambda^+_{k}(t-r)}- e^{-\lambda^+_{k}(s-r)} \right)\,\phi^+_{k}(x)\big| \\
&& + \dfrac{1}{2}\Big|\int_{s-r}^{t-r}{\cal G}^+_{\theta}\Big(\,[v^+_t(r+\theta)+v^-_t(r+\theta)]\,u_-(r+\theta)\,\Big)(x)\,d\theta\Big|\\
&&+ \dfrac{1}{2}\Big|\int_{0}^{s-r}
{\cal G}^+_{\theta}\Big(\,[(v^+_t-v^+_s)(r+\theta)+(v^-_t-v^-_s)(r+\theta)]\,u_-(r+\theta)\,\Big)(x)
\,d\theta\Big|.\\
\end{eqnarray*}
A similar inequality holds for $|v^-_t(r,y)-v^-_s(r,y)|$, which has 2 terms instead of 3 terms on the RHS. View $0\leq s\leq t$ as fixed and define, for $r\in[0,s]$,
\begin{equation*}
f(r) := \big\|Q_{(r,t)}(\phi^+_{k},0)\,-\,Q_{(r,s)}(\phi^+_{k},0)\big\|=\|(v^+_t-v^+_s)(r)\|+\|(v^-_t-v^-_s)(r)\|.
\end{equation*}
Then the above estimates, together with (\ref{E:Surface_integral_boundedness}) and (\ref{E:ContractionQNQ_AnnihilatingSystem}), implies that
\begin{equation}\label{E:Bound_QtMinusQs_AnnihilatingSystem}
f(r) \leq A + B \int_0^{s-r}\dfrac{f(r+\theta)}{\sqrt{\theta}}\,d\theta \quad \text{for } r\in[0,s],
\end{equation}
where $B=(\|u^+_0\|\vee\|u^-_0\|)\,(c^++c^-)$ and $$A=\lambda^+_{k}\,\|\phi^+_{k}\|\,(t-s)+\widehat{c}\,(c^++c^-)\,(\|u^+_0\|\vee\|u^-_0\|)\,\|\phi^+_{k}\|\,(t-s)^{1/2}.$$
Iterating (\ref{E:Bound_QtMinusQs_AnnihilatingSystem}), we have
\begin{eqnarray*}
f(r) &\leq&  A+ AB\int_{w_1=r}^s\frac{1}{\sqrt{w_1-r}} + AB^2\int_{w_1=r}^s\int_{w_2=w_1}^s\frac{1}{\sqrt{(w_1-r)(w_2-w_1)}} \\
&& \; +AB^3\int_{w_1=r}^s\int_{w_2=w_1}^s\int_{w_3=w_2}^s\frac{1}{\sqrt{(w_1-r)(w_2-w_1)(w_3-w_2)}}+ \cdots  \\
&=& A\,\sum_{k=0}^{\infty} B^k\,a_k\,(s-r)^{k/2}, \quad \text{where}\quad a_k=\frac{\pi^{k/2}}{\Gamma((k+2)/2)}\\
&\leq& \frac{\bar{c}}{2}\,A\,\sum_{k=0}^{\infty} B^k\,(s-r)^{k/2} \quad \text{for some absolute constant }\bar{c}>0\\
&\leq& \bar{c}\,A \qquad \text{if}\quad |B\sqrt{s-r}|\leq 1/2.
\end{eqnarray*}
Hence,
\begin{equation}\label{E:Bound_QtMinusQs_AnnihilatingSystem2}
\sup_{r\in \left[0\vee s-1/(4B^2),\,s\right]}f(r) \,\leq\,\bar{c}\,A.
\end{equation}
(The case $B=0$ is trivial.)
We can then extend (\ref{E:Bound_QtMinusQs_AnnihilatingSystem2}) to take care of the case $0\leq r<s-1/(4B^2)$. Namely, by the evolution property and (\ref{E:ContractionQNQ_AnnihilatingSystem}), we have
\begin{eqnarray*}
&& \big\|Q_{(r,t)}(\phi^+_{k},0) \,-\, Q_{(r,s)}(\phi^+_{k},0)\big\|\\
&=& \big\|Q_{(r,\,s-1/(4B^2))}\left(Q_{(s-1/(4B^2),\,t)}(\phi^+_{k},0) \,-\, Q_{(s-1/(4B^2),\,s)}(\phi^+_{k},0)\right)\big\|\\
&\leq& \widehat{c}\,\big\| Q_{(s-1/(4B^2),\,t)}(\phi^+_{k},0)\,-\,Q_{(s-1/(4B^2),\,s)}(\phi^+_{k},0) \big\| \leq \hat{c}\,\bar{c}\,A.
\end{eqnarray*}
The proof is complete.
\end{pf}

Due to the annihilation between two kinds of particles, unlike the case considered in
\cite{zqCwtF14c}, we need to analyze the correlation functions more deeply. This will be developed in the next two sections.

\subsection{Asymptotic expansion for correlation functions $F^{N,(n,m)}_t$}

\begin{definition}\label{Def:CorrelationFcn}
    Fix $N\in \mathbb{N}$ and consider the annihilating diffusion system. For $n,m\in \mathbb{N}$ and $t\geq 0$, we define the $(n,m)$\textbf{-correlation function at time} $t$, $F^{(n,m)}_t=F^{N,(n,m)}_t$, by
    \begin{equation*}
    \int_{D_+^n\times D_-^m}\Phi(\vec{x},\vec{y})\,F^{(n,m)}_t(\vec{x},\vec{y})\,d(\vec{x},\vec{y})= \E\left[\Phi_{(n,m)}(t)\right]
    \end{equation*}
    for all $\Phi\in C(\bar{D}_+^n\times \bar{D}_-^m)$, where
    \begin{equation}\label{E:Def_Phi_nm}
    \Phi_{(n,m)}(t):= \dfrac{1}{N^{(n)}\,N^{(m)}}\sum_{\substack{i_1,\cdots i_n\\ \text{distinct}}}^{\sharp_t}\,\sum_{\substack{j_1,\cdots j_m\\ \text{distinct}}}^{\sharp_t}\,\Phi(X^{i_1}_t,\cdots, X^{i_n}_t,\,Y^{j_1}_t,\cdots,Y^{j_m}_t),
    \end{equation}
    $\sharp_t$ is the number of particles alive at time $t\in[0,\infty)$ in each of $\bar{D}_{\pm}$ and $N^{(n)}:= N(N-1)\cdots (N-n+1)$ is the number of permutations of $n$ objects chosen from $N$ objects
with $N^{(0)} := 1$.
\end{definition}

\begin{example}
For example, we have
\begin{equation*}
\E[\<\phi,\X^{N,+}_t\>]=\int_{D_+}\phi(x)\,F^{(1,0)}_{t}(x) \quad\text{and}\quad
\E[\<\ell,\,\X^{N,+}_t\otimes \X^{N,-}_t\>]=\int_{D_+\times D_-} \ell(x,y)\,F^{(1,1)}_{t}(x,y).
\end{equation*}
\end{example}

Intuitively, if we randomly pick $n$ and $m$
indistinguishable
living particles in $D_+$ and $D_-$ respectively at time $t$, then $F^{(n,m)}_t(\vec{x},\vec{y})$ is the probability joint density function for their positions. Note that $F^{(n,m)}_t$ is defined for almost all $(\vec{x},\vec{y})\in D_+^n\times D_-^m$, and that it depends on both $N$ and the initial configurations $(\X^{N,+}_0,\,\X^{N,-}_0)$. We will see, via the BBGKY hierarchy (\ref{E:BBGKY_nm_correlation_t}) which will be proved below, that $F^{(n,m)}_t\in C(\bar{D}_+^n\times \bar{D}_-^m)$ for $t>0$. We can also replace $N^{(n)}$ by $N^n$ (cf. Dittrich \cite{pD88b} and Lang and Xanh \cite{LX80}).
This  is because we are interested in the behavior of $F^{n,m}$ as $N\to\infty$, and for each fixed $n$,
$$
\frac{N^{(n)}}{N^n}=(1-\frac{1}{N})(1-\frac{2}{N})\cdots (1-\frac{n-1}{N})\nearrow 1 \text{  as }N\to\infty.
$$

It is natural, base on the annihilating random walk model in \cite{zqCwtF13a}, to expect that we have \textbf{propagation of chaos}, which says that when the number of particles tends to infinity, the particles will appear to be independent from each other. More precisely, we expect to have
\begin{equation}
    \lim_{N\to\infty}F^{(n,m)}_t(\vec{x},\vec{y})= \prod_{i=1}^n u_+(t,x_i)\,\prod_{j=1}^m u_-(t,y_j).
\end{equation}
This will be implied by a more exact asymptotic behavior of the $F^{(n,m)}$, namely Theorem \ref{T:Asymptotic_nm_correlation_t}, which is a key ingredient for the study of fluctuation.
Our method is motivated by the approach of \cite{pD88b}.

\begin{thm}\label{T:Asymptotic_nm_correlation_t}
Suppose that $F_0^{(n,m)}(\vec{x},\vec{y})=\prod_{i=1}^nu^+_0(x_i)\,\prod_{j=1}^mu^-_0(y_j)$ 
(this implies that the $N$ particles are initially independently distributed) 
and that $C_0:= \|u^{+}_0\|\vee \|u^{-}_0\|$. Then for all $T>0$, there exists $C=C(D_+,D_-,T)>0$ and an integer $N_0=N_0(D_+,D_-)$ such that for $0\leq t \leq T\wedge(C_0\,C)^{-2}$ and $N\geq N_0$, the correlation function has decomposition
\begin{equation}
F^{(n,m)}_t=A_t^{(n,m)}+\frac{1}{N}\,B_t^{(n,m)}+\frac{1}{N}\,C_t^{(n,m)}
\end{equation}
with
\begin{equation}\label{E:Asymptotic_nm_correlation_t}
\|A^{(n,m)}_t-F^{(n,m)}_t\|_{(n,m)} \leq \frac{(C_0\,C)^{n+m}\,\sqrt{t}}{N\,\delta_N^d}\, \quad\text{and }\quad
\|C^{(n,m)}_t\|_{(n,m)} \leq \frac{(C_0\,C)^{n+m}\,t}{N\,\delta_N^{2d}},
\end{equation}
where
\begin{eqnarray*}
A_t^{(n,m)}(\vec{x},\vec{y})&:=&\prod_{i=1}^nf^+_t(x_i)\,\prod_{j=1}^mf^-_t(y_j)\,,\\
B_t^{(n,m)}(\vec{x},\vec{y})&:=& -\,A_t^{(n,m)}(\vec{x},\vec{y})\,\bigg(
\sum_{i=1}^n\frac{g^+_t(x_i)}{f^+_t(x_i)}+\sum_{j=1}^m\frac{g^-_t(y_j)}{f^-_t(y_j)}+\sum_{i=1}^n\sum_{j=1}^m\frac{G_t(x_i,y_j)}{f^+_t(x_i)f^-_t(y_j)}\\
&& \qquad\qquad\qquad
+\sum_{i<p}^{n}\frac{G^+_t(x_i,x_p)}{f^+_t(x_i)f^+_t(x_p)}+\sum_{j<q}^{m}\frac{G^-_t(y_j,y_q)}{f^-_t(y_j)f^-_t(y_q)}
\bigg).
\end{eqnarray*}
\end{thm}

\begin{pf}
The key point of our method is to compare three hierarchies (\ref{E:BBGKY_nm_correlation_t},) (\ref{E:BBGKY_nm_A_t}) and (\ref{E:BBGKY_nm_alpha_t}) in Step 1 below:

\textbf{Step 1: BBGKY hierarchy for the correlation functions. }

Apply Dynkin's formula to (see \cite[Corollary 7.8]{zqCwtF14b}) the functional
$$(s,\,(\X^{N,+}_s,\X^{N,-}_s)) \mapsto
\dfrac{1}{N^{(n)}\,N^{(m)}}
\sum_{\substack{\text{distinct} \\ i_1,\cdots, i_n=1}}^{\sharp_s}
\,\sum_{\substack{\text{distinct} \\j_1,\cdots, j_m=1 }}^{\sharp_s}
\,P_{t-s}\Phi(X^{i_1}_s,\cdots, X^{i_n}_s,\,Y^{j_1}_s,\cdots,Y^{j_m}_s)
\;\;,\,s\in[0,t]$$
yields
\begin{equation}\label{E:BBGKY_nm_correlation_t}
F^{(n,m)}_t=P^{(n,m)}_t F^{(n,m)}_0 - \int_0^t P^{(n,m)}_{t-s}\left( V^+F^{(n,m+1)}_s + V^-F^{(n+1,m)}_s + \frac{Q}{N}F^{(n,m)}_s \right)\;ds,
\end{equation}
where $V^+=\sum_{i=1}^nV_{+_i}$,  $V^-=\sum_{j=1}^mV_{-_j}$ are operators, $V_{+_i}F^{(n,m+1)}$, $V_{-_j}F^{(n+1,m)}$ and $QF^{(n,m)}$ are functions on $\bar{D}_+^n \times \bar{D}_-^m$ defined by
\begin{eqnarray*}
V_{+_i}F^{(n,m+1)}(\vec{x},\vec{y})&:=& \int_{D_-}\ell(x_i,y)\,F^{(n,m+1)}(\vec{x},\,(\vec{y},y))\,dy \\
V_{-_j}F^{(n+1,m)}(\vec{x},\vec{y})&:=& \int_{D_+}\ell(x,y_j)\,F^{(n+1,m)}((\vec{x},x),\,\vec{y})\,dx \\
QF^{(n,m)}(\vec{x},\vec{y})&:=& \left(\sum_{i=1}^n\sum_{j=1}^m \ell(x_i,y_j)\right)\,F^{(n,m)}(\vec{x},\vec{y}).
\end{eqnarray*}
Note that $Q$ is a multiplication operator, so it is natural to denote $Q^{(n,m)}$ to be the function $Q^{(n,m)}(\vec{x},\vec{y})=\sum_{i=1}^n\sum_{j=1}^m \ell(x_i,y_j)$. Note also that the above is a finite sum since $F^{(n,m)}=0$ when $n\vee m>N$. The system of equation (\ref{E:BBGKY_nm_correlation_t}) is usually called \emph{BBGKY hierarchy}.
\footnote{We can also view (\ref{E:BBGKY_nm_correlation_t}) as the `variation of constant' and
$F^{(n,m)}$ as the probabilistic solution (cf. \cite[Proposition 2.19]{zqCwtF13a}) for the following heat equation on $D_+^n\times D_-^m$ with Neumann boundary condition:
\begin{equation}
\dfrac{\partial F^{(n,m)}_t}{\partial t}= \frac{1}{2}\,\Delta F^{(n,m)}_t-
\left( V^+F^{(n,m+1)}_t + V^-F^{(n+1,m)}_t + \frac{Q}{N}F^{(n,m)}_t \right).
\end{equation}
}

On the other hand, it can be easily verified that $A^{(n,m)}_t$ solves
\begin{equation}\label{E:BBGKY_nm_A_t}
A^{(n,m)}_t=P^{(n,m)}_t A^{(n,m)}_0 - \int_0^t P^{(n,m)}_{t-s}\left( V^+A^{(n,m+1)}_s + V^-A^{(n+1,m)}_s \right)\;ds,
\end{equation}
and that we have chosen $B^{(n,m)}_t$ in such a way that $\alpha_t^{(n,m)}:= A^{(n,m)}_t+\dfrac{B^{(n,m)}_t}{N}$ solves
\begin{equation}\label{E:BBGKY_nm_alpha_t}
\alpha^{(n,m)}_t=P^{(n,m)}_t F^{(n,m)}_0 - \int_0^t P^{(n,m)}_{t-s}\left( V^+\alpha^{(n,m+1)}_s + V^-\alpha^{(n+1,m)}_s + \frac{Q}{N}A^{(n,m)}_s \right)\;ds.
\end{equation}

\textbf{Step 2: Duhamel expansion for $N(A^{(n,m)}_t-F^{(n,m)}_t)$ in terms of a tree. }

Since $F^{(n,m)}_0=A^{(n,m)}_0$ by assumption, by repeatedly iterating (\ref{E:BBGKY_nm_correlation_t}) and (\ref{E:BBGKY_nm_A_t}), we have
\begin{eqnarray*}
&& N(A^{(n,m)}_t-F^{(n,m)}_t)\\
&=& \;\;-\int_{t_2=0}^t P^{(n,m)}_{t-t_2}QF^{(n,m)}_{t_2}\\
&&+ \int_{t_2=0}^t\int_{t_3=0}^{t_2} P^{(n,m)}_{t-t_2}\left(\sum_{i=1}^nV_{+_i}P^{(n,m+1)}_{t_2-t_3}QF^{(n,m+1)}_{t_3} +
\sum_{j=1}^mV_{-_j}P^{(n+1,m)}_{t_2-t_3}QF^{(n+1,m)}_{t_3}
\right)\\
&&-\cdots\\
&&+ (-1)^{M}
\int_{t_2=0}^t\int_{t_3=0}^{t_2}\cdots \int_{t_{M+1}=0}^{t_M}   \sum_{\vec{\theta}\in \,\mathbb{T}^{(n,m)}_{M-1}}\,
P^{(n,m)}_{t-t_2}\,V_{\theta_1}\,P^{l_1(\vec{\theta})}_{t_2-t_3}\,V_{\theta_2}\,P^{l_2(\vec{\theta})}_{t_3-t_4}\,V_{\theta_3}\cdots P^{l_{M-1}(\vec{\theta})}_{t_M-t_{M+1}}\,Q\,F^{l_{M-1}(\vec{\theta})}_{t_{M+1}}\\
&&+\cdots,
\end{eqnarray*}
where $\mathbb{T}^{(n,m)}_{M-1}$ and $(l_1(\vec{\theta}),\,l_2(\vec{\theta}),\cdots,\,l_N(\vec{\theta}))$ are the tree and the labels defined in Subsection 3.5 in \cite{zqCwtF13a}.

Replacing $F^{l_{M-1}(\vec{\theta})}_{t_{M+1}}$ by the constant function $1$ in the $M-$th iterated integral above, we define the following function on $\bar{D}_+^n \times \bar{D}_-^m$:
\begin{equation}
\Theta^{(n,m)}_M(t) :=
\int_{t_2=0}^t\int_{t_3=0}^{t_2}\cdots \int_{t_{M+1}=0}^{t_M} \sum_{\vec{\theta}\in \,\mathbb{T}^{(n,m)}_{M-1}}\,
P^{(n,m)}_{t-t_2}\,V_{\theta_1}\,P^{l_1(\vec{\theta})}_{t_2-t_3}\,V_{\theta_2}\,P^{l_2(\vec{\theta})}_{t_3-t_4}\,V_{\theta_3}\cdots P^{l_{M-1}(\vec{\theta})}_{t_M-t_{M+1}}\,Q\,1.
\end{equation}

\textbf{Step 3: Bounding $\|\Theta^{(n,m)}_M(t)\|_{(n,m)}$. }

We now bound $\|\Theta^{(n,m)}_M(t)\|_{(n,m)}$ by employing our method developed in Subsection 3.5 in \cite{zqCwtF13a}. For the convenience of the reader, we summarize the key steps.

Note that $\Theta^{(n,m)}_M(t)$ is a sum of $(n+m)(n+m+1)\cdots (n+m+M-2)$ terms of multiple integrals. Following Subsection 3.5, we simplify (or telescope) each integrand by Chapman-Kolmogorov equation, and then apply \eqref{E:Surface_integral_boundedness} to obtain
\begin{eqnarray*}
\|\Theta^{(n,m)}_M(t)\|_{(n,m)} &\leq & \frac{1}{\delta_N^d}\,C^M\,(n+m+M-1) \\
&& \quad \int_{t_2=0}^{t}\cdots \int_{t_{M+1}=0}^{t_M}
\sum_{\vec{\upsilon}\,\in \mathbb{S}^{(n,m)}_M}\frac{1}{\sqrt{(t_{\upsilon_1}-t_2)\,(t_{\upsilon_2}-t_3)\cdots(t_{\upsilon_M}-t_{M+1})}}
\end{eqnarray*}
for all $0\leq t \leq T$ and $N\geq N_0(D_+,D_-)$, where $C=C(D_+,D_-,T)>0$ and $\mathbb{S}^{(n,m)}_M$ is a relabeled tree of $\mathbb{T}^{(n,m)}_M$ defined in Subsection 3.5 of \cite{zqCwtF13a}. Lemma 3.9 and Lemma 3.10 in \cite{zqCwtF13a}
give
\begin{eqnarray*}
&& \int_{t_2=0}^{t}\cdots \int_{t_{M+1}=0}^{t_M} \sum_{\vec{\upsilon}\,\in \mathbb{S}^{(n,m)}_M}\frac{1}{\sqrt{(t_{\upsilon_1}-t_2)\,(t_{\upsilon_2}-t_3)\cdots(t_{\upsilon_M}-t_{M+1})}}\\
&\leq& \frac{(n+m)^{(n+m)}}{(n+m)!}\,2^M\int_{t_2=0}^{t}\cdots \int_{t_{M+1}=0}^{t_M} \prod_{i=2}^{M+1}\bigg(\sum_{j=1}^{i-1}\frac{1}{\sqrt{t_j-t_i}}\bigg) \\
&\leq& c^{n+m+M}\,t^{M/2},
\end{eqnarray*}
where $c$ is an absolute constant. Therefore we have
\begin{eqnarray}\label{E:Bound_I_M_mn}
\|\Theta^{(n,m)}_M(t)\|_{(n,m)}\leq \frac{1}{\delta_N^d}\,C^{n+m+M}\,t^{M/2}
\end{eqnarray}
for all $0\leq t \leq T$ and $N\geq N_0(D_+,D_-)$, where $C=C(D_+,D_-,T)>0$.

\textbf{Step 4: Upper bound for $N\,\|A^{(n,m)}_t-F^{(n,m)}_t\|_{(n,m)} $. }

Since $\|F^{(p,q)}\|_{(p,q)}\leq C_0^{p+q}$, and since the sum of the two components in $l_{M-1}(\vec{\theta})$ is $n+m+M-1$ for any $\vec{\theta}\in \mathbb{T}^{(n,m)}_{M-1}$, Step 2 and Step 3 yields
\begin{eqnarray}\label{E:Bound_A_minus_F}
N\,\|A^{(n,m)}_t-F^{(n,m)}_t\|_{(n,m)} &\leq& \sum_{M=1}^{\infty}  C_0^{n+m+M-1}\,\|\Theta^{(n,m)}_M(t)\|_{(n,m)} \notag\\
&\leq& \frac{1}{\delta_N^d}\,\frac{C_0^{n+m}\,C^{n+m+1}\sqrt{t}}{1-C_0\,C\sqrt{t}} \notag\\
&\leq& \frac{1}{\delta_N^d}\,(C_0\,C)^{n+m}\,\sqrt{t}
\end{eqnarray}
for all $t\in[0,\,T\wedge(C_0\,C)^{-2})$ and $N\geq N_0$.

\textbf{Step 5: Upper bound for $C^{(n,m)} := \,N\,\left(F^{(n,m)}_t-A^{(n,m)}_t-\frac{B^{(n,m)}_t}{N}\right)$. }

Iterating $C^{(n,m)}_t := N\left(F^{(n,m)}_t-A^{(n,m)}_t-\dfrac{B^{(n,m)}_t}{N}\right)$ as in Step 2, we have
\begin{eqnarray*}
C^{(n,m)}_t
&=& \;\;\int_{t_2=0}^t P^{(n,m)}_{t-t_2}Q(A^{(n,m)}_{t_2}-F^{(n,m)}_{t_2})\\
&&+ \int_{t_2=0}^t\int_{t_3=0}^{t_2} P^{(n,m)}_{t-t_2}\bigg(\sum_{i=1}^nV_{+_i}P^{(n,m+1)}_{t_2-t_3}Q(A^{(n,m+1)}_{t_3}-F^{(n,m+1)}_{t_3}) \\
&& \qquad\qquad\qquad\qquad\qquad +
\sum_{j=1}^mV_{-_j}P^{(n+1,m)}_{t_2-t_3}Q(A^{(n+1,m)}_{t_3}-F^{(n+1,m)}_{t_3})
\bigg)
 +\cdots\\
&=& \sum_{M=1}^{\infty}
\int_{t_2=0}^t\int_{t_3=0}^{t_2}\cdots \int_{t_{M+1}=0}^{t_M} \\
&&\qquad \qquad \qquad \sum_{\vec{\theta}\in \,\mathbb{T}^{(n,m)}_{M-1}}\,
P^{(n,m)}_{t-t_2}\,V_{\theta_1}\,P^{l_1(\vec{\theta})}_{t_2-t_3}\,V_{\theta_2}\,P^{l_2(\vec{\theta})}_{t_3-t_4}\,V_{\theta_3}\cdots P^{l_{M-1}(\vec{\theta})}_{t_M-t_{M+1}}\,Q\,(A^{l_{M-1}(\vec{\theta})}_{t_{M+1}}-F^{l_{M-1}(\vec{\theta})}_{t_{M+1}}).
\end{eqnarray*}
Hence
\begin{eqnarray*}
\|C^{(n,m)}_t\|_{(n,m)} &\leq&
\sum_{M=1}^{\infty}
\frac{2C_0^{n+m+M-1}\,C_1^{n+m+M}\sqrt{t}}{N\,\delta_N^d}\,\|\Theta^{(n,m)}_M(t)\|_{(n,m)} \quad \text{by }(\ref{E:Bound_A_minus_F})\\
&\leq&  \sum_{M=1}^{\infty}
\frac{2C_0^{n+m+M-1}\,C_1^{n+m+M}\sqrt{t}}{N\,\delta_N^d}\,\frac{1}{\delta_N^d}\,C_1^{n+m+M}\,t^{M/2} \quad \text{by }(\ref{E:Bound_I_M_mn})\\
&\leq& \frac{(C_0\,C)^{n+m}\,t}{N\,\delta_N^{2d}}
\end{eqnarray*}
for all $t\in [\,0,\;(C_0\,C)^{-2}\,]$ and $N\geq N_0(D_{\pm})$.
\end{pf}

It follows from Theorem \ref{T:Asymptotic_nm_correlation_t} that we have
\begin{cor}\label{cor:Asymptotic_nm_correlation_t_1}

\begin{enumerate}
\item[\rm (i)]  (\textbf{Propagation of chaos})
Suppose Assumption \ref{A:ShrinkingRateLLN} holds. Then for any $T>0$ and any $(n,m)$, we have
$$\lim_{N\to\infty}\|(A^{(n,m)}_t-F^{(n,m)}_t)\|_{(n,m)}=0$$
uniformly for $t\in [0,\,T]$.
\item[\rm (ii)]   Suppose, furthermore, that Assumption \ref{A:ShrinkingRateCLT} holds. Then for any $T>0$ and any $(n,m)$, we have
$$\lim_{N\to\infty}\|C^{(n,m)}_t\|_{(n,m)}=0$$
uniformly for $t\in [0,\,T\wedge (C_0\,C)^{-2}]$.
\end{enumerate}
\end{cor}

\begin{pf}
Suppose $\liminf_{N\to\infty} N\,\delta_N^d \in (0,\infty]$. We have shown that the following upper bound of the series expansion of $\|(A^{(n,m)}_t-F^{(n,m)}_t)\|_{(n,m)}$ (in Step 2 of the proof of Theorem \ref{T:Asymptotic_nm_correlation_t}) converges uniformly in $N$.
\begin{eqnarray*}
&& \|A^{(n,m)}_t-F^{(n,m)}_t\|\\
&\leq& \;\;\int_{t_2=0}^t \frac{1}{N}\Big\|P^{(n,m)}_{t-t_2}QF^{(n,m)}_{t_2}\Big\|\\
&&+ \int_{t_2=0}^t\int_{t_3=0}^{t_2} \frac{1}{N}\Big\|P^{(n,m)}_{t-t_2}\left(\sum_{i=1}^nV_{+_i}P^{(n,m+1)}_{t_2-t_3}QF^{(n,m+1)}_{t_3} +
\sum_{j=1}^mV_{-_j}P^{(n+1,m)}_{t_2-t_3}QF^{(n+1,m)}_{t_3}\Big\|
\right)\\
&&+\cdots\\
\end{eqnarray*}
We can check that the integrand (w.r.t. $dt_2\,dt_3\,\cdots $) for each term converges to zero by Lemma \ref{L:MinkowskiContent_I_k_CLT}. Hence each term converges to zero as $N\to\infty$. Therefore, the whole series converges to zero and we obtained part (i).

The proof for part (ii) is the same, using the series expansion of $C^{(n,m)}_t$ in Step 4 of the proof of Theorem \ref{T:Asymptotic_nm_correlation_t}.
\end{pf}

\begin{remark}\rm
From Corollary \ref{cor:Asymptotic_nm_correlation_t_1}(i), we have
\begin{equation}\label{E:Asymptotic_nm_correlation_t_1b}
F^{(n,m)}_t=A_t^{(n,m)}+\frac{1}{N}\,B_t^{(n,m)}+\frac{o(N)}{N},
\end{equation}
where $o(N)$ is a term which tends to zero uniformly for $t\in  [0,\,T\wedge (C_0\,C)^{-2}]$.
\end{remark}

The following corollary of Theorem \ref{T:Asymptotic_nm_correlation_t} gives us a pointwise bound for the difference between $F^{(n+p,m+q)}$ and $F^{(n,m)}\cdot F^{(p,q)}$.
\begin{cor}\label{cor:Asymptotic_nm_correlation_t_2}
For any $T>0$ and any non-negative integers $n,\,m,\,p,\,q$, we have
\begin{eqnarray*}
&& N\,\Big|\,F^{(n+p,m+q)}_t(\vec{x},\vec{z},\,\vec{y},\vec{w})\,-\,F^{(n,m)}_t(\vec{x},\vec{y})\cdot F^{(p,q)}_t(\vec{z},\vec{w})\,\Big|\\
&\leq& C_0^{n+m+p+q-2}\,\Big|
\sum_{i,\,l}G_t(x_i,w_l)+\sum_{k,\,j}G_t(z_k,y_j)+\sum_{i,\,k}G^+_t(x_i,z_k)+\sum_{j,\,l}G^-_t(y_j,w_l)
\Big|\\
&& + \dfrac{(C_0\,C)^{n+m+p+q}\,t}{N\,\delta^{2d}}
\end{eqnarray*}
whenever $0\leq t \leq T\wedge(C_0\,C)^{-2}$ and $N\geq N_0(D_{+},D_-)$.
\end{cor}

\begin{pf}
By Theorem \ref{T:Asymptotic_nm_correlation_t} and the shorthand $F^{(n,m)}=F^{(n,m)}_t(\vec{x},\vec{y})$, we have
\begin{eqnarray*}
&& N\,\left[\,F^{(n+p,m+q)}-F^{(n,m)}\cdot F^{(p,q)} \,\right]\\
&=&N\,\left[\,\left(A+\frac{B+C}{N}\right)^{(n+p,m+q)}-\,\left(A+\frac{B+C}{N}\right)^{(n,m)}\cdot
\left(A+\frac{B+C}{N}\right)^{(p,q)}\,\right]\\
&=&\left( B^{(n+p,m+q)}-A^{(n,m)}B^{(p,q)}-A^{(p,q)}B^{(n,m)}\right)\\
&&\,+\,\left(C^{(n+p,m+q)}-A^{(n,m)}C^{(p,q)}-A^{(p,q)}C^{(n,m)}-\dfrac{(B+C)^{(n,m)}(B+C)^{(p,q)}}{N}\right).
\end{eqnarray*}
It is remarkable that all terms involving $g^+$ and $g^-$ cancel out in $B^{(n+p,m+q)}-A^{(n,m)}B^{(p,q)}-A^{(p,q)}B^{(n,m)}$ and we have control over all the remaining terms via the bounds (\ref{E:Asymptotic_nm_correlation_t}) in Theorem \ref{T:Asymptotic_nm_correlation_t}. In fact,
\begin{eqnarray}\label{E:Asymptotic_nm_correlation_t_2}
&&B^{(n+p,m+q)}(\vec{x},\vec{z},\,\vec{y},\vec{w})\,-\,A^{(n,m)}(\vec{x},\vec{y})B^{(p,q)}(\vec{z},\vec{w})-A^{(p,q)}(\vec{z},\vec{w})B^{(n,m)}(\vec{x},\vec{y}) \notag\\
&=&\,-\,A^{(n+p,m+q)}(\vec{x},\vec{z},\,\vec{y},\vec{w})\,\bigg(
\sum_{i=1}^{n}\sum_{l=1}^{q}\frac{G(x_i,w_l)}{f^+(x_i)f^-(w_l)}
+\sum_{k=1}^{p}\sum_{j=1}^{m}\frac{G(z_k,y_j)}{f^+(z_k)f^-(y_j)} \notag\\
&& \qquad \qquad \qquad \qquad +\sum_{i=1}^{n}\sum_{k=1}^{p}\frac{G^+(x_i,z_k)}{f^+(x_i)f^+(z_k)}
+\sum_{j=1}^{m}\sum_{l=1}^{q}\frac{G^-(y_j,w_l)}{f^-(y_j)f^-(w_l)}
\bigg).
\end{eqnarray}
The result now follows from the fact that $\|f^{\pm}\|\leq \|u^{\pm}_0\| \leq C_0$ and (\ref{E:Asymptotic_nm_correlation_t}).
\end{pf}

\begin{remark}\rm
(Generalizing to the case $\|F^{(n,m)}_0-A^{(n,m)}_0\|_{(n,m)}\neq 0$) In Theorem \ref{T:Asymptotic_nm_correlation_t}, we have assumed the initial error $e^{(n,m)}_{N}:= \|F^{(n,m)}_0-A^{(n,m)}_0\|_{(n,m)}$ to be zero for all $n,m$ and $N$. In fact we can weaken this condition by requiring $e^{(n,m)}_{N} \to 0$ fast enough as $N\to\infty$. This can be quantified by taking into account the contributions of the terms $F^{(n,m)}_0-A^{(n,m)}_0$ in the difference between (\ref{E:BBGKY_nm_correlation_t}) and (\ref{E:BBGKY_nm_A_t}) in Step 2. \qed
\end{remark}

\subsection{Generalized correlation functions $F^{N,(n,m),(p,q)}_{s,t}$}

The proof for Step 2 (Tightness) and Step 6 (Boltzman-Gibbs Principle) for Theorem \ref{T:Convergence_AnnihilatingSystem} require analysis not only for the correlation function at a fixed time $t$, but also for the joint probability distributions of the particles at two different times $s < t$.

\begin{definition}\label{Def:GeneralizedCorrelationFcn}
    For $n,m,p,q\in \mathbb{N}$ and $0\leq s\leq t$, we define the \textbf{generalized correlation functions} $F^{(n,m),(p,q)}_{s,t}=F^{N,(n,m),(p,q)}_{s,t}$ by
    \begin{equation}
    \int\,\Phi(\vec{x},\vec{y})\,\Psi(\vec{z},\vec{w})\,F^{(n,m),(p,q)}_{s,t}(\vec{x},\vec{y},\;\vec{z},\vec{w})\,d(\vec{x},\vec{y},\vec{z},\vec{w})= \E\left[\Phi_{(n,m)}(s)\;\Psi_{(p,q)}(t)\right]
    \end{equation}
    for all $\Phi\in C(\bar{D}_+^n\times \bar{D}_-^m)$ and $\Psi\in C(\bar{D}_+^p\times \bar{D}_-^q)$. Here $\Phi_{(n,m)}$ is defined in (\ref{E:Def_Phi_nm}) and $\Psi_{(p,q)}$ is defined in the same way.
\end{definition}

\begin{example}
For example, we have
\begin{equation*}
\E[\<\phi,\X^{N,+}_s\>\,\<\psi,\X^{N,+}_t\>]=\int_{D_+^2}\phi(x)\psi(z)\,F^{(1,0)(1,0)}_{s,t}(x,z)\,d(x,z)\quad \text{ and}
\end{equation*}
\begin{equation*}
\E\left[
\<\ell,\,\X^{N,+}_u\otimes \X^{N,-}_u\>\,\<\ell,\,\X^{N,+}_v\otimes \X^{N,-}_v\>
\right]
= \int \ell(x,y)\ell(z,w)\,F^{(1,1)(1,1)}_{u,v}(x,y,z,w)\,d(x,y,z,w)
\end{equation*}
\end{example}

To compare $F^{(n,m),(p,q)}_{u,u+r}$ and $F^{(n,m)}_u\cdot F^{(p,q)}_{u+r}$, we also define
\begin{equation}\label{Def:Enmpq}
E^{(n,m),(p,q)}_{u,r}(\vec{x},\vec{y},\;\vec{z},\vec{w}) :=
F^{(n,m),(p,q)}_{u,u+r}(\vec{x},\vec{y},\;\vec{z},\vec{w})- F^{(n,m)}_{u}(\vec{x},\vec{y})\cdot F^{(p,q)}_{u+r}(\vec{z},\vec{w})
\end{equation}

\subsubsection{A technical lemma towards tightness}

The following lemma is the key and hardest part towards the proof of the tightness
result (Theorem \ref{T:Tightness_Z}) for $\ZZ^N$.

\begin{lem}\label{L:Tightness_Y_phi_2}
Suppose Assumption \ref{A:ShrinkingRateCLT} holds. For any $T>0$, there exists $C=C(D_+,D_-,T)>0$ and $N_0=N_0(D_+,D_-)$ so that we have
\begin{equation*}
 \E\left[\left(\sqrt{N}\int_{a}^{b}
\<\ell\,\varphi_r,\;\X^{N,+}_r\otimes \X^{N,-}_r\>- \E[\<\ell\,\varphi_r,\;\X^{N,+}_r\otimes \X^{N,-}_r\>]\,dr
\right)^2\right] \leq C\,\|\varphi\|^2\,(b-a)^{3/2}
\end{equation*}
whenever $0\leq a\leq b\leq \,T_0:=T\wedge (C_0\,C)^{-2}$,  for any $N>N_0$ and any bounded function $\varphi_t(x,y)$ on $[0,\,T_0]\times D_+\times D_-$ with uniform norm $\|\varphi\|$.
\end{lem}

A direct calculation suggests that the $L^2(\P)$ norm of
$$\sqrt{N}\left(\,\<\ell,\,\X^{N,+}_r\otimes \X^{N,-}_r\>- \E[\<\ell,\,\X^{N,+}_r\otimes \X^{N,-}_r\>]\,\right)$$
blows up in the order of $1/\delta$ (for $r>0$), due to the fact that $\int \ell(x,y_1)\ell(x,y_2)\,d(x,y_1,y_2)$ is of order $1/\delta$. Hence we need to look into the generalized correlation functions.

\medskip

\begin{pf}
\textbf{Step 1: Write LHS in terms of the generalized correlation functions. }

Note that $(\int_a^b f(r)\,dr)^2= 2\int_{u=a}^b \int_{v=u}^b f(u)f(v)= 2\int_{u=a}^b \int_{t=0}^{b-u} f(u)f(t)$ by Fubinni's Theorem followed by the change of variable $t=v-u$. Hence Lemma \ref{L:Tightness_Y_phi_2} is implied by
\begin{equation}\label{e:6.34a}
\int_{u=a}^b \int_{t=0}^{b-u}\,\int_{\substack{(x,y)\in D_+\times D_-\\(\tilde{x},\tilde{y})\in D_+\times D_-}} N\,\ell(x,y)\ell(\tilde{x},\tilde{y})\,E^{(1,1)(1,1)}_{u,t}(x,y,\,\tilde{x},\tilde{y}) \;\leq C\,(b-a)^{3/2},
\end{equation}
where $E^{(1,1)(1,1)}_{u,t}$ is defined in (\ref{Def:Enmpq}). The ideas is to first obtain a `variation of constant' formula for $E_{u,t}$ via the Dynkin's formula; then iterate the formula to obtain a series expansion of $E_{u,t}$ in terms of $E_{u,0}$;
and finally
estimate $E_{u,0}$ and each term of the series.

\textbf{Step 2: Estimate $|E^{(1,1)}_{t}|$ in terms of $\{E^{(p,q)}_0\}$. }

Applying Dynkin's formula as in (\ref{E:BBGKY_nm_correlation_t}) yields
\begin{eqnarray}\label{E:BBGKY_nmpq_correlation}
F^{(n,m),(p,q)}_{u,u+r} &=& P^{(p,q)}_r F^{(n,m),(p,q)}_{u,u} \\
&& \,- \int_0^r P^{(p,q)}_{r-\theta}\left( V^+F^{(n,m),(p,q+1)}_{u,u+\theta} + V^-F^{(n,m),(p+1,q)}_{u,u+\theta} + \frac{Q}{N}F^{(n,m),(p,q)}_{u,u+\theta} \right)\;d\theta, \notag
\end{eqnarray}
where $P^{(p,q)}_t$, $V^+$, $V^-$ and $Q$ are operators defined as before and act on the $(\vec{z},\vec{w})$ variables.

Fix $u\geq 0$, $(n,m)$ and $(\vec{x},\vec{y})\in \bar{D}_+^n\times \bar{D}_-^m$, and  write
\begin{equation*}
E^{(p,q)}_{r}(\vec{z},\vec{w}) :=
E^{(n,m),(p,q)}_{u,r}(\vec{x},\vec{y},\;\vec{z},\vec{w}) \quad\text{for notational simplicity.}
\end{equation*}
Then (\ref{E:BBGKY_nmpq_correlation}) yields
\begin{equation}\label{E:BBGKY_Enmpq}
E^{(p,q)}_{r}=P^{(p,q)}_r E^{(p,q)}_{0} - \int_0^r P^{(p,q)}_{r-\theta}\left( V^+E^{(p,q+1)}_{\theta} + V^-E^{(p+1,q)}_{\theta} + \frac{Q}{N}E^{(p,q)}_{\theta} \right)\;d\theta,
\end{equation}
where $P^{(p,q)}_t$, $V^+$, $V^-$ and $Q$ are operators defined before, acting on the $(\vec{z},\vec{w})$ variables. In other words,
$(t,\,(\vec{z},\vec{w}))\mapsto E^{(p,q)}_{t}(\vec{z},\vec{w})$ is the probabilistic solution of
    \begin{equation*}
        \left\{\begin{aligned}
        \dfrac{\partial E}{\partial t} &= \frac{1}{2}\Delta E - \frac{Q}{N}E - \left( V^+E^{(p,q+1)} + V^-E^{(p+1,q)}\right)
         &&\qquad  \text{on }(0,\infty)\times D_+^p\times D_-^q , \\
        \dfrac{\partial E}{\partial \vec{n}} &= 0 &&\qquad \text{on } (0,\infty)\times\, \partial (D_+^p\times D_-^q), \\
        E_0(\cdot)&= F^{(n,m),(p,q)}_{u,u}(\vec{x},\vec{y},\;\cdot\,)- F^{(n,m)}_{u}(\vec{x},\vec{y})\cdot F^{(p,q)}_{u}(\cdot)
        &&\qquad \text{on }  D_+^p\times D_-^q.
        \end{aligned}\right.
    \end{equation*}
It can be shown (see Proposition 2.19 in \cite{zqCwtF13a} for a proof) that the following probabilistic representation holds true for $E=E^{(p,q)}$:
\begin{equation}
E_t(\vec{z},\vec{w})=\E^{\vec{z},\vec{w}}\left[E_0(X_t)e^{-\int_0^tk(X_s)ds}-\int_0^t g(t-\theta,X_{\theta})e^{-\int_0^{\theta}k(X_s)ds}\,d\theta\right],
\end{equation}
where $k=\frac{Q^{(p,q)}}{N}$, $g(t)=V^+E_t^{(p,q+1)} + V^-E_t^{(p+1,q)}$ and $X_t$ is the RBM in $D_+^p\times D_-^q$ starting at $(\vec{z},\vec{w})$.
From this, the triangle inequality and the non-negativity of $k=\frac{Q^{(p,q)}}{N}$, we have
\begin{equation*}
\Big|E^{(p,q)}_{t}\Big| \leq P^{(p,q)}_t(|E^{(p,q)}_{0}|)
+ \int_0^t P^{(p,q)}_{t-t_2}\left(\Big| V^+E^{(p,q+1)}_{t_2} + V^-E^{(p+1,q)}_{t_2} \Big|\right)\;dt_2.
\end{equation*}
It then follows that almost everywhere in $D_+^p \times D_-^q$, we have
\begin{equation}\label{E:Iterate_E_pq}
\Big|E^{(p,q)}_{t}\Big| \leq P^{(p,q)}_t(|E^{(p,q)}_{0}|)
+ \int_0^t P^{(p,q)}_{t-t_2}\left( V^+ |E^{(p,q+1)}_{t_2}| + V^-|E^{(p+1,q)}_{t_2}|\right)\;dt_2.
\end{equation}
Now we iterate (\ref{E:Iterate_E_pq}) to obtain
\begin{eqnarray}\label{E:Iterate_E_pq_2}
\Big|E^{(p,q)}_{t}\Big| &\leq& \;\; P^{(p,q)}_t |E^{(p,q)}_{0}|
+  \int_{t_2=0}^t P^{(p,q)}_{t-t_2}\left( V^+P^{(p,q+1)}_{t_2}|E^{(p,q+1)}_0| + V^-P^{(p+1,q)}_{t_2}|E^{(p+1,q)}_0| \right) \notag\\
&&+
 \int_{t_2=0}^t \int_{t_3=0}^{t_2}
P^{(p,q)}_{t-t_2}\bigg(
V^+P^{(p,q+1)}_{t_2-t_3}\,\Big(\,V^+P^{(p,q+2)}_{t_3} |E_0| + V^-P^{(p+1,q+1)}_{t_3} |E_0|\,\Big) \notag\\
&&\qquad \qquad \qquad \qquad \quad
+ V^-P^{(p+1,q)}_{t_2-t_3}\,\Big(\,V^+P^{(p+1,q+1)}_{t_3} |E_0| + V^-P^{(p+2,q)}_{t_3} |E_0|\,\Big)\,\bigg) \notag\\
&& + \cdots \notag\\
&=&
\sum_{M=0}^{\infty}\int_{t_2=0}^t\int_{t_3=0}^{t_2}\cdots \int_{t_{M+1}=0}^{t_M} \notag\\
&&\qquad \; \sum_{\vec{\theta}\in \,\mathbb{T}^{(n,m)}_{M}}\,
P^{(n,m)}_{t-t_2}\,V_{\theta_1}\,P^{l_1(\vec{\theta})}_{t_2-t_3}\,V_{\theta_2}\,P^{l_2(\vec{\theta})}_{t_3-t_4}\,V_{\theta_3}\cdots P^{l_{M-1}(\vec{\theta})}_{t_M-t_{M+1}}\,V_{\theta_N}\,|E^{l_{M}(\vec{\theta})}_{0}|.
\end{eqnarray}

From this inequality and the triangle inequality, we have, for any $u\geq 0$, $(n,m)=(1,1)$ and $(x,y)\in D_+\times D_-$,
\begin{eqnarray}\label{E:Tightness_Y_phi_2_Step2}
&& \int\ell\,\Big|E^{(1,1)}_{t}\Big| := \int_{(z,w)\in D_+\times D_-} \ell(z,w)\,\Big|E^{(1,1)}_{t}(x,y,z,w)\Big| \notag\\
&\leq& \int \Psi^{(root)}\Big|E^{(1,1)}_0\Big| \notag\\
&&\,+ \int_{0}^t\,\left( \int\Psi^{(+_1)}\Big|E^{(1,2)}_0\Big|+\int\Psi^{(-_1)}\Big|E^{(2,1)}_0\Big|\right)\,dt_2 \notag\\
&&\,+ \int_{0}^t\int_{0}^{t_2}\,\bigg(
\int\Psi^{(+_1,+_2)}\Big|E^{(1,3)}_0\Big|
+\int\Psi^{(+_1,-_1)}\Big|E^{(2,2)}_0\Big|
+\int\Psi^{(+_1,-_2)}\Big|E^{(2,2)}_0\Big| \notag\\
&& \qquad \qquad
+\int\Psi^{(-_1,+_1)}\Big|E^{(2,2)}_0\Big|
+\int\Psi^{(-_1,+_2)}\Big|E^{(2,2)}_0\Big|
+\int\Psi^{(-_1,-_2)}\Big|E^{(3,1)}_0\Big|
\bigg)\,dt_3\,dt_2 \notag\\
&& \,+\, \cdots \notag\\
&=& \sum_{M=0}^{\infty}\int_{t_2=0}^t\int_{t_3=0}^{t_2}\cdots\int_{t_{M+1}=0}^{t_M} \left( \sum_{\vec{\theta}\in \,\mathbb{T}^{(1,1)}_M}  \int\Psi^{\vec{\theta}}\,\Big|E^{l_M(\vec{\theta})}_0\Big|\right)\,dt_{M+1}\,\cdots\,dt_3\,dt_2,
\end{eqnarray}
where the integral sign for $E^{(p,q)}_0$ is on the set $D_+^p\times D_-^q$,
\begin{eqnarray*}
\Psi^{(root)}(z,w) &:=& P^{(1,1)}_t\ell(z,w)\\
\Psi^{(+_1)}(z,w_1,w_2)&:=& P^{(1,2)}_{t_2}\left( (P^{(1,1)}_{t_1-t_2}\ell)(a_1,b_1)\cdot \ell(a_1,b_2)\right)(z,w_1,w_2)\\
\Psi^{(-_1)}(z_1,z_2,w)&:=& P^{(2,1)}_{t_2}\left( (P^{(1,1)}_{t_1-t_2}\ell)(a_1,b_1)\cdot \ell(a_2,b_1)\right)(z_1,z_2,w).
\end{eqnarray*}
Inductively, $\Psi^{(\vec{\theta},+_i)} \in C(\bar{D}_+^p\times \bar{D}_-^{q+1})$ and $\Psi^{(\vec{\theta},-_j)} \in C(\bar{D}_+^{p+1}\times \bar{D}_-^{q})$ are obtained from $\Psi^{\vec{\theta}}$ as follows: if $\Psi^{\vec{\theta}}(\vec{z},\vec{w})= P^{(p,q)}_{t_{M+1}}F(\vec{z},\vec{w})$, then
\begin{eqnarray*}
\Psi^{(\vec{\theta},+_i)}(\vec{z},(\vec{w},w_{q+1})) &:=&
P^{(p,q+1)}_{t_{M+2}}\left( (P^{(p,q)}_{t_{M+1}-t_{M+2}}\,F)(\vec{a},\vec{b})\; \ell(a_i,b_{q+1}) \right)(\vec{z},(\vec{w},w_{q+1})) \quad\text{and}\\
\Psi^{(\vec{\theta},-_j)}((\vec{z},z_{p+1}),\vec{w}) &:=&
P^{(p+1,q)}_{t_{M+2}}\left( (P^{(p,q)}_{t_{M+1}-t_{M+2}}\,F)(\vec{a},\vec{b})\; \ell(a_{p+1},b_{j}) \right)((\vec{z},z_{p+1}),\vec{w}).
\end{eqnarray*}

\textbf{Step 3: Estimate $E^{(p,q)}_0=F^{(1,1),(p,q)}_{u,u}- F^{(1,1)}_u\cdot F^{(p,q)}_u$. }

For any $\Psi\in C(\bar{D}_+^p\times \bar{D}_-^q)$, by Definition \ref{Def:GeneralizedCorrelationFcn} we have
\begin{eqnarray}\label{E:FuuFu}
&&\int_{D_+^{p+1}\times D_-^{q+1}}\,\ell(x,y)\,\Psi(\vec{z},\vec{w})\,F^{(1,1),(p,q)}_{u,u}(x,y,\,\vec{z},\vec{w})\,d(x,y,\,\vec{z},\vec{w}) \notag\\
&=&\E\left[
\bigg(\dfrac{1}{N^2}\sum_{i}^{\sharp_u}\sum_{j}^{\sharp_u}\ell(X^{i}_u,\,Y^{j}_u)\bigg)\,
\bigg(\dfrac{1}{N^{(p)}\,N^{(q)}}
\sum_{\substack{\text{distinct} \\ k_1,\cdots k_p =1}}^{\sharp_u}\,\sum_{\substack{\text{distinct} \\
l_1,\cdots l_q =1  }}^{\sharp_u}
\,\Psi(X^{k_1}_u,\cdots, X^{k_p}_u,\,Y^{l_1}_u,\cdots,Y^{l_q}_u)\bigg)
\right] \notag\\
&=& \,\dfrac{N^{(p+1)}N^{(q+1)}}{N^2N^{(p)}N^{(q)}}\,\int_{D_+^{p+1}\times D_-^{q+1}}\ell(x,y)\,\Psi(\vec{z},\vec{w})\,F^{(p+1,q+1)}_{u}((x,\vec{z}),\,(y,\vec{w})) \notag\\
&&\,+\,\dfrac{N^{(q+1)}}{N^2\,N^{(q)}}\,\sum_{i=1}^p \int_{D_+^{p}\times D_-^{q+1}}\ell(z_i,y)\,\Psi(\vec{z},\vec{w})\,F^{(p,q+1)}_{u}(\vec{z},\,(y,\vec{w})) \notag\\
&&\,+\,\dfrac{N^{(p+1)}}{N^2\,N^{(p)}}\,\sum_{j=1}^q \int_{D_+^{p+1}\times D_-^{q}}\ell(x,w_j)\,\Psi(\vec{z},\vec{w})\,F^{(p+1,q)}_{u}((x,\vec{z}),\,\vec{w}) \notag\\
&&\,+\,\dfrac{1}{N^2}\,\sum_{i=1}^p\sum_{j=1}^q \int_{D_+^{p}\times D_-^{q}}\ell(z_i,w_j)\,\Psi(\vec{z},\vec{w})\,F^{(p,q)}_{u}(\vec{z},\,\vec{w}).
\end{eqnarray}
This connects $F^{(1,1),(p,q)}_{u,u}$ to $F^{(p+1,q+1)}_{u}$ and
we know more about the latter  (such as Theorem \ref{T:Asymptotic_nm_correlation_t}).
Furthermore, we use the simple fact that $\int f\,|g|= \int\tilde{f}\,g$
where
$\tilde{f}(x)=
\begin{cases}
f(x),  &\text{ if } g(x)\geq 0 \\
-f(x), &\text{ if } g(x)<0.
\end{cases}$
Therefore, for any $\Psi\in C_+(\bar{D}_+^p\times \bar{D}_-^q)$, we have
\begin{eqnarray*}
&&\int_{D_+^{p+1}\times D_-^{q+1}}\,\ell(x,y)\,\Psi(\vec{z},\vec{w})\,\Big|E^{(p,q)}_{0}(x,y,\vec{z},\vec{w})\Big|\,d(x,y,\,\vec{z},\vec{w})\\
&\leq&  \dfrac{(N-p)(N-q)}{N^2}\int_{D_+^{p+1}\times D_-^{q+1}}\ell\,\Psi\,\Big|F^{(p+1,q+1)}_{u}-F^{(1,1)}_u\cdot F^{(p,q)}_u\Big|\\
&& + \Big| \dfrac{(N-p)(N-q)}{N^2}-1 \Big|\,\left(\int_{D_+\times D_-}\ell\,F^{(1,1)}_u\right)\,\left(\int_{D_+^{p}\times D_-^{q}} \Psi\,F^{(p,q)}_u \right)\\
&&+\,\dfrac{N-q}{N^2}\,\sum_{i=1}^p \int_{D_+^{p}\times D_-^{q+1}}\ell(z_i,y)\,\Psi(\vec{z},\vec{w})\,F^{(p,q+1)}_{u}(\vec{z},\,(y,\vec{w})) \\
&&+\,\dfrac{N-p}{N^2}\,\sum_{j=1}^q \int_{D_+^{p+1}\times D_-^{q}}\ell(x,w_j)\,\Psi(\vec{z},\vec{w})\,F^{(p+1,q)}_{u}((x,\vec{z}),\,\vec{w}) \\
&&+\,\dfrac{1}{N^2}\,\sum_{i=1}^p\sum_{j=1}^q \int_{D_+^{p}\times D_-^{q}}\ell(z_i,w_j)\,\Psi(\vec{z},\vec{w})\,F^{(p,q)}_{u}(\vec{z},\,\vec{w})\\
&\leq&  \int_{D_+^{p+1}\times D_-^{q+1}}\ell\,\Psi\,\Big|F^{(p+1,q+1)}_{u}-F^{(1,1)}_u\cdot F^{(p,q)}_u\Big|\\
&& + \dfrac{p+q}{N}\,\left(\int_{D_+\times D_-}\ell\,F^{(1,1)}_u\right)\,\left(\int_{D_+^{p}\times D_-^{q}} \Psi\,F^{(p,q)}_u \right)\\
&&+\,\dfrac{1}{N}\,\sum_{i=1}^p \int_{D_+^{p}\times D_-^{q+1}}\ell(z_i,y)\,\Psi(\vec{z},\vec{w})\,F^{(p,q+1)}_{u}(\vec{z},\,(y,\vec{w})) \\
&&+\,\dfrac{1}{N}\,\sum_{j=1}^q \int_{D_+^{p+1}\times D_-^{q}}\ell(x,w_j)\,\Psi(\vec{z},\vec{w})\,F^{(p+1,q)}_{u}((x,\vec{z}),\,\vec{w}) \\
&&+\,\dfrac{1}{N^2}\,\sum_{i=1}^p\sum_{j=1}^q \int_{D_+^{p}\times D_-^{q}}\ell(z_i,w_j)\,\Psi(\vec{z},\vec{w})\,F^{(p,q)}_{u}(\vec{z},\,\vec{w}).
\end{eqnarray*}

On the other hand, by Corollary \ref{cor:Asymptotic_nm_correlation_t_2},
\begin{eqnarray}
&&N\,\Big|F^{(p+1,q+1)}_{u}(x,y,\vec{z},\vec{w})-F^{(1,1)}_u(x,y)\cdot F^{(p,q)}_u(\vec{z},\vec{w})\Big|\\
&\leq&  C_0^{p+q}\,\left( \sum_{i=1}^p G_u(z_i,y)+\sum_{j=1}^q G_u(x,w_j)+\sum_{i=1}^p G^+_u(x,z_i)+\sum_{j=1}^q G^-_u(y,w_j)\right)
+ \frac{(C_0\,C)^{p+q+2}\,u}{N\,\delta_N^{2d}}. \notag
\end{eqnarray}

Combining with the calculation just before the
proceeding inequality, we obtain

\begin{eqnarray}\label{E:Epq_0}
&& N\,\int_{D_+^{p+1}\times D_-^{q+1}}\,\ell(x,y)\,\Psi(\vec{z},\vec{w})\,\Big|E^{(p,q)}_{0}(x,y,\vec{z},\vec{w})\Big|\,\,d(x,y,\,\vec{z},\vec{w})\\
&\leq&  C_0^{p+q}\,\int_{D_+^{p+1}\times D_-^{q+1}}\ell(x,y)\,\Psi(\vec{z},\vec{w})\,\bigg( \sum_{i=1}^p G_u(z_i,y)+\sum_{j=1}^q G_u(x,w_j)+\sum_{i=1}^p G^+_u(x,z_i)+\sum_{j=1}^q G^-_u(y,w_j)\bigg) \notag\\
&& +\,\frac{(C_0\,C)^{p+q+2}\,u}{N\,\delta_N^{2d}} \int_{D_+^{p}\times D_-^{q}}\Psi(\vec{z},\vec{w})
 +\,C_0^{p+q+1}\,\int_{D_+^{p}\times D_-^{q+1}} \sum_{i=1}^p \ell(z_i,y)\,\Psi(\vec{z},\vec{w}) \notag\\
&&+\,C_0^{p+q+1}\,\int_{D_+^{p+1}\times D_-^{q}} \sum_{j=1}^q \ell(x,w_j)\,\Psi(\vec{z},\vec{w})
+\,\dfrac{C_0^{p+q}}{N}\,\int_{D_+^{p}\times D_-^{q}} \sum_{i=1}^p\sum_{j=1}^q\ell(z_i,w_j)\,\Psi(\vec{z},\vec{w}).  \notag
\end{eqnarray}

\textbf{Step 4: Final estimates. }

We now put $\Psi=\Psi^{\vec{\theta}}$ into inequality (\ref{E:Epq_0}) for each $\Psi^{\vec{\theta}}$ that appears in (\ref{E:Tightness_Y_phi_2_Step2}) at the end of Step 2. Specifically, by (\ref{E:Tightness_Y_phi_2_Step2}) and (\ref{E:Epq_0}) respectively, we have
\begin{eqnarray}\label{E:Tightness_Y_phi_2_Step4}
&& N\int_{(x,y)\in D_+\times D_-}\int_{(\tilde{x},\tilde{y})\in D_+\times D_-}  \ell(x,y)\,\ell(\tilde{x},\tilde{y})\,\Big|E^{(1,1)}_{t}(x,y,\,\tilde{x},\tilde{y})\Big| \notag\\
&\leq& \sum_{M=0}^{\infty}\int_{t_2=0}^t\cdots\int_{t_{M+1}=0}^{t_M}
 \sum_{\vec{\theta}\in \,\mathbb{T}^{(1,1)}_M} \bigg[\, N\int_{(x,y)} \ell(x,y)\, \int_{(\vec{z},\vec{w})\in D_+^{p}\times D_-^{q}} \Psi^{\vec{\theta}}(\vec{z},\vec{w})\,\Big|E^{(1,1),\,l_M(\vec{\theta})}_{u,0}(x,y,\vec{z},\vec{w})\Big|\,\bigg]\notag\\
&\leq& \sum_{M=0}^{\infty}\int_{t_2=0}^t\cdots\int_{t_{M+1}=0}^{t_M}
 \sum_{\vec{\theta}\in \,\mathbb{T}^{(1,1)}_M} \bigg[\,\sum_{i=1}^5 \Theta_i^{l_M(\vec{\theta})}\Big(\Psi^{\vec{\theta}}\Big)\,\bigg],
\end{eqnarray}
where in the first inequality, the integration over the variables $(\vec{z},\vec{w})$ is on $D_+^{p}\times D_-^{q}$ where $l_M(\vec{\theta})=(p,q)$; in the second inequality, $\Theta_i^{l_M(\vec{\theta})}\Big(\Psi\Big)$ is the $i$-th term that appear on the RHS of (\ref{E:Epq_0}).

We will estimate each of the five terms ($i=1,2,3,4,5$) on the RHS of (\ref{E:Tightness_Y_phi_2_Step4}) separately. The arguments are the same for all of them. We first consider the term for $i=2$. This term is
\begin{eqnarray}\label{E:Tightness_Y_phi_2_Step4_2}
&& \sum_{M=0}^{\infty}\int_{t_2=0}^t\cdots\int_{t_{M+1}=0}^{t_M}
\sum_{\vec{\theta}\in \,\mathbb{T}^{(1,1)}_M} \bigg[\,\Theta_2^{l_M(\vec{\theta})}\Big(\Psi^{\vec{\theta}}\Big)\,\bigg] \notag\\
&=& \sum_{M=0}^{\infty}\int_{t_2=0}^t\cdots\int_{t_{M+1}=0}^{t_M}
\sum_{\vec{\theta}\in \,\mathbb{T}^{(1,1)}_M} \bigg[\,\frac{(C_0\,C)^{M+4}\,u}{N\,\delta_N^{2d}} \int_{D_+^{p}\times D_-^{q}}\Psi^{\vec{\theta}}(\vec{z},\vec{w})\,\bigg],
\end{eqnarray}
where we have used the fact that the sum of the two components of $l_M(\vec{\theta})$ is $M+2$ (i.e. $p+q=M+2$). Using the same argument of Step 3 in the proof of Theorem \ref{T:Asymptotic_nm_correlation_t}, we have, for each $M\geq 1$,
\begin{eqnarray*}
&&\int_{t_2=0}^t\cdots\int_{t_{M+1}=0}^{t_M}\sum_{\vec{\theta}\in \,\mathbb{T}^{(1,1)}_M} \int_{D_+^{p}\times D_-^{q}}\Psi^{\vec{\theta}}(\vec{z},\vec{w}) \\
&\leq& \int_{t_2=0}^t\cdots\int_{t_{M+1}=0}^{t_M}\sum_{\vec{\upsilon}\,\in \mathbb{S}^{(1,1)}_M}\frac{C^M}{\sqrt{(t_{\upsilon_1}-t_2)\,(t_{\upsilon_2}-t_3)\cdots(t_{\upsilon_M}-t_{M+1})}} \\
&\leq& C^M\,t^{M/2}
\end{eqnarray*}
for $N\geq N(D_+,D_-)$, where $C=C(D_+,D_-,T)>0$. This inequality implies that (\ref{E:Tightness_Y_phi_2_Step4_2}) is at most
$$\frac{u}{N\,\delta_N^{2d}}\,(C_0\,C)^{4}\sum_{M=0}^{\infty}(C_0\,C)^{M}\,t^{M/2} \leq  \frac{C_0^4\,C\,u}{N\delta_N^{2d}}$$
when  $0\leq t \leq (C_0\,C)^{-2}$ and $N$ is large enough, where $C=C(D_+,D_-,T)>0$.

For $i=1$, we only need to invoke Lemma \ref{L:BoundsForGGG} and then use the same argument for $i=2$. The term on the RHS of (\ref{E:Tightness_Y_phi_2_Step4}) for $i=1$ is at most
$$\frac{C_0^4\,C\,\sqrt{u}}{\sqrt{t}}+\frac{C_0^5\,C\,u}{\sqrt{t}}.$$

For $i=3$, the term on the RHS of (\ref{E:Tightness_Y_phi_2_Step4}) is equal to
\begin{eqnarray}\label{E:Tightness_Y_phi_2_Step4_3}
&& \sum_{M=0}^{\infty}\int_{t_2=0}^t\cdots\int_{t_{M+1}=0}^{t_M}
\sum_{\vec{\theta}\in \,\mathbb{T}^{(1,1)}_M} \bigg[\,C_0^{M+3}\,\int_{D_+^{p}\times D_-^{q+1}} \sum_{i=1}^p \ell(z_i,y)\,\Psi^{\vec{\theta}}(\vec{z},\vec{w})
\,\bigg].
\end{eqnarray}
By the same argument as that for $i=2$, we have, for each $M\geq 1$,
\begin{eqnarray*}
&&\int_{t_2=0}^t\cdots\int_{t_{M+1}=0}^{t_M} \sum_{\vec{\theta}\in \,\mathbb{T}^{(1,1)}_M} \int_{D_+^{p}\times D_-^{q+1}} \sum_{i=1}^p \ell(z_i,y)\,\Psi^{\vec{\theta}}(\vec{z},\vec{w}) \\
&\leq& C^M\int_{t_2=0}^t\cdots\int_{t_{M}=0}^{t_{M-1}}\int_{t_{M+1}=0}^{t_M}\Big(\frac{M+1}{\sqrt{t_{M+1}}}\Big) \sum_{\vec{\upsilon}\,\in \mathbb{S}^{(1,1)}_M}
\frac{1}{\sqrt{(t_{\upsilon_1}-t_2)\,(t_{\upsilon_2}-t_3)\cdots(t_{\upsilon_M}-t_{M+1})}}\\
&\leq& C^M\int_{t_2=0}^t\cdots\int_{t_{M}=0}^{t_{M-1}}(M+1)\sum_{\vec{\upsilon}\,\in \mathbb{S}^{(1,1)}_{M-1}}
\frac{1}{\sqrt{(t_{\upsilon_1}-t_2)\,(t_{\upsilon_2}-t_3)\cdots(t_{\upsilon_{M-1}}-t_{M})}}\\
&&\qquad\qquad\qquad \qquad \cdot\bigg(\int_{t_{M+1}=0}^{t_M}\frac{M+1}{\sqrt{(t_M-t_{M+1})\,t_{M+1}}}\,dt_{M+1}\bigg)\\
&=& C^M\,(M+1)^2\,\pi\int_{t_2=0}^t\cdots\int_{t_{M}=0}^{t_{M-1}}\sum_{\vec{\upsilon}\,\in \mathbb{S}^{(1,1)}_{M-1}}
\frac{1}{\sqrt{(t_{\upsilon_1}-t_2)\,(t_{\upsilon_2}-t_3)\cdots(t_{\upsilon_{M-1}}-t_{M})}}\\
&\leq& C^M\,t^{(M-1)/2},
\end{eqnarray*}
where we have used the facts that $\int_{0}^{t_M}\frac{1}{\sqrt{s(t_M-s)}}\,ds = \pi$ and that $\upsilon_M\leq M$. The extra factor $(M+1)$ in the second inequality comes from the number of children (in $\mathbb{S}^{(1,1)}_{M}$) for each leaf in $\mathbb{S}^{(1,1)}_{M-1}$. Therefore, (\ref{E:Tightness_Y_phi_2_Step4_3}) is at most
$$C_0^3\sum_{M=0}^{\infty}(C_0\,C)^{M}\,t^{(M-1)/2} \leq  \frac{C_0^3\,C}{\sqrt{t}}.$$

The term for $i=4$ is symmetric to that of $i=3$, hence the upper bound is of the same form.

Finally, the term for $i=5$ can be compared to the term for $i=2$ directly, since $\sup_{(x,y)}\frac{\ell(x,y)}{N} \leq \delta^{d-1}\leq 1$ under Assumption \ref{A:ShrinkingRateCLT} and hence we can ignore the factor $\ell(z_i,w_j)$. Therefore, the term for $i=5$ is at most $C_0^4\,C$.

From the above five estimates for the RHS of (\ref{E:Tightness_Y_phi_2_Step4}), it follows that
\begin{eqnarray*}
&& N\int_{(x,y)\in D_+\times D_-}\int_{(\tilde{x},\tilde{y})\in D_+\times D_-}  \ell(x,y)\,\ell(\tilde{x},\tilde{y})\,\Big|E^{(1,1)}_{t}(x,y,\,\tilde{x},\tilde{y})\Big| \\
&\leq&   C\,\bigg( \frac{C_0^4\,\sqrt{u}}{\sqrt{t}}+\frac{C_0^5\,u}{\sqrt{t}}+\frac{C_0^4\,u}{N\delta_N^{2d}}+\frac{C_0^3}{\sqrt{t}}+C_0^4 \bigg)
\end{eqnarray*}
for $N\geq N(D_+,D_-)$, where $C=C(D_+,D_-,T)>0$.

This proves \eqref{e:6.34a} and hence the lemma.
\end{pf}

\subsubsection{A technical lemma towards Boltzman-Gibbs principle}

The goal for this subsection is to prove the following lemma, which is an indicator of the validity of the Boltzman-Gibbs principle for our annihilating diffusion model. It is instructive to compare the statement of Lemma \ref{L:BGprinciple} below with that of Lemma \ref{L:Tightness_Y_phi_2}.

\begin{lem}\label{L:BGprinciple}
Suppose Assumption \ref{A:ShrinkingRateCLT} holds. For any $T>0$, there exists $C=C(D_+,D_-,T)>0$, $N_0=N_0(D_+,D_-)$ and positive constants $\{C_N\}$ satisfying $\lim_{N\to\infty}C_N=0$ such that
\begin{equation*}
\E\Big[\,\Big|\int_0^t \ZZ^{N}_s\Big(\<\ell\varphi_s,\,f^-_{s}\>_-,\;\<\ell\varphi_s,\,f^+_{s}\>_+\Big)\,-\,\sqrt{N}\Big(\,\<\ell\varphi_s,\,\otimes_s\>- \E[\<\ell\varphi_s,\,\otimes_s\>]\,\Big)\,ds\Big|^2\,\Big] \leq C_N\,\|\varphi\|^2\,t^{3/2}
\end{equation*}
whenever $0\leq t\leq\,T_0:=T\wedge (C_0\,C)^{-2}$, for any $N>N_0$ and any bounded function $\varphi_t(x,y)$ on $[0,\,T_0]\times D_+\times D_-$ with uniform norm $\|\varphi\|$. Here $\otimes_s:=\X^{N,+}_s\otimes \X^{N,-}_s$ in abbreviation.
\end{lem}

\begin{pf}
The proof follows from the same argument that we used for Lemma \ref{L:Tightness_Y_phi_2}. Namely, we first write the LHS in terms of the generalized correlation functions (more specifically in terms of $E_{u,r}=E^{(n,m),(p,q)}_{u,r}$ defined in (\ref{Def:Enmpq})); we then bound $E_{u,r}$ in terms of $E_{u,0}$ via (\ref{E:Iterate_E_pq_2}); finally we estimate $E_{u,0}$. However, unlike Lemma \ref{L:Tightness_Y_phi_2}, the LHS here \emph{vanishes} in the limit due to a 'magical cancelations' of the first two terms in the asymptotic expansion of the correlation functions. See (\ref{E:MagicCancel_1}) and (\ref{E:MagicCancel_2}) in the proof below.

\textbf{Step 1: Abbreviations and notations. }

To avoid unnecessary complications, we assume $\varphi_s=1$ in the proof. The general case follows from a routine modification. By the fact (which follows from Fubinni's theorem and a change of variable $r=v-u$)
\begin{equation*}
\Big(\int_0^t h(s)\,ds\Big)^2= 2\int_{u=0}^t \int_{v=u}^t h(u)h(v)\,dv\,du = 2\int_{u=0}^t \int_{r=0}^{t-u} h(u)h(r)\,dr\,du,
\end{equation*}
we have
\begin{eqnarray}\label{E:LHS_BGprinciple}
&&\E\Big[\,\Big|\int_0^t \ZZ^{N}_s\Big(\<\ell,\,f^-_{s}\>_-,\;\<\ell,\,f^+_{s}\>_+\Big)\,-\,\sqrt{N}\Big(\,\<\ell,\,\otimes_s\>- \E[\<\ell,\,\otimes_s\>]\,\Big)\,ds\Big|^2\,\Big] \notag\\
&=& N\,\E\Big[\,\Big|\int_0^t  \<\alpha_s,\,\X^{N,+}_s\>+ \<\beta_s,\,\X^{N,-}_s\>-\<\ell,\,\otimes_s\> \,-\,\E[\<\alpha_s,\,\X^{N,+}_s\>+ \<\beta_s,\,\X^{N,-}_s\>-\<\ell,\,\otimes_s\>] \,ds\Big|^2\,\Big] \notag\\
&=& N\,\E\Big[\,\Big|\int_0^t (\eta_s - \xi_s)\,-\,\E[\eta_s-\xi_s] \,ds\Big|^2\,\Big] \notag\\
&=& 2N\,\int_0^t \int_u^t \E[(\eta_u - \xi_u)(\eta_v - \xi_v)]\,-\,\E[\eta_u - \xi_u]\cdot\E[\eta_v - \xi_v]\,dv\,du,
\end{eqnarray}
where we have used
the abbreviations
\begin{equation}\label{E:Abbreviate_Alpha_Beta}
\alpha_s:=\<\ell,\,f^-_{s}\>_-,\;\beta_s:=\<\ell,\,f^+_{s}\>_+,\;\eta_s:=\<\alpha_s,\,\X^{N,+}_s\>+ \<\beta_s,\,\X^{N,-}_s\>\, \text{ and }\,\xi_s:=\<\ell,\,\otimes_s\>.
\end{equation}
Note that we have, for example,
\begin{eqnarray*}
\E[\eta_s]&=& \int_{D_+} \alpha_s(x)\,F^{(1,0)}_s(x)\,dx + \int_{D_-} \beta_s(y)\,F^{(0,1)}_s(y)\,dy \\
&=& \int_{D_+\times D_-} \ell(x,y)\,f^-_{s}(y)\,F^{(1,0)}_s(x) \,+\,\ell(x,y)\,f^+_{s}(x)\,F^{(0,1)}_s(y)\,dx\,dy.
\end{eqnarray*}

\textbf{Step 2: Write LHS in terms of correlation functions. }

Direct calculation yields
\begin{eqnarray*}
&& \E[(\eta_u - \xi_u)(\eta_v - \xi_v)]= \E[\eta_u\eta_v  - \eta_v\xi_u - \eta_u\xi_v + \xi_u\xi_v]\\
&=& \int_{D_+^2} \alpha_u(x_1)\alpha_v(x_2)\,F^{(10)(10)}_{u,v}(x_1,x_2)
    + \int_{D_+\times D_-} \alpha_u(x_1)\beta_v(y_2)\,F^{(10)(01)}_{u,v}(x_1,y_2)\\
&&  + \int_{D_+\times D_-} \alpha_u(x_2)\beta_v(y_1)\,F^{(01)(10)}_{u,v}(x_2,y_1)
    + \int_{D_-^2} \beta_u(y_1)\beta_v(y_2)\,F^{(01)(01)}_{u,v}(y_1,y_2) \\
&&  - \int_{D_+^2\times D_-} \alpha_v(x_2)\,\ell(x_1,y_1)\,F^{(11)(10)}_{u,v}((x_1,y_1),x_2)
    - \int_{D_+\times D_-^2} \beta_v(y_2)\,\ell(x_1,y_1)\,F^{(11)(01)}_{u,v}((x_1,y_1),y_2)\\
&&  - \int_{D_+^2\times D_-} \alpha_u(x_1)\,\ell(x_2,y_2)\,F^{(10)(11)}_{u,v}(x_1,(x_2,y_2))
    - \int_{D_+\times D_-^2} \beta_u(y_1)\,\ell(x_2,y_2)\,F^{(01)(11)}_{u,v}(y_1,(x_2,y_2))\\
&&  + \int_{D_+^2\times D_-^2} \ell(x_1,y_1)\,\ell(x_2,y_2)\,F^{(11)(11)}_{u,v}((x_1,y_1),(x_2,y_2)).
\end{eqnarray*}
Computing $\E[\eta_u - \xi_u]\cdot\E[\eta_v - \xi_v]$ in the same way, then using the definition of $\alpha_s$ and $\beta_s$ in (\ref{E:Abbreviate_Alpha_Beta}), we can rewrite the integrand in (\ref{E:LHS_BGprinciple}) as follows.
\begin{eqnarray}\label{E:Integrand_LHS}
&& \E[(\eta_u - \xi_u)(\eta_v - \xi_v)]-\E[\eta_u - \xi_u]\cdot\E[\eta_v - \xi_v]\\
&=& \int_{D_+^2\times D_-^2}\ell(x_1,y_1)\,\ell(x_2,y_2)\,\bigg\{\,
f^-_u(y_1)f^-_v(y_2)\,\Big[F^{(10)(10)}_{u,v}(x_1,x_2)- F^{(10)}_u(x_1)F^{(10)}_v(x_2)\Big] \notag\\
&& \qquad \qquad \qquad\qquad\qquad  +f^-_u(y_1)f^+_v(x_2)\,\Big[F^{(10)(01)}_{u,v}(x_1,y_2)- F^{(10)}_u(x_1)F^{(01)}_v(y_2)\Big] \notag\\
&& \qquad \qquad \qquad\qquad\qquad  +f^+_u(x_1)f^-_v(y_2)\,\Big[F^{(01)(10)}_{u,v}(y_1,x_2)- F^{(01)}_u(y_1)F^{(10)}_v(x_2)\Big] \notag\\
&& \qquad \qquad \qquad\qquad\qquad  +f^+_u(x_1)f^+_v(x_2)\,\Big[F^{(01)(01)}_{u,v}(y_1,y_2)- F^{(01)}_u(y_1)F^{(01)}_v(y_2)\Big] \notag\\
&& \qquad \qquad \qquad\qquad\qquad  -f^-_v(y_2)\,\Big[F^{(11)(10)}_{u,v}((x_1,y_1),x_2)- F^{(11)}_u(x_1,y_1)F^{(10)}_v(x_2)\Big] \notag\\
&& \qquad \qquad \qquad\qquad\qquad  -f^+_v(x_2)\,\Big[F^{(11)(01)}_{u,v}((x_1,y_1),y_2)- F^{(11)}_u(x_1,y_1)F^{(01)}_v(y_2)\Big] \notag\\
&& \qquad \qquad \qquad\qquad\qquad  -f^-_u(y_1)\,\Big[F^{(10)(11)}_{u,v}(x_1,(x_2,y_2))- F^{(10)}_u(x_1)F^{(11)}_v(x_2,y_2)\Big] \notag\\
&& \qquad \qquad \qquad\qquad\qquad  -f^+_u(x_1)\,\Big[F^{(01)(11)}_{u,v}(y_1,(x_2,y_2))- F^{(01)}_u(y_1)F^{(11)}_v(x_2,y_2)\Big] \notag\\
&& \qquad \qquad \qquad\qquad\qquad  +\,\Big[F^{(11)(11)}_{u,v}((x_1,y_1),(x_2,y_2))- F^{(11)}_u(x_1,y_1)F^{(11)}_v(x_2,y_2)\Big]\,\bigg\}. \notag
\end{eqnarray}
Note that each of the nine terms can be written in terms of
\begin{equation*}
E^{(n,m),(p,q)}_{u,r}(\vec{x},\vec{y},\;\vec{z},\vec{w}) :=
F^{(n,m),(p,q)}_{u,u+r}(\vec{x},\vec{y},\;\vec{z},\vec{w})- F^{(n,m)}_{u}(\vec{x},\vec{y})\cdot F^{(p,q)}_{u+r}(\vec{z},\vec{w})
\end{equation*}
defined in (\ref{Def:Enmpq}), where $r=v-u$. We split these nine terms into three groups $\Lambda_1(u,v) + \Lambda_2(u,v)+ \Lambda_3(u,v)$, where $\Lambda_1(u,v)$ consists of the first, third and fifth terms; $\Lambda_2(u,v)$ consists of the second, forth and sixth terms; and $\Lambda_3(u,v)$ consists of the last three terms. That is,
\begin{eqnarray}\label{E:Integrand_LHS_Gp1}
\Lambda_1(u,v) &:=& \int_{D_+^2\times D_-^2}\ell(x_1,y_1)\,\ell(x_2,y_2)\,\Big\{\,
f^-_u(y_1)f^-_v(y_2)\,E^{(10)(10)}_{u,r}(x_1,x_2) \notag\\
&&   +f^+_u(x_1)f^-_v(y_2)\,E^{(01)(10)}_{u,r}(y_1,x_2)    -f^-_v(y_2)\,E^{(11)(10)}_{u,r}((x_1,y_1),x_2) \,\Big\},
\end{eqnarray}
\begin{eqnarray}\label{E:Integrand_LHS_Gp2}
\Lambda_2(u,v) &:=& \int_{D_+^2\times D_-^2}\ell(x_1,y_1)\,\ell(x_2,y_2)\,\Big\{\,
f^-_u(y_1)f^-_v(y_2)\,E^{(10)(10)}_{u,r}(x_1,x_2) \notag\\
&&  +f^-_u(y_1)f^+_v(x_2)\,E^{(10)(01)}_{u,r}(x_1,y_2)   +f^+_u(x_1)f^+_v(x_2)\,E^{(01)(01)}_{u,r}(y_1,y_2) \notag\\
&&    -f^+_v(x_2)\,E^{(11)(01)}_{u,r}((x_1,y_1),y_2)\,\Big\}
\end{eqnarray}
and
\begin{eqnarray}\label{E:Integrand_LHS_Gp3}
 \Lambda_3(u,v) &:=& \int_{D_+^2\times D_-^2}\ell(x_1,y_1)\,\ell(x_2,y_2)\,\Big\{\,
-f^-_u(y_1)\,E^{(10)(11)}_{u,r}(x_1,(x_2,y_2)) \notag\\
&&    -f^+_u(x_1)\,E^{(01)(11)}_{u,r}(y_1,(x_2,y_2)) \notag\\
&&  +\,E^{(11)(11)}_{u,r}((x_1,y_1),(x_2,y_2))\,\Big\}.
\end{eqnarray}

\textbf{Step 3: Cancelations. }
To illustrate the `magical cancelations' mentioned at the beginning of the proof, we first provide details of these cancelations for $\Lambda_3$.

Note that we can bound $E_{u,r}$ in terms of $E_{u,0}$ via (\ref{E:Iterate_E_pq_2}). Consider the first among the three terms in $\Lambda_3$ with $E_{u,r}$ replaced $E_{u,0}$. We apply (\ref{E:FuuFu}) to write $F^{(10)(11)}_{u,u}$ in terms of $F^{(21)}_u$ plus a lower order term. This gives
\begin{eqnarray*}
&& -\int_{D_+^2\times D_-^2}\ell(x_1,y_1)\,\ell(x_2,y_2)\,f^-_u(y_1)\,E^{(10)(11)}_{u,0}(x_1,(x_2,y_2))\\
&=& -\int_{D_+^2\times D_-^2}\ell(x_1,y_1)\,\ell(x_2,y_2)\,f^-_u(y_1)\,\Big(F^{(21)}_{u}(x_1,x_2,y_2)-F^{(10)}_{u}(x_1)F^{(11)}_{u}(x_2,y_2)\Big)\\
&& -\,\frac{1}{N}\,\int_{D_+\times D_-^2}\ell(x,y_1)\,\ell(x,y_2)\,f^-_u(y_1)\,F^{(11)}_{u}(x,y_2).
\end{eqnarray*}

Similarly, when $r=0$, the second term and the third term in $\Lambda_2$ are, respectively,
\begin{eqnarray*}
&& -\int_{D_+^2\times D_-^2}\ell(x_1,y_1)\,\ell(x_2,y_2)\,f^+_u(x_1)\,E^{(01)(11)}_{u,0}(x_2,(y_1,y_2))\\
&=& -\int_{D_+^2\times D_-^2}\ell(x_1,y_1)\,\ell(x_2,y_2)\,f^+_u(x_1)\,\Big(F^{(12)}_{u}(x_2,y_1,y_2)-F^{(01)}_{u}(y_1)F^{(11)}_{u}(x_2,y_2)\Big)\\
&& -\,\frac{1}{N}\,\int_{D_+^2\times D_-}\ell(x_1,y)\,\ell(x_2,y)\,f^+_u(x_1)\,F^{(11)}_{u}(x_2,y)
\end{eqnarray*}
and
\begin{eqnarray*}
&& \int_{D_+^2\times D_-^2}\ell(x_1,y_1)\,\ell(x_2,y_2)\,E^{(11)(11)}_{u,0}((x_1,y_1),(x_2,y_2))\\
&=& \int_{D_+^2\times D_-^2}\ell(x_1,y_1)\,\ell(x_2,y_2)\,\Big(F^{(22)}_{u}(x_1,x_2,y_1,y_2)-F^{(11)}_{u}(x_1,y_1)F^{(11)}_{u}(x_2,y_2)\Big)\\
&& +\,\frac{1}{N}\,\int_{D_+\times D_-^2}\ell(x,y_1)\,\ell(x,y_2)\,F^{(12)}_{u}(x,y_1,y_2)\\
&& +\,\frac{1}{N}\,\int_{D_+^2\times D_-}\ell(x_1,y)\,\ell(x_2,y)\,F^{(21)}_{u}(x_1,x_2,y)\\
&& +\,\frac{1}{N^2}\,\int_{D_+\times D_-}\ell^2(x,y)\,F^{(11)}_{u}(x,y).
\end{eqnarray*}
Now we add up the three equations above. The sum of the lower order terms is, by Theorem \ref{T:Asymptotic_nm_correlation_t} or (\ref{E:Asymptotic_nm_correlation_t_1b}), of order $o(N)/N$ (i.e. a term which tends to zero even if we multiply it by $N$) uniformly for $u\in [0,t]$. On other hand, the sum of the leading terms is, by Theorem \ref{T:Asymptotic_nm_correlation_t} again, equal to
\begin{eqnarray}
&&\int_{D_+^2\times D_-^2}\ell(x_1,y_1)\,\ell(x_2,y_2)\,\Big\{\, -\,f^-_u(y_1)\,\Big(F^{(21)}_{u}(x_1,x_2,y_2)-F^{(10)}_{u}(x_1)F^{(11)}_{u}(x_2,y_2)\Big)
\notag\\
&& \qquad\qquad\qquad\qquad\qquad -\,f^+_u(x_1)\,\Big(F^{(12)}_{u}(x_2,y_1,y_2)-F^{(01)}_{u}(y_1)F^{(11)}_{u}(x_2,y_2)\Big)  \notag \\
&& \qquad\qquad\qquad\qquad\qquad +\,F^{(22)}_{u}(x_1,x_2,y_1,y_2)-F^{(11)}_{u}(x_1,y_1)F^{(11)}_{u}(x_2,y_2)\;\Big\} \notag \\
&=& \frac{1}{N}\int_{D_+^2\times D_-^2}\ell(x_1,y_1)\,\ell(x_2,y_2)\,\Big\{\, \notag\\
&& \qquad\quad -\,f^-_u(y_1)\,\Big(B^{(21)}_{u}(x_1,x_2,y_2)-A^{(10)}_{u}(x_1)B^{(11)}_{u}(x_2,y_2)-B^{(10)}_{u}(x_1)A^{(11)}_{u}(x_2,y_2)\Big) \notag\\
&& \qquad\quad -\,f^+_u(x_1)\,\Big(B^{(12)}_{u}(x_2,y_1,y_2)-A^{(01)}_{u}(y_1)B^{(11)}_{u}(x_2,y_2)-B^{(01)}_{u}(y_1)A^{(11)}_{u}(x_2,y_2)\Big) \notag\\
&& \qquad\quad +\,B^{(22)}_{u}(x_1,x_2,y_1,y_2)-A^{(11)}_{u}(x_1,y_1)B^{(11)}_{u}(x_2,y_2)-B^{(11)}_{u}(x_1,y_1)A^{(11)}_{u}(x_2,y_2)\;\Big\} \notag\\
&& \,+\, o(N)/N \notag\\
&=& \frac{-1}{N}\int_{D_+^2\times D_-^2}\ell(x_1,y_1)\,\ell(x_2,y_2)\,\Big\{\,   -\,f^-_u(y_1)\,\Big(\,G_u(x_1,y_2)f^+_u(x_2)+G^+_u(x_1,x_2)f^-_u(y_2)\,\Big) \notag\\
&& \qquad -\,f^+_u(x_1)\,\Big(\,G_u(x_2,y_1)f^+_u(x_1)+G^-_u(y_1,y_2)f^+_u(x_2)\,\Big) \notag\\
&& \qquad +\,\Big(\,G_u(x_1,y_2)f^+_u(x_2)f^-_u(y_1)+G_u(x_2,y_1)f^+_u(x_1)f^-_u(y_2)\,\Big) \notag\\
&& \qquad +\,\Big(\,G^+_u(x_1,x_2)f^-_u(y_1)f^-_u(y_2)+G^-_u(y_1,y_2)f^+_u(x_1)f^+_u(x_2)\,\Big)\;\Big\} \notag\\
&& \,+\, o(N)/N  \quad \text{by }(\ref{E:Asymptotic_nm_correlation_t_2}) \notag \\
&=& \,o(N)/N.  \label{E:MagicCancel_1}
\end{eqnarray}
The two $o(N)/N$ terms are the same and can be kept track of via the computation in the proof of Corollary \ref{cor:Asymptotic_nm_correlation_t_2}. Note that {\it  all}  terms involving $G,\,G^+,\,G^-$ cancel out in the last equality. The cancelation in (\ref{E:MagicCancel_1}), together with the cancelation for the lower order terms, are the `magical cancelations' mentioned at the beginning of the proof.

The same type of `magical cancelations' occur for each of $\Lambda_1$ and $\Lambda_2$ by the same reasons. In short, applying (\ref{E:FuuFu}) and (\ref{E:Asymptotic_nm_correlation_t_2}) to each of the six terms in $\Lambda_1+\Lambda_2$, we see that the sum of these six terms when $r=0$ is, up to an additive error of order $o(N)/N$ which is uniform for $u\in [0,t]$, equal to
\begin{eqnarray*}
&& \frac{1}{N}\int_{D_+^2\times D_-^2}\ell(x_1,y_1)\,\ell(x_2,y_2)\,\Big\{\, \notag\\
&& \qquad\quad \;\;\; f^-_u(y_1)f^-_v(y_2)\,\Big(B^{(20)}_{u}(x_1,x_2)-A^{(10)}_{u}(x_1)B^{(10)}_{u}(x_2)-B^{(10)}_{u}(x_1)A^{(10)}_{u}(x_2)\Big) \notag\\
&& \qquad\quad +\,f^-_u(y_1)f^+_v(x_2)\,\Big(B^{(11)}_{u}(x_1,y_2)-A^{(10)}_{u}(x_1)B^{(01)}_{u}(y_2)-B^{(10)}_{u}(x_1)A^{(01)}_{u}(y_2)\Big)\notag\\
&& \qquad\quad  +\,f^+_u(x_1)f^-_v(y_2)\,\Big(B^{(11)}_{u}(x_2,y_1)-A^{(01)}_{u}(y_1)B^{(10)}_{u}(x_2)-B^{(01)}_{u}(y_1)A^{(10)}_{u}(x_2)\Big) \notag\\
&& \qquad\quad  +\,f^+_u(x_1)f^+_v(x_2)\,\Big(B^{(02)}_{u}(y_1,y_2)-A^{(01)}_{u}(y_1)B^{(01)}_{u}(y_2)-B^{(01)}_{u}(y_1)A^{(01)}_{u}(y_2)\Big) \notag\\
&& \qquad\quad  -\,f^-_v(y_2)\,\Big(B^{(21)}_{u}(x_1,x_2,y_1)-A^{(11)}_{u}(x_1,y_1)B^{(10)}_{u}(x_2)-B^{(11)}_{u}(x_1,y_1)A^{(10)}_{u}(x_2)\Big) \notag\\
&& \qquad\quad  -\,f^+_v(x_2)\,\Big(B^{(12)}_{u}(x_1,y_1,y_2)-A^{(11)}_{u}(x_1,y_1)B^{(01)}_{u}(y_2)-B^{(11)}_{u}(x_1,y_1)A^{(01)}_{u}(y_2)\Big) \;\Big\}
\end{eqnarray*}

\begin{eqnarray}\label{E:MagicCancel_2}
&=& \frac{-1}{N}\int_{D_+^2\times D_-^2}\ell(x_1,y_1)\,\ell(x_2,y_2)\,\Big\{\, f^-_u(y_1)f^-_v(y_2)\,G^+_u(x_1,x_2)
 +\,f^-_u(y_1)f^+_v(x_2)\,G_u(x_1,y_2) \notag\\
&&   +\,f^+_u(x_1)f^-_v(y_2)\,G_u(x_2,y_1)    +\,f^+_u(x_1)f^+_v(x_2)\,G^-_u(y_1,y_2)  \notag\\
&&    -\,f^-_v(y_2)\,\Big(\,G_u(x_2,y_1)f^+_u(x_1)+G^+_u(x_1,x_2)f^-_u(y_1)\,\Big) \notag\\
 &&   -\,f^+_v(x_2)\,\Big(\,G_u(x_1,y_2)f^-_u(y_1)+G^-_u(y_1,y_2)f^+_u(x_1)\,\Big) \;\Big\}  \notag\\
&=& \,0.
\end{eqnarray}

Observe that on the RHS of $\Lambda_3(u,v)$ in (\ref{E:Integrand_LHS_Gp3}), if we view $u$ and $(x_1,y_1)$ as fixed variables, then
\begin{eqnarray*}
\Upsilon_r^{(p,q)}(x_2,y_2) &:=& -f^-_u(y_1)\,E^{(10)(pq)}_{u,r}(x_1,(x_2,y_2)) -f^+_u(x_1)\,E^{(01)(pq)}_{u,r}(y_1,(x_2,y_2)) \\ && +\,E^{(11)(pq)}_{u,r}((x_1,y_1),(x_2,y_2))
\end{eqnarray*}
satisfies
\begin{equation}\label{E:BBGKY_Upsilon_nmpq}
\Upsilon^{(p,q)}_{r}=P^{(p,q)}_r \Upsilon^{(p,q)}_{0} - \int_0^r P^{(p,q)}_{r-\theta}\left( V^+\Upsilon^{(p,q+1)}_{\theta} + V^-\Upsilon^{(p+1,q)}_{\theta} + \frac{Q}{N}E^{(p,q)}_{\theta} \right)\;d\theta
\end{equation}
since $E^{(p,q)}_r$ satisfies (\ref{E:BBGKY_Enmpq}). That is, $\{\Upsilon^{(p,q)}\}$ and $\{E^{(p,q)}\}$ solve the same hierarchy of equations, but the initial condition  $\Upsilon^{(p,q)}_{0}$ is of smaller order of magnitude $o(N)/N$, by the above cancelations. Following the same argument that we used for Lemma \ref{L:Tightness_Y_phi_2}, with $\Upsilon^{(p,q)}_{r}$ in place of $E^{(p,q)}_r$, while keeping track of these $o(N)$ terms, we obtain
\begin{equation}\label{E:FinalBoundLambda3}
N\,\int_0^t \int_u^t \Lambda_3(u,v)\,dv\,du \leq o(N)\,t^{3/2}
\end{equation}
whenever $0\leq t\leq\,T_0:=T\wedge (C_0\,C)^{-2}$ and $N>N_0$. By the same argument, (\ref{E:FinalBoundLambda3}) holds with $\Lambda_3$ replaced by either $\Lambda_1$ or $\Lambda_2$.

Recall that the integrand of (\ref{E:LHS_BGprinciple}) is $\Lambda_1+\Lambda_2+\Lambda_3$. The proof is complete.
\end{pf}

\subsection{Proof of main theorem}

With all the results developed in the previous sections, the proof of Theorem \ref{T:Convergence_AnnihilatingSystem} is ready to be presented in this section. Recall Steps 1-6 in the outline of proof at the end of Section 5.
 We will establish tightness of $\{\ZZ^N\}$ (which is Step 2) and then identify any subsequential limit through Steps 1, 3, 4, 5 and 6. Note that for Steps 1, 3, 4 and 5, we do not need to go into the analysis of correlation functions;
 the results for these steps are for arbitrary time interval rather than for a short time interval as in Steps 2 and 6.

The following is Step 2 in the outline of proof for Theorem \ref{T:Convergence_AnnihilatingSystem}. Note that we do not need any estimate about the evolution systems  $\mathbf{U}^N_{(t,s)}$ and $\mathbf{U}_{(t,s)}$ for this step. The key of the proof is Lemma \ref{L:Tightness_Y_phi_2}.

\begin{thm}\label{T:Tightness_Z}
(Step 2: Tightness)
Suppose Assumption \ref{A:ShrinkingRateCLT} holds and $\alpha>d\vee (d/2+2)$. For any $T>0$, there exists $C=C(D_+,D_-,T)>0$ such that $\{\ZZ^N\}$ is tight in $D([0,\,T_0],\,\mathbf{H}_{-\alpha})$, where $T_0:=T\wedge (\|u^+_0\|\vee \|u^-_0\|)^{-2}\,C$. Moreover, any subsequential limit has a continuous version.
\end{thm}

\begin{pf}
We first prove the following one dimensional tightness result: For any $\phi_{\pm}\in C(\bar{D}_{\pm})$ fixed (such as eigenfunctions), $\{(\Y^{N,+}(\phi_+),\,\Y^{N,-}(\phi_-))\}_N$ is tight in $D([0,T_0],\R^2)$. For this, it suffices to show $\{Z_N:= \Y^{N,+}(\phi_+)+\Y^{N,-}(\phi_-)\}_N$ is tight in $D([0,T],\R)$ for any fixed $\phi_{\pm}\in \H_{\pm}$ (cf. Problem 22 in Chapter 3 of \cite{EK86}). By Prohorov's Theorem. It suffices to show that
  \begin{itemize}
  \item[(i)]    for all $t\in[0,T_0]$ and $\epsilon_0>0$, there exists $K\in(0,\infty)$ s.t.
     $\varlimsup_{N\to\infty}\P\left(|Z_N(t)|> K\right)<\epsilon_0$;  and that
  \item[(ii)]    for all $\epsilon_0>0$, we have $$\lim_{\delta\to 0}\varlimsup_{N\to\infty}\P\left(\sup_{\substack{|t-s|<\delta\\0\leq s,t\leq T}}\big| Z_N(t)-Z_N(s) \big| >\epsilon_0\right)=0.$$
  \end{itemize}

By (\ref{E:Ch6_KeyStep_1}) and (\ref{E:QuadVarMN_CLT_N}), we have $Z_N:= \Y^{N,+}(\phi_+)+\Y^{N,-}(\phi_-)$ satisfies
\begin{eqnarray}\label{E:Tightness_ZNt-ZNs}
    Z_N(t)-Z_N(s) &=& \int_s^t \Y^{N,+}_r(\frac{1}{2}\Delta\phi_+)+\Y^{N,-}_r(\frac{1}{2}\Delta\phi_-)\,dr \notag\\
    && -\sqrt{N}\int_s^t\<\ell(\phi_++\phi_-),\,\X^{N,+}_r\otimes \X^{N,-}_r\>- \E[\<\ell(\phi_++\phi_-),\,\X^{N,+}_r\otimes \X^{N,-}_r\>]\,dr\notag \\
    && + M_{N}(t)-M_{N}(s)
\end{eqnarray}
for $0\leq s\leq t$, where $M_{N}(t)$ is a real valued $\F^{(\X^{N,+},\,\X^{N,-})}_t$-martingale with quadratic variation
\begin{equation}\label{E:QuadVarMN_Tightness_N}
    \int_0^t \<|\nabla \phi_+|^2 ,\,\X^{N,+}_s\> + \<|\nabla \phi_-|^2,\,\X^{N,-}_s\> +\<\ell(\phi_++\phi_-)^2,\,\X^{N,+}_s\otimes\X^{N,-}_s\>\;ds.
\end{equation}

(i) is implied by the fact that $\sup_{N\geq N_0(D)}\E[(Z_N(t))^2]<\infty$ for all $t\in [0,T_0]$ and $\alpha>d$. This fact can be proved as follows: By definition of the correlation functions, the covariance
\begin{eqnarray}\label{E:L2BoundCov+t+t}
&& \E\left[\Y^{N,+}_t(\phi)\,\Y^{N,+}_t(\psi)\right] \\
&=& N\left(\E[\<\phi,\X^{N,+}_t\>\<\psi,\X^{N,+}_t\>]-\E[\<\phi,\X^{N,+}_t\>]\E[\<\psi,\X^{N,+}_t\>]\right) \notag\\
&=& N\bigg(\int_{D_+^2}\phi(x_1)\psi(x_2)F^{(2,0)}_t(x_1,x_2)\,dx_1dx_2 + \frac{1}{N}\int_{D_+}(\phi\psi)(x)F^{(1,0)}_t(x)dx \notag\\
&& \quad - \int_{D_+}\phi(x_1)F^{(1,0)}_t(x_1)dx_1\cdot \int_{D_+}\psi(x_2)F^{(1,0)}_t(x_2)dx_2\bigg) \notag\\
&=& \int_{D_+} (\phi\psi)(x)f^+_t(x)dx + N\,\int_{D_+^2}\phi(x_1)\psi(x_2)\,\Big(F^{(2,0)}_t(x_1,x_2)-F^{(1,0)}_t(x_1)F^{(1,0)}_t(x_2)\Big)\,dx_1dx_2. \notag
\end{eqnarray}
By Theorem \ref{T:Asymptotic_nm_correlation_t} and Lemma \ref{L:BoundsForGGG}, the absolute value of the last quantity in (\ref{E:L2BoundCov+t+t}) is bounded above by
\begin{eqnarray*}
&& C_0\int_{D_+} |(\phi\psi)(x)|dx + \|\phi\|\,\|\psi\|\,\Big(\int_{D_+^2} |G^+(x_1,x_2)|\,dx_1dx_2 + \frac{(C_0\,C)^2\,t}{N\,\delta^{2d}}\Big) \notag \\
&\leq& C\,\|\phi\|\,\|\psi\|\,\Big(C_0+C_0^3\,t^{3/2}+ \frac{C_0^2\,t}{N\,\delta^{2d}}\Big) \notag\\
&\leq& C\,\|\phi\|\,\|\psi\|\,(C_0 \vee 1) \notag
\end{eqnarray*}
for all $0\leq t \leq T\wedge(C_0\,C)^{-2}$ and $N\geq N_0(D)$, where $C_0:=\|u^+_0\|\vee \|u^-_0\|$ and $C=C(D_+,D_-,T)$. In particular, $\E[(\Y^{N,+}_t(\phi^+_k))^2] \leq C\,\|\phi^+_k\|^2$. Similarly, we have $\E[(\Y^{N,-}_t(\phi^-_k))^2] \leq C\,\|\phi^-_k\|^2$. Therefore, when $\alpha>d$, we have $\E[(Z_N(t))^2]<\infty$ for all $t\in [0,T_0]$ and $N\geq N_0(D)$ (as in the proof of Lemma \ref{L:StateSpace_AnnihilatingSystem}). Hence (i) is satisfied.

It remains to show that (ii) holds with $Z_N(t)-Z_N(s)$ replaced by each of the three terms on the RHS of (\ref{E:Tightness_ZNt-ZNs}). For the first term, (2) holds by Chebyshev's inequality, Holder's inequality and (\ref{E:L2BoundCov+t+t}). For the second term, (2) holds by Lemma \ref{L:Tightness_Y_phi_2}. For the third term, namely $M_{N}(t)-M_{N}(s)$, we have (ii) holds upon applying Chebyshev's inequality, Doob's maximal inequality and the explicit expression for the quadratic variation (\ref{E:QuadVarMN_Tightness_N}). Hence we have one dimensional tightness for fixed $\phi_{\pm}\in C(\bar{D}_{\pm})$.

Following the same proof of \cite[Theorem 4.7]{zqCwtF14c}, we complete the proof by using the definition (\ref{E:norm_MinusAlpha_AnnihilatingSystem}) of the metric of $\mathbf{H}_{-\alpha}$ and the condition on $\alpha$.
\end{pf}

We identify any subsequential limit of $\{\ZZ^N\}$ for the rest of this section. Steps 1, 3, 4 and 5 follow from the method developed in \cite{zqCwtF14c}, via the estimates for $\mathbf{U}^N_{(t,s)}$ and $\mathbf{U}_{(t,s)}$ that we developed. We now present the precise statements that we obtain.

Using  Lemma \ref{L:KeyMtgAnnihilatingDiffusionModel}, we can follow the proof of \cite[Theorem 4.3]{zqCwtF14c} to obtain the following.

\begin{thm}\label{T:3.3_CLT} (Step 1)
Suppose $\alpha >d\vee (d/2+1)$. For all $N$ large enough, there exists a $c\grave{a}dl\grave{a}g$ square integrable $\mathbf{H}_{-\alpha}$-valued $\F^{N}_t$-martingale $M^N=(M^N_t)_{t\geq 0}$ such that
    \begin{equation}\label{E:EvolutionSol_Annihilation}
        \ZZ^{N}_t= \mathbf{U}^N_{(t,0)}\ZZ^{N}_0+ \int_0^t \mathbf{U}^N_{(t,s)}\,dM^N_s + \int_0^t \mathbf{U}^N_{(t,s)}(\mathbf{B}^N_s\ZZ^N_s- K^N_s)\,ds  \qquad\text{for }t\geq 0,\quad \P-a.s,
    \end{equation}
where $\mathbf{U}^N_{(t,s)}$ is defined in Definition \ref{Def:U^N_ts} and $\F^N_t$ is the natural filtration of the annihilating diffusion process. Moreover, $M^N$ has bounded jumps and predictable quadratic variation given by (\ref{E:QuadVarMN_CLT_N}).
\end{thm}

As a remark, equation (\ref{E:EvolutionSol_Annihilation}) is equivalent to (\ref{E:Ch6_KeyStep_1}) by variation of constant (see Section 2.1.2 of \cite{GT95}). For Step 4, it can be checked that we have the following, as in \cite[Theorem 4.8]{zqCwtF14c}.

 \begin{lem}\label{L:ConvergenceOfUNZN0} (Step 4)
For $\alpha>d+2$ and $T>0$, we have
\begin{equation*}\label{E:ConvergenceOfUNZN0}
\mathbf{U}^N_{(t,0)}\ZZ^{N}_0 \toL \mathbf{U}_{(t,0)}\ZZ_0 \quad \text{in  }D([0,T],\mathbf{H}_{-\alpha}).
\end{equation*}
Moreover, $\mathbf{U}_{(t,0)}\ZZ_0$ has a version in $C^{\gamma}([0,T],\,\mathbf{H}_{-\alpha})$ for any $\gamma\in(0,1/2)$.
\end{lem}

By Theorem \ref{T:Conjecture_delta_N_CLT} and Lemma \ref{L:MinkowskiContent_I_k_CLT}, we can check that the quadratic variation of $M^N$ converges in probability to the deterministic quantity (\ref{E:QuadVarMN_CLT}). Hence, by a standard functional central limit theorem for semi-martingales (see, e.g., \cite{LS80}), we have for any $\phi_{\pm}\in Dom^{Feller}(\A^{\pm})$ fixed, $\{M^N(\phi_+,\phi_-)\}$ converges in distribution in $D([0,T],\R)$ to a continuous Gaussian martingale with independent increments and covariance functional (\ref{E:QuadVarMN_CLT}). In fact, following the proof of \cite[Theorem 4.6]{zqCwtF14c}, we obtain Step 3.

\begin{thm}\label{T:ConvergenceOfM^N_Annihilation} (Step 3)
When $\alpha >d \,\vee (d/2+1)$, the square-integrable martingale $\{M^N\}$ in Theorem \ref{T:3.3_CLT} converges to $M$ in distribution in $D([0,T],\mathbf{H}_{-\alpha})$ for any $T>0$, where $M$ is the (unique in distribution) continuous square-integrable $\mathbf{H}_{-\alpha}$-valued Gaussian martingale with independent increments and covariance functional characterized by (\ref{E:QuadVarMN_CLT}).
\end{thm}

With Lemma \ref{L:MinkowskiContent_I_k_CLT}, we can check, as in \cite{zqCwtF14c}, that the expression $\int_0^{t} \mathbf{U}_{(t,s)}dM_{s}$ is well-defined. That is $\mathbf{U}_{(t,s)}$ (for $s\in[0,t]$) lies within the class of integrands with respect to  $M$. Furthermore, following the same proof for \cite[Theorem 4.9]{zqCwtF14c}, we obtain the following.

 \begin{thm}\label{T:ConvergenceStochInt_CLT} (Step 5)
For $\alpha>d+2$ and $T>0$, we have
\begin{equation}
\int_0^{t} \mathbf{U}^N_{(t,s)}dM^N_{s} \toL \int_0^{t} \mathbf{U}_{(t,s)}dM_{s} \quad \text{in  }D([0,\,T],\mathbf{H}_{-\alpha})
\end{equation}
Moreover,  $\int_0^{t} \mathbf{U}_{(t,s)}dM_{s}$ has a version in $C^{\gamma}([0,T],\,\mathbf{H}_{-\alpha})$ for any $\gamma\in(0,1/2)$.
\end{thm}

Finally, we complete Step 6 towards the proof of Theorem \ref{T:Convergence_AnnihilatingSystem}. The key is Lemma \ref{L:BGprinciple}.

 \begin{thm}\label{T:BGprinciple}(Step 6: Boltzman-Gibbs principle)
Suppose $\alpha>d+2$ and Assumption \ref{A:ShrinkingRateCLT} holds. For any $T>0$, there exists $C=C(D_+,D_-,T)>0$ such that
\begin{equation}\label{E:BGprinciple}
\int_0^t \mathbf{U}^N_{(t,s)}(\mathbf{B}^N_s\ZZ^N_s- K^N_s)\,ds \toL 0 \quad\text{in }D([0,\,T_0],\,\mathbf{H}_{-\alpha}),
\end{equation}
where $T_0:=T\wedge (C_0\,C)^{-2}$, the operator $\mathbf{U}^N_{(t,s)}$ is defined in (\ref{Def:U^N_ts}),
\begin{eqnarray*}
\mathbf{B}^N_s\mu(\phi_+,\phi_-)&:=& \mu \big(\,\<\ell(\phi_+ + \phi_-),\,f^-_s\>_{-},\; \<\ell(\phi_+ + \phi_-),\,f^+_s\>_{+} \,\big) \quad\text{and}\\
K^N_s(\phi_+,\phi_-) &:=&  \sqrt{N}\Big(\,\<\ell_{\delta_N}(\phi_++\phi_-),\,\X^{N,+}_s\otimes \X^{N,-}_s\>- \E[\<\ell(\phi_++\phi_-),\,\X^{N,+}_s\otimes \X^{N,-}_s\>]\,\Big).
\end{eqnarray*}
\end{thm}

\begin{pf}
Observe that $\alpha>d$ (we will need $\alpha>d+2$ later in the proof) guarantees, base on Weyl's law (\ref{E:WeylLaw_CLT}) and (\ref{E:EigenfcnUpperBound_CLT}), that
$$\sum_{k\geq 1} \left( \frac{\|\phi^{+}_{k}\|^2}{(1+\lambda^{+}_{k})^{\alpha}} +
\frac{\|\phi^{-}_{k}\|^2}{(1+\lambda^{-}_{k})^{\alpha}} \right) <\infty.$$
Using the definition of the norm $|\,\cdot\,|_{-\alpha}$ is defined in (\ref{E:norm_MinusAlpha_AnnihilatingSystem}), the uniform bound (\ref{E:ContractionQNQ_AnnihilatingSystem}) and Lemma \ref{L:BGprinciple}, we have the following: For any $T>0$, there exists a constant $C=C(D_+,D_-,T)>0$,  an integer $N_0=N_0(D_+,D_-)$ and positive constants $\{C_N\}$ satisfying $\lim_{N\to\infty}C_N=0$ such that
\begin{equation}\label{E:BGprinciple_2}
   \E\Big[\,\Big|\int_0^t \mathbf{U}^N_{(t,s)}(\mathbf{B}^N_s\ZZ^N_s- K^N_s)\,ds\Big|^2_{-\alpha} \,\Big]\leq C_N\,t^{3/2}
\end{equation}
whenever $0\leq t\leq \,T_0:=T\wedge (C_0\,C)^{-2}$ and $N>N_0$.  In particular, we have, for $\alpha>d$,
\begin{equation}\label{E:BGprinciple_3}
\lim_{N\to\infty}\sup_{t\in[0,\,T_0]}\E\Big[\,\Big|\int_0^t \mathbf{U}^N_{(t,s)}(\mathbf{B}^N_s\ZZ^N_s- K^N_s)\,ds\Big|^2_{-\alpha}\,\Big] = 0.
\end{equation}
On other hand, the process $\mathbf{e}_N(t):=\int_0^t \mathbf{U}^N_{(t,s)}(\mathbf{B}^N_s\ZZ^N_s- K^N_s)\,ds$ is tight in $D([0,\,T_0],\,\mathbf{H}_{-\alpha})$. This can be verified by the same argument that we used for $\ZZ$ in the proof of Theorem \ref{T:Tightness_Z}. Precisely, by (\ref{E:EvolutionSol_Annihilation}), we have almost surely,
    \begin{equation*}
       \mathbf{e}_N(t)= \ZZ^{N}_t- \mathbf{U}^N_{(t,0)}\ZZ^{N}_0- \int_0^t \mathbf{U}^N_{(t,s)}\,dM^N_s  \qquad\text{for }t\geq 0.
    \end{equation*}
Each of the three terms on the RHS is $C$-tight (i.e. has only continuous limits) in $D([0,\,T_0],\,\mathbf{H}_{-\alpha})$ by Theorem \ref{T:Tightness_Z}, Lemma \ref{L:ConvergenceOfUNZN0} and Theorem \ref{T:ConvergenceStochInt_CLT} respectively, provided that $\alpha>d+2$. Hence $\mathbf{e}_N$ is tight in $D([0,\,T_0],\,\mathbf{H}_{-\alpha})$. Now Theorem \ref{T:BGprinciple} follows from (\ref{E:BGprinciple_3}).
\end{pf}

\medskip

This completes the proof of Theorem \ref{T:Convergence_AnnihilatingSystem}.

\vspace{5mm}
    \textbf{Zhen-Qing Chen}

    Department of Mathematics, University of Washington, Seattle, WA 98195, USA

    Email: zqchen@uw.edu
\vspace{2mm}

    \textbf{Wai-Tong (Louis) Fan}

Department of Mathematics, University of Wisconsin,  Madison, WI 53706, USA

Email: louisfan@math.wisc.edu

\end{document}